\documentclass{svmult}

\usepackage{makeidx}
\usepackage{graphicx}

\usepackage{multicol}
\usepackage[centertags]{amsmath}
\usepackage{amsfonts}
\usepackage{euscript}
\usepackage{citesort}
\allowdisplaybreaks[4]

\makeindex

\begin{document}

\title*{Fictitious Fluid Approach and Anomalous Blow-up of the Dissipation Rate in a 2D Model of Concentrated Suspensions}
\titlerunning{Anomalous blow-up of the Viscous Dissipation Rate}

\author{Leonid Berlyand\inst{1}\and
Yuliya Gorb \inst{2} \and Alexei Novikov \inst{2} }

\institute{Department of Mathematics
              \& Materials Research Institute,
              Pennsylvania State University, McAllister Bld.,
              University Park, PA 16802, USA, \texttt{berlyand@math. psu.edu}
 \and Department of Mathematics,
              Pennsylvania State University,
              University Park, PA 16802, USA,
              \texttt{gorb@math.psu.edu}, \texttt{anovikov@math.psu.edu}}

\renewcommand{\thesection}{\arabic{section}}
            \renewcommand{\thesubsection}{\arabic{subsection}}
                    \renewcommand{\theequation}{\thesection.\arabic{equation}}
                      \renewcommand{\thetheorem}{\thesection.\arabic{theorem}}
                     \renewcommand{\thefigure}{\thesection.\arabic{figure}}
                     \renewcommand{\thetable}{\thesection.\arabic{table}}
                     \renewcommand{\thelemma}{\thesection.\arabic{lemma}}
                     \renewcommand{\theproposition}{\thesection.\arabic{proposition}}
                     \renewcommand{\theremark}{\arabic{remark}}
                    \renewcommand{\theremark}{\thesection.\arabic{remark}}
                     \renewcommand{\theexample}{\arabic{example}}
                    \renewcommand{\theexample}{\thesection.\arabic{example}}

\newcommand{\Section}[1]{\section{#1} \setcounter{equation}{0} \setcounter{figure}{0} \setcounter{theorem}{0} \setcounter{proposition}{0}
\setcounter{lemma}{0} \setcounter{remark}{0}
\setcounter{example}{0}}

\maketitle

\begin{abstract}
We present a two-dimensional (2D) mathematical model of a highly
concentrated suspension or a thin film of the rigid inclusions in
an incompressible Newtonian fluid. Our objectives are two-fold:
$(i)$ to obtain all singular terms in the asymptotics of the
overall viscous dissipation rate as the interparticle distance
parameter $\delta$ tends to zero, $(ii)$ to obtain a qualitative
description of a microflow between neighboring inclusions in the
suspension.

Due to reduced analytical and computational complexity, 2D models
are often used for a description of 3D suspensions. Our analysis
provides the limits of validity of 2D models for 3D problems and
highlights novel features of 2D physical problems (e.g. thin
films). It also shows that the Poiseuille type microflow
contributes to a singular behavior of the dissipation rate. We
present examples in which this flow results in anomalous rate of
blow up of the dissipation rate in 2D. We show that this anomalous
blow up has no analog in 3D.

While previously developed techniques allowed to derive and
justify the leading singular term only for special symmetric
boundary conditions, a \textit{fictitious fluid approach},
developed in this paper, captures \textit{all} singular terms in
the asymptotics of the dissipation rate for generic boundary
conditions. This approach seems to be an appropriate tool for
rigorous analysis of 3D models of suspensions as well as various
other models of highly packed composites.
\end{abstract}


\Section{Introduction} \label{S:introduction}

Many classical and novel engineering processes involving
multiphase flow require to capture the overall behavior of
suspensions. The problem of the behavior of  suspensions is
important  in geophysics (mud-flow and debris flow rheology),
pharmacology (drugs design), ceramics processing among others.
Wide range of experimental (see, e.g.,
\cite{ATGAFMS,ASNS1,ASNS2,CSI,HMGT,KLH,LA2,LL1,LL2,PABGA}), and
numerical (see, e.g., \cite{JGP,NB,bB,sb,SIMAG,ladd}) results are
available.

Two-dimensional (2D) models of three-dimensional (3D) suspensions
are often used in numerical simulations (e.g. \cite{mau,hvmf,dwz}), because they require less computational
effort. In addition, 2D suspensions can describe biological thin films \cite{wu_l} and certain types of 3D concentrated suspensions
of uniaxial thin rods (see e.g. \cite{larson}).

We consider a 2D mathematical model of a non-colloidal (Browning
motion can be neglected) concentrated suspension of neutrally
buoyant rigid particles (inclusions) in a Newtonian fluid. The
suspension occupies a 2D domain $\Omega$, and rigid inclusions are
modeled by disks of equal radii. The main objective is to
characterize in a rigorous mathematical framework the dependence
of the effective rheological properties (e.g. effective
viscosity, effective permeability, effective viscous dissipation rate) of such suspension on the %the irregular
geometry of inclusions array and applied boundary conditions on
$\partial\Omega$. We focus on highly packed suspensions when the
concentration of inclusions is close to maximal, which means that
the distance between neighboring inclusions (\textit{interparticle
distance}) is much smaller than their sizes. We consider an
irregular (non-periodic) array of disks and our analysis takes
into account the variable distances between adjacent inclusions.

The main features of the problem under consideration are the high
concentration of the inclusions and the irregular geometry of
their spatial distribution. The key quantity of interest in
describing the effective rheological properties of suspensions is
the rate of viscous dissipation of energy $\displaystyle
\widehat{W}=\displaystyle \widehat{W}(\boldsymbol{u})$, where
$\boldsymbol{u}$ is the velocity of the incompressible fluid (see
the precise definition in Section~\ref{S:formulation}).

Initially our interest was motivated by the problem of dewatering
process which led to the sedimentation due to gravity problem
\cite{burger} arising in the disposal of solid waste in industrial
and municipal pollution management. In this problem particles are
close to touching. In particular, shear external boundary
condition was brought up by experimental and numerical studies
\cite{glow,GO,GOE} and it was addressed mathematically by analysis of
2D model in \cite{bp}. This led us to the question of limits of
validity of 2D models for description of 3D suspension.

Note that the problem of sedimentation is quite difficult in
mathematical content and far from being settled. Even in the case
of the dilute limit this problem is not completely understood
mathematically and we mention here the number of recent
mathematical studies in this direction \cite{otto,bw}.

For highly concentrated suspensions of rigid inclusions,
$\widehat{W}$ exhibits a singular behavior (see e.g.
\cite{fa,graham,sb,bbp,bp,kn}) and its understanding is a
fundamental issue. A formal asymptotic analysis of the singular
behavior of the viscous dissipation rate in a thin gap between a
single pair of two closely spaced spherical inclusions in a
Newtonian fluid was performed in \cite{fa}. In this work only
translational motions of inclusions but not rotational were
considered and the asymptotics of the form $C\delta^{-1}+O(\ln
1/\delta)$, where $\delta$ is the distance between two spheres,
was obtained. Based on the analysis of a single pair of spheres
the authors of \cite{fa} suggested that the asymptotics of the
effective viscosity of the periodic array of inclusions is of the
same form, that is, the main singular behavior is of
$O(\delta^{-1})$. Next, a periodic array of inclusions in a
Newtonian fluid was considered in \cite{kn}. Under the assumption
that all inclusions follow the shear motion of the fluid (formula
(5) in \cite{kn}) it was shown that $\displaystyle \widehat{W}=
O(\delta^{-1})$. This assumption is analogous to the well-known
Cauchy-Born hypothesis in solid state physics, which is known to
be {\it not} always true \cite{FT}. Indeed, in the case of
suspensions it was shown in \cite{bp} that for shear external
boundary conditions the inclusions may not follow the shear
motion. Moreover, it was shown in  numerical studies of \cite{sb}
that the asymptotics $O(\delta^{-1})$ may or may not hold for
suspensions of a large number of inclusions with generic boundary
conditions. There it was observed numerically that while in some
cases the asymptotics of the effective viscosity is of order
$1/\delta$, in other cases it is of order $\ln 1/\delta$. Also the
problem of the exact analytical form of the singular behavior for
generic suspensions was posed in \cite{sb} (p. 140) which
motivated subsequent studies of \cite{bbp,bp} and the present
paper, where a network approximation approach was developed for
this type of problems.

In \cite{bbp} a formal asymptotic analysis of the effective
viscosity in 3D for a disordered array of inclusions was
performed. In a view of discrepancies between predictions of the
formal asymptotic analysis  \cite{fa} and numerics
 \cite{sb} such formal asymptotics requires a mathematical
justification. In \cite{bbp} for special (extensional) boundary
conditions the leading term of the asymptotics of the effective
viscosity as $\delta\rightarrow 0$, where $\delta$ is the
characteristic spacing between neighboring inclusions, was
justified in a 2D model. In subsequent work \cite{bp} this leading
term, that exhibits a so-called \textit{strong blow up} of order
$\delta^{-3/2}$, was analyzed. It turned out that in many
important cases, e.g. when shear boundary conditions are applied,
this term degenerates, so the next term of order $\delta^{-1/2}$,
that exhibits a so-called \textit{weak blow up}, becomes the
leading term of the asymptotics in many physical situations.

However, the techniques of \cite{bbp} are only capable to capture
the strong blow up but not the other singular terms in the
asymptotic expansion of the effective viscosity. By contrast, a
\textit{fictitious fluid approach} introduced in this paper allows
to derive the complete asymptotic expansion of the overall viscous
dissipation rate in which all singular terms are captured and justified. In
particular, we ruled out  singular terms other than presented in
Theorem \ref{T:main-thm3}. Previously \cite{fa,graham,sb,bbp,bp,kn,bgmp} inclusions in a
dense suspension were characterized by the sets of their
translational and rotational velocities, so-called discrete
variables. Our analysis shows that in order to obtain the complete
asymptotics of singular behavior it is necessary use an additional
set of discrete variables, permeation constants. To the best of
our knowledge these discrete variables have not been used in
previous studies of dense suspensions. We now explain the physical
consequences of this asymptotics.

The key feature of rheology of concentrated suspensions is that
the dominant contribution to the overall or effective viscous
dissipation rate (and therefore, to the effective viscosity) comes
from thin gaps (lubrication regions) between closely spaced
neighboring inclusions \cite{leal}. The mathematical techniques
introduced in the above mentioned works \cite{fa,kn,bbp} utilized
this observation. More specifically, they took into account
certain types of relative movements of inclusions which resulted
in the corresponding microflows in the gaps. The formal
asymptotics in \cite{fa} was based on analysis of the squeeze motion,
when two inclusions move toward each other along the line joining
their centers (see Fig. \ref{F:bc1}$c$) but did not provide the
detailed analysis of other motions, which was sufficient for
certain type of boundary conditions (e.g. extensional boundary
conditions) but not sufficient for others. A justification of the
formal asymptotics was not considered in \cite{fa} (see also
\cite{graham} where similar results were obtained).

In \cite{bbp} four types of relative motions of neighboring
inclusions were considered which result in a singular behavior of
dissipation rate: the squeeze (Fig. \ref{F:bc1}$c$), the shear
(Fig. \ref{F:bc1}$b$) and two rotations (Fig. \ref{F:bc2}). While
it was sufficient for the analysis of the leading term of the
overall viscous dissipation rate (in a suspension of free particles) which was the goal of \cite{bbp},
in the present paper we observed that this analysis does not
provide a complete picture of microflows. Indeed, the Poiseuille
flow in 2D also contributes into the singular behavior. Examples in the present
paper suggest that when an external field is applied to inclusions
in a suspension, this Poiseuille flow may give rise to an
anomalous strong rate of blow up (called a \textit{superstrong
blow up}, of order $\delta^{-5/2}$) of the viscous dissipation
rate, whereas for suspensions of free inclusions there is at most strong blow up (of order $\delta^{-3/2}$).

As explained above the previous studies \cite{fa,kn,bbp,graham}
provided only a partial analysis of the single behavior of the
viscous dissipation rate and therefore did not provide the
complete physical picture of microflows. In this work we obtain
the complete asymptotic description of the singular behavior of
viscous dissipation rate which led us to the complete description
of microflows (Fig. \ref{F:bc1}-\ref{F:bc3}).

The techniques of \cite{bbp,bp} and the present paper are based on
the discrete network approximation. Discrete networks have been
used as analogues of the continuum problems in various areas of
physics and engineering for a long time (see, for example,
\cite{bk} and references therein, and see also the recent review
\cite{newman} for various applications of networks in social and
biological studies). However, the fundamental issue of
relationship between a continuum problem and the corresponding
discrete network was not addressed until recently.

The pioneering study of electro-magnetic properties of high
contrast materials arising in imaging using a rigorous network
approximation was done in\cite{borcea,bor_pap,bobp}. There the local
resistivity of a periodic medium with continuously distributed
properties was given by the Kozlov's function (\cite{koz}):
$\rho(\boldsymbol{x})=\mbox{Const}\,
e^{-S(\boldsymbol{x})/\varepsilon^2}$, where $S(\boldsymbol{x})$
is a periodic smooth phase function, $\varepsilon$ is an aspect
ratio of the material properties, a small parameter of the problem
(as $\varepsilon\rightarrow 0$, $\rho(\boldsymbol{x})$ describes a
resistivity of a high contrast
composite). The effective resistivity
obtained and justified there was given in terms of the principal
curvatures $\kappa^+$, $\kappa^-$ of the function
$S(\boldsymbol{x})$ at the saddle point of $S(\boldsymbol{x})$.

The discrete network approximation for a medium with piecewise
constant characteristics, which correspond to particle-filled
(particulate) composites, was developed in \cite{bk}. In this
case, the principle curvatures and saddle points of the phase
function are not determined. In \cite{bk} an infinite contrast
material ($\varepsilon=0$) with irregularly or non-periodically
distributed circular fibers was considered and asymptotics of the
effective conductivity of such a composite as $\delta \rightarrow
0$ was derived and justified. In \cite{bk} the Keller's asymptotic
solution \cite{keller} for two closely spaced disks was used as a
building block for construction of the network for a large number
of inclusions (disks). The network model is a linear system on a
graph with appropriate boundary conditions. The sites or vertices
of the graph correspond to the centers of inclusions
$\{\boldsymbol{x}_i, \,\, i=1,\ldots,N\}$ and the values of the
local fluxes between adjacent inclusions $\boldsymbol{x}_i$ and
$\boldsymbol{x}_j$ described by the Keller's formula:
$\pi\sqrt{\frac{R}{\delta_{ij}}}$, where $R$ is the radius of
disks, $\delta_{ij}$ is the distance between the closely spaced
neighbors, are assigned to the edges of the graph. The notion of
the neighbors is introduced via the Voronoi tessellation. The
asymptotic closeness of the effective conductivity of continuum
and discrete problems was justified by using the direct and dual
variational bounds matching to the leading order which was
obtained to be of order $\delta^{-1/2}$ as $\delta\rightarrow 0$,
where $\delta=\max \delta_{ij}$.

In \cite{bn} the notion of the characteristic interparticle
distance was generalized for a broad class of 2D geometrical
patterns (highly non-uniform arrays of inclusions). For this class
of arrays an explicit error estimate for the discrete network
approximation was obtained. A scalar conductivity problem
analogous to \cite{bk,bn} in 3D was considered in \cite{bgn} where
both asymptotics of the effective conductivity and its error
estimated were obtained. While the formulation of the problem is
analogous to the 2D one considered in \cite{bk,bn}, the
connectivity patterns of conducting inclusions are much more
complex in 3D. That is why the justification of this approximation
and derivation of its error estimate required introduction of new
techniques (e.g., the central projection partition). By contrast,
our work shows that 2D vectorial problems may exhibit features
which have no analogs in similar 3D  problems.

Next the method of the discrete network approximation was extended to
a vectorial problem \cite{bbp}, where formal asymptotics for 3D
suspensions was derived. The approach developed in \cite{bk,bn}
for the justification of the formal asymptotics in the scalar
problem could not be readily applied to the vectorial problems,
which exhibit new features and require new technical tools. In
particular, the key point of discrete network approach is a
construction of trial functions in the variational formulation of
the problem for determining the effective properties of
composites. In vectorial problems (unlike in scalar ones) such trial functions (for the direct variational principle) must satisfy the divergence free
condition. Also the dual variational formulation involves integral
constraints - balances of  linear and angular momenta. The balance
equations on an inclusion involve all its neighbors that have
their own neighbors etc. This accounts for a long-range
hydrodynamic interactions and shows that the pairwise (local)
analysis (such as \cite{fa,graham}) is not sufficient and a global
analysis is necessary for the whole array of inclusions. That is why a
generalization of scalar techniques from \cite{bk,bn} was done for
the leading term only in the asymptotics of the effective
viscosity in the 2D model \cite{bbp}. It seems that a
straightforward generalization of the approach in \cite{bbp}
becomes increasingly difficult and probably impossible due to the
global nature of above mentioned constraints.

The objectives of this work are two-fold: ($i$) develop a new
approach, called \textit{a fictitious fluid approach}, that allows
to effectively deal with global constraints (balance equations and the divergence free condition), ($ii$) apply the
fictitious fluid approach to derive and justify \textit{all}
singular terms in the asymptotics of the overall viscous
dissipation rate. Such an asymptotic formula results in a
\textit{complete} description of microflows in suspensions, as
oppose to the partial description of the microflow in
\cite{fa,graham,kn,bbp}.

The techniques in \cite{bk,bn,bgn} for scalar problems in both 2D
and 3D were based on a direct (one-step) discretization of the
original continuum problem. In vectorial problems \cite{bbp} the
direct discretization allowed to obtain only the leading term for
the extensional viscosity (strong blow up) but not for the shear
viscosity (weak blow up) due to the global constraints. In the
present work we use a two-step discretization procedure based on
the fictitious fluid approach to derive and justify both strong
and weak blow up terms.

We briefly describe now the idea behind the fictitious fluid
approach. As mentioned before, it consists of two steps. In step
$1$ we introduce a fictitious fluid domain which comprises thin
gaps between neighboring inclusions. The dissipation rate
restricted to this domain is denoted by $\widehat{W}_{\boldsymbol{\Pi}}$. We show
that for generic Dirichlet boundary conditions it describes the
singular behavior in the following sense:
\[
\widehat{W}=\widehat{W}_{\boldsymbol{\Pi}}+h.o.t. \quad \mbox{as} \quad
\delta\rightarrow 0.
\]

In step $2$ we perform a discretization of $\widehat{W}_{\boldsymbol{\Pi}}$, that
is, the continuum problem for $\widehat{W}_{\boldsymbol{\Pi}}$ is reduced to an
algebraic problem on a graph, called a \textit{network problem}.
By using the fictitious fluid approach in this step most of the difficulties are eliminated.
The network problem is a minimization of a quadratic form whose coefficients depend on $R$, $\mu$ and boundary data.
The quadratic form on the minimizing set of discrete variables is called the \textit{discrete viscous dissipation rate} and denoted by $\mathcal{I}$.
Unknowns of this problem are vectors $\mathbb{U}=\{\boldsymbol{U}^i\}_{i=1}^{N}$,
$\boldsymbol{\omega}=\{\omega^i\}_{i=1}^{N}$, the translational
and angular velocities of inclusions, respectively, and a
collection of numbers $\boldsymbol{\beta}=\{\beta_{ij}\}$ characterizing the
Poiseuille microflow between a pair of inclusions. The discrete
dissipation rate $\mathcal{I}$ is given by:
\begin{equation}  \label{E:intro}
\mathcal{I}=\mathcal{I}_{1}(\boldsymbol{\beta})\delta^{-5/2}+
\mathcal{I}_{2}(\mathbb{U},\boldsymbol{\omega},\boldsymbol{\beta})\delta^{-3/2}+
\mathcal{I}_{3}(\mathbb{U},\boldsymbol{\omega},\boldsymbol{\beta})\delta^{-1/2},
\mbox{ as } \delta\rightarrow 0
\end{equation}
where $\mathcal{I}_{k}$, $k=1,2,3$, are explicitly computable
quadratic polynomials of $\mathbb{U}$, $\boldsymbol{\omega}$,
$\boldsymbol{\beta}$.

As a result we obtain the following asymptotic formula for the
generic Dirichlet boundary conditions:
\[
\frac{|\widehat{W}-\mathcal{I}|}{\mathcal{I}}=O(\delta^{1/2})
\quad \mbox{as} \quad \delta\rightarrow 0.
\]
In fact, we prove the following result about the error term:
\[
|\widehat{W}-\mathcal{I}|\leq
\mu\left(\sum_{i,j}\mathbb{E}_1(\beta_{ij})+\mathbb{E}_2(\boldsymbol{U}^i-\boldsymbol{U}^j)+\mathbb{E}_3(\omega^i)\right)
\]
where $\mathbb{E}_1$, $\mathbb{E}_2$, $\mathbb{E}_3$ are quadratic
polynomials of $\beta_{ij}$, difference
$\boldsymbol{U}^i-\boldsymbol{U}^j$ and $\omega^i$, respectively,
whose coefficients are independent of $\delta$.

Finally, we analyze the physical ramification of the obtained
asymptotic formula by presenting several examples. We construct an
a example of a suspension in a strong ``pinning'' external field,
where $\widehat{W}$ is of order $\delta^{-5/2}$ (superstrong blow
up). We also show an example of the superstrong blow up due to the
boundary layer effect. Note that to the best of our knowledge this
rate of blow up was not observed before and we call it an
\textit{anomalous rate}. For generic suspensions (free particles
or a weak external field) we expect that
$\boldsymbol{\beta}=O(\delta)$ and, therefore, the first and the
third terms of \eqref{E:intro} are of the same order
$O(\delta^{-1/2})$. For a hexagonal array of inclusions we prove
that $\widehat{W}$ exhibits the strong blow up (of order
$\delta^{-3/2}$) and $\boldsymbol{\beta}=\boldsymbol{0}$. Note that a typical close packing array in 2D is ``approximately'' hexagonal.

The paper is organized as follows.
In Section \ref{S:formulation} we give a mathematical formulation of the problem (Subsection \ref{SS:form}), describe the fictitious fluid approach,
and present our main results in Theorems \ref{T:main-thm1}, \ref{T:main-thm4} (Subsection \ref{SS:FFA}).
In Subsection \ref{S:constr-micro} we construct our discrete network and discuss how local flows in thin gaps between neighbors (microflows) contribute the overall
viscous dissipation rate and state the theorem about a representation of the error term of the discrete approximation.
In Section \ref{S:discussions} we discuss main and present examples.
Section \ref{SS:FFP} is devoted to the the fictitious fluid problem. In Section \ref{S:dn} we present results related to our discrete network. In section \ref{S:proof}
prove two lemmas about approximation of the overall dissipation rate by the dissipation rate of the fictitious fluid and about approximation of the
fictitious fluid (continuum) dissipation rate by the discrete one. Conclusions are presented in Section \ref{S:conclusions}.
The proofs of some auxiliary facts are given in Appendices.

\section*{Acknowledgements}
The authors thank A.G. Kolpakov and A. Panchenko for careful
reading of the manuscript and useful suggestions. The work of L.
Berlyand was supported by NSF Grant No. DMS-0204637. The work of
A. Novikov was supported by grants BSF-2005133 and DMS-0604600.

\Section{Formulation of the Problem and Main Results} \label{S:formulation}

\subsection{Mathematical Formulation of the Problem} \label{SS:form}

\begin{figure}[!ht]
  \centering
  \includegraphics[scale=1.05]{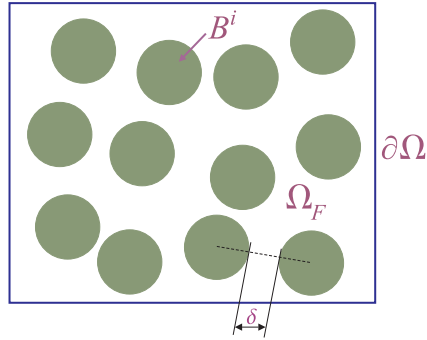}
  \caption{Domain $\Omega_F$ occupied by the fluid of viscosity $\mu$, and disordered array of closely spaced inclusions $B^i$}\label{F:domain}
\end{figure}

Consider an irregular or non-periodic array of $N$ identical circular disks $B^i$, of the radius $R$ distributed in a rectangular domain
$\Omega$. Denote by $\displaystyle \Omega_F=\Omega\setminus\bigcup_{i=1}^{N}B^i$ the {\it fluid domain} which is occupied by incompressible fluid with viscosity
$\mu$ (see Fig. \ref{F:domain}). Disks $B^i$ represent absolutely rigid inclusions. Inertia of both fluid and inclusions is neglected.
In the fluid domain $\Omega_{F}$ consider the following boundary value problem:
\begin{equation}   \label{E:v_form}
\left\{
\begin{array}{r l l}
(a) & \displaystyle \mu\triangle \boldsymbol{u}=\nabla p, & \boldsymbol{x} \in \Omega_F \\[5pt]
(b) & \displaystyle \nabla\cdot \boldsymbol{u}=0, &  \boldsymbol{x} \in \Omega_F\\[5pt]
(c) & \displaystyle \boldsymbol{u}=\boldsymbol{U}^{i}+R\omega^i(n_1^i\boldsymbol{e}_2-n_2^i\boldsymbol{e}_1), & \boldsymbol{x} \in \partial B^{i}, \quad i=1\ldots N \\[5pt]
(d) & \displaystyle \int_{\partial B^{i}}\boldsymbol{\sigma}(\boldsymbol{u})\boldsymbol{n}^i ds=\boldsymbol{0} & i=1\ldots N \\[10pt]
(e) & \displaystyle \int_{\partial B^{i}}\boldsymbol{n}^i\times\boldsymbol{\sigma}(\boldsymbol{u})\boldsymbol{n}^ids=0, & i=1\ldots N \\[5pt]
(f) & \displaystyle \boldsymbol{u}=\boldsymbol{f}, &  \boldsymbol{x} \in \partial\Omega
\end{array}
\right.
\end{equation}
where $\boldsymbol{u}(\boldsymbol{x})$ is the velocity field at a point $\boldsymbol{x} \in \Omega_F$, $p(\boldsymbol{x})$ is the
pressure, $\boldsymbol{\sigma}(\boldsymbol{u})=2\mu\boldsymbol{D}(\boldsymbol{u})-p\boldsymbol{I}$ is the stress tensor,
$\displaystyle D_{ij}(\boldsymbol{u})=\frac{1}{2}\left(\frac{\partial u_i}{\partial x_j}+\frac{\partial u_j}{\partial x_i}\right)$,
$i,j=1,2$, is the rate of strain which satisfies the incompressibility condition: $\displaystyle \mbox{tr}\,\boldsymbol{D}(\boldsymbol{u})=0$,
another form of (\ref{E:v_form}$b$). The vector $\boldsymbol{n}^i=(n_1^i,n_2^i)$ is the outer normal to $B^i$.
Constant vectors $\boldsymbol{U}^{i}=(U^i_1,U^i_2)$ and scalars $\omega^{i}$, $i=1,\ldots, N $, which are translational and
angular velocities of the inclusion $B^i$, respectively, are to be found in the course of solving the problem.

Here $N$ is closed to maximal packing number $N_{max}=N_{max}(\Omega,R)$. This number is finite and $|N-N_{max}|$ depends on the small parameter $\delta$ called
interpaticle distance rigorously defined below in Definition \ref{D:delta}.

\begin{figure}[!ht]
  \centering
  \includegraphics[scale=0.85]{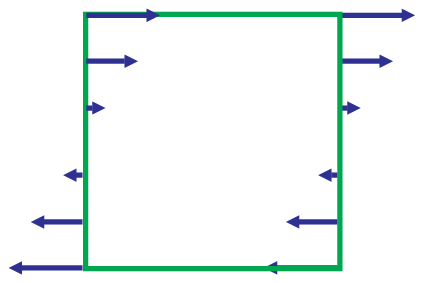} \hspace{2cm} \includegraphics[scale=0.85]{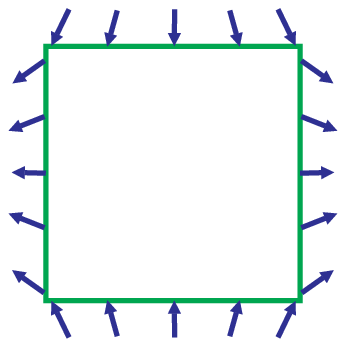}\\
  $(a)$ shear \hspace{3.5cm}$(b)$ extensional
  \caption{Shear and extensional external boundary conditions}\label{F:bc}
\end{figure}

We consider the linear external boundary conditions of the form:
\begin{equation}   \label{E:ext-bc}
\boldsymbol{f}=A\boldsymbol{x}=
\begin{pmatrix} a & b \\ c & -a \end{pmatrix}\begin{pmatrix} x \\ y \end{pmatrix}
\end{equation}
where the components $a$, $b$, $c$ of the matrix $A$ are given constants. Note that the most general form of the linear boundary conditions is
$\boldsymbol{f}=\boldsymbol{f}_0+A\boldsymbol{x}$ where $\boldsymbol{f}_0$ is a constant vector.
Observe that when $a=0$ and $b=-1/c$ the vector $A\boldsymbol{x}$ corresponds to a rotation and $\boldsymbol{f}_0$ to a translation of the
boundary, hence, $\boldsymbol{f}$ describes the rigid body motion of $\partial\Omega$.
Hereafter, we exclude this trivial motion from our consideration assuming that $a=0$ and $b=-1/c$ does not hold simultaneously in \eqref{E:ext-bc}
and $\boldsymbol{f}_0=\boldsymbol{0}$.
We use such boundary conditions for two reasons:
$a)$ for technical simplicity, which does not lead to the loss of generality; $b)$ they include the \textit{shear} (when $a=c=0$, $b=1$) and \textit{extensional} ($a=1$, $b=c=0$) boundary conditions
(see e.g. \cite{bp}) which model two basic types of viscometric measurements (see Fig. \ref{F:bc}).
It is possible to extend our results to arbitrary Dirichlet boundary conditions $\boldsymbol{f} \in H^{1/2}(\partial \Omega)$ satisfying
$\displaystyle \int_{\partial \Omega} \boldsymbol{f}\cdot\boldsymbol{n}ds=0$.

For an arbitrary set $A \subseteq \Omega_{F}$ consider  the following integral:
\begin{equation}   \label{E:W}
\begin{array}{r l r}
W_{A}(\boldsymbol{v})= & \displaystyle \frac{1}{2}\int_{A}\boldsymbol{\sigma}(\boldsymbol{v}):\boldsymbol{D}(\boldsymbol{v})=
\mu \int_{A} \boldsymbol{D}(\boldsymbol{v}):\boldsymbol{D}(\boldsymbol{v})d\boldsymbol{x}\\
=& \displaystyle \mu\int_{A} \left[\left(\frac{\partial v_1}{\partial x}\right)^{2}+ \frac{1}{2}\left(\frac{\partial v_1}{\partial y}+
\frac{\partial v_2}{\partial x}\right)^{2}+\left(\frac{\partial v_2}{\partial y}\right)^{2}\right]d\boldsymbol{x}
\end{array}
\end{equation}
where $\boldsymbol{v}=(v_1,v_2)$.
 Then the variational formulation of  \eqref{E:v_form} is:
\begin{equation}   \label{E:v_form-var}
\mbox{Find } \,\, \boldsymbol{u} \in V, \,\, \mbox{ such that } W_{\Omega_F}(\boldsymbol{u})=\min_{\boldsymbol{v} \in V}W_{\Omega_F}(\boldsymbol{v}),
\end{equation}
where the set of admissible vector fields $V$ is defined by
\begin{equation}   \label{E:V}
\begin{array}{l l}
\displaystyle V =& \displaystyle \left\{ \boldsymbol{v} \in \boldsymbol{H}^{1}(\Omega_F):\,\,
\nabla\cdot\boldsymbol{v}=0 \mbox{ in } \Omega_F,\,\,\boldsymbol{v}=\boldsymbol{f} \mbox{ on } \partial\Omega,\right.\\[5pt]
& \displaystyle \left. \boldsymbol{v}=\boldsymbol{U}^{i}+\boldsymbol{\omega}^{i}\times(\boldsymbol{x}-\boldsymbol{x}^{i}).
\,\,\,\boldsymbol{x}\in \partial B^{i}, \quad i=1,\ldots, N \right\},
\end{array}
\end{equation}
$W_{\Omega_F}(\boldsymbol{u})$ is called the (continuum) viscous dissipation rate \cite{landau} and it is the principal quantity of interest in the
study of overall properties of suspensions. We will use the following notation:
\begin{equation}   \label{E:eff-visc}
\widehat{W}:=W_{\Omega_F}(\boldsymbol{u}).
\end{equation}
The key feature of our problem is that we study suspensions where  concentration of inclusions is close to its maximum.
Therefore, the domain $\Omega_F$ depends on the characteristic interparticle distance parameter $\delta$.
Our main objective is to derive and justify an asymptotics of $\widehat{W}$ as $\delta\to 0$. We will show that the coefficients of
this asymptotic formula are determined by the solution to a discrete network problem, which determine the discrete viscous dissipation rate.
As explained in Introduction, $\widehat{W}$ determines a number of important measurable properties of suspensions.

\subsection{The Fictitious Fluid Approach and Discretization} \label{SS:FFA}

\begin{figure}[ht]
  \begin{center}
 \includegraphics[scale=0.99]{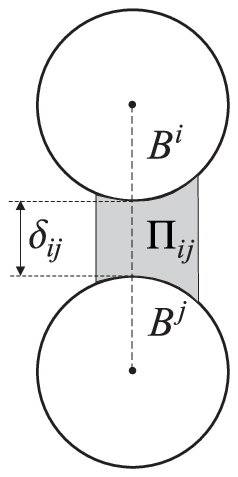} \hspace{2cm} \includegraphics[scale=0.99]{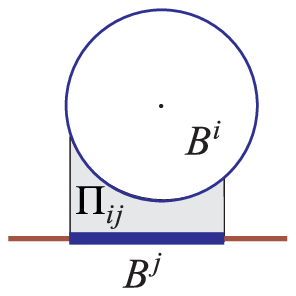}\\
 $(a)$ \hspace{4.5cm}
 $(b)$
  \end{center}
  \caption{$(a)$ neck connecting two disks;
 $(b)$ neck connecting disk $B^i$ and quasidisk $B^j$} \label{F:connect}
\end{figure}

\begin{figure}[ht]
  \begin{center}
 \includegraphics[scale=0.79]{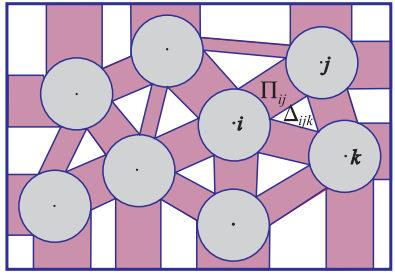}
  \end{center}
  \caption{The decomposition of the original fluid domain $\Omega_F$ into the fictitious  fluid domain and the set of triangles} \label{F:f_domain}
\end{figure}

In Appendix \ref{A:constr-in} we show that using the notion of Voronoi tessellation we can decompose the domain $\Omega_F$ into \textit{necks} $\boldsymbol{\Pi}$ and
\textit{triangles} $\boldsymbol{\Delta}$: $\Omega_F=\boldsymbol{\Pi}\cup\boldsymbol{\Delta}$ (see Fig. \ref{F:f_domain}). Necks connect either two disks (Fig. \ref{F:connect}$a$) or a disk and a part of the boundary $\partial\Omega$ called a \textit{quasidisk}
(see Fig. \ref{F:connect}$b$), that is, necks connect \textit{neighbors}. The velocities of quasidisks are given by the prescribed boundary conditions \eqref{E:ext-bc}.
Note that near the boundary when quasidisks are involved the ``triangles'' are actually trapezoids.
With slight abuse of terminology, we also call them \textit{triangles}.

We distinguish boundary disks (quasidisks) and interior disks and introduce two sets of the corresponding indices. For indices of interior disks we use the notation
$\mathbb{I}=\{1,\ldots,N\}$. If the disk $B^i$ centered at $\boldsymbol{x}_i$ is a quasidisk
then $i$ belongs to the set $\mathbb{B}$ of the indices of quasidisks. Also denote by $\mathcal{N}_i$ the set of indices of all neighbors of $B^i$.

%\begin{definition} \label{D:network}
For a given array of the disks and quasidisks $B^{i}$ centered at $\boldsymbol{x}_{i}$, the \textit{\textbf{discrete network}} is the graph
$\mathcal{G}=(\mathcal{X},\,\mathcal{E})$, with set of vertices $\mathcal{X}=\{\boldsymbol{x}_{i}:\,i \in \mathbb{I}\cup\mathbb{B}\}$ and set of edges
$\mathcal{E}=\{e_{ij}:\,i\in \mathbb{I},\,j \in \mathcal{N}_i\}$ with $e_{ij}$ connecting neighbors $B^i$ and $B^j$.
%\end{definition}

As mentioned in Introduction our main approach in study of the asymptotics of $\widehat{W}$ as $\delta\rightarrow 0$ consists of two steps.
This two-step approach allows to separate geometric construction of the network and its subsequent asymptotic analysis.

In the first step, we show that the minimization problem \eqref{E:v_form-var} in the fluid domain $\Omega_F$ can be
approximated by a ``fictitious fluid" problem in which fluid is assumed to occupy necks $\displaystyle \boldsymbol{\Pi}=\bigcup_{i\in \mathbb{I}, j \in \mathcal{N}_i}\Pi_{ij}$
between closely spaced neighboring inclusions (the shadowed region in Fig. \ref{F:f_domain}). We call $\boldsymbol{\Pi}$ the \textit{fictitious fluid domain}.
On the boundary of the complementary part of the domain (triangles in Fig. \ref{F:f_domain}) the relaxed incompressibility conditions:
\begin{equation} \label{E:weak-incompres}
\int_{\partial \triangle_{ijk}} \boldsymbol{v}\cdot\boldsymbol{n}ds=0, \quad i \in \mathbb{I}, \quad j,k \in \mathcal{N}_i,
\end{equation}
are imposed.
This reflects a well-known physical fact that for densely packed suspensions the dominant contribution to
the viscous dissipation rate over the fluid domain comes from those necks.

Below we show that the functional $\widehat{W}$ is decomposed as follows:
\begin{equation}   \label{E:main-decomp}
\widehat{W}=\widehat{W}_{\boldsymbol{\Pi}}+\widehat{W}_\Delta,
\end{equation}
where $\widehat{W}_{\boldsymbol{\Pi}}$ is the overall viscous dissipation rate of the fictitious fluid defined below in \eqref{E:W-fict} and $\widehat{W}_\Delta$ is the remainding
contribution from the domain $\boldsymbol{\Delta}$.

Consider the problem of minimization of the functional $W_{\boldsymbol{\Pi}}$, defined by \eqref{E:W} over the fictitious fluid domain $\boldsymbol{\Pi}$,
in the following class of functions:
\begin{equation}   \label{E:V_F}
\begin{array}{ l l}
& \displaystyle V_{\boldsymbol{\Pi}} = \left\{ \boldsymbol{v} \in \boldsymbol{H}^{1}(\boldsymbol{\Pi}): \right.
\nabla\cdot\boldsymbol{v}=0\,\,\mbox{in}\,\,\boldsymbol{\Pi},\,\,\int_{\partial\triangle_{ijk}}\boldsymbol{v}\cdot\boldsymbol{n}ds=0\,\,
\mbox{for all}\,\,\triangle_{ijk}\in \boldsymbol{\Delta},\\[10pt]
&  \displaystyle \left. \boldsymbol{v}=\boldsymbol{U}^{i}+ R\omega^i(n_1^i\boldsymbol{e}_2-n_2^i\boldsymbol{e}_1), \,
\boldsymbol{x} \in \partial B^{i}, \, i=1,\ldots, N, \,\,\boldsymbol{v}=\boldsymbol{f} \mbox{ on }\partial\boldsymbol{\Pi}\cap\partial\Omega\right\},
\end{array}
\end{equation}
where $\boldsymbol{U}^{i}$ and
$\omega^i$, $i=1,\ldots,N$, are arbitrary constant vectors and arbitrary constants, respectively.
The zero-flux condition through $\partial\triangle_{ijk}$ \eqref{E:weak-incompres} is inherited from the problem in the original fluid domain $\Omega_F$
due to incompressibility condition of the fluid in triangles. Such a condition is a necessary (but not sufficient) for $\nabla\cdot\boldsymbol{v}=0$ in the
triangle $\triangle_{ijk}$.

We define the overall dissipation rate of the fictitious fluid by
\begin{equation}   \label{E:W-fict}
\widehat{W}_{\boldsymbol{\Pi}}=\min_{\boldsymbol{v}\in V_{\boldsymbol{\Pi}}}W_{\boldsymbol{\Pi}}(\boldsymbol{v}).
\end{equation}

The first principal result of this paper is that the dissipation rate $\widehat{W}$ can be approximated by the rate $\widehat{W}_{\boldsymbol{\Pi}}$ and $\widehat{W}_\Delta$ can be neglected.
To show this we need to introduce a small parameter of the problem, which is a characteristic interparticle distance $\delta$.

For each pair of neighbors $\boldsymbol{x}_{i}$, $\boldsymbol{x}_{j}$ define
\begin{equation} \label{E:delta_ij}
\delta_{ij}=\left\{
\begin{array}{l l}
\displaystyle |\boldsymbol{x}_{i}-\boldsymbol{x}_{j}|-2R, \quad &  \mbox{when } i,\,j \in \mathbb{I}, \\
\displaystyle |\boldsymbol{x}_{i}-\boldsymbol{x}_{j}|-R,  \quad &  \mbox{when either } i \mbox{ or } j \in \mathbb{B}, \\
\displaystyle |\boldsymbol{x}_{i}-\boldsymbol{x}_{j}|,  \quad &  \mbox{when } i,\,j \in \mathbb{B}.
\end{array}
\right.
\end{equation}
As mentioned above, we study domains with closely spaced neighboring disks. More precisely, for all pairs of neighbors we
assume that following \textit{close-packing condition} holds.

\begin{definition} \label{D:delta}
Write the minimal distance $\delta_{ij}$ (see Fig. \ref{F:connect}) between any two neighboring disks $B^i$ and $B^j$ in the form
$\delta_{ij}=\delta d_{ij}$, where $d_{ij}$ is such that $0<c_1<d_{ij}<c_2$ for some absolute constants $c_1$, $c_2$. If the parameter $\delta$,
called the \textbf{\textit{characteristic interparticle distance}}, tends to zero,
then $\Omega_F$ is said to satisfy the \textbf{\textit{close-packing condition}}.
\end{definition}

We remark that this definition describes uniformly dense arrays of disks. A more general definition which covers a notion of a ``hole'' corresponding to the void space in the
composite is introduced and discussed e.g. in \cite{bn,bgn}.

Hereafter we call the array of inclusions under consideration a \textit{quasi-hexagonal array} (e.g. in \cite{bn} it is referred to as ``randomized hexagonal").
Recall, that for such arrays all neighbors are closely spaced and a typical number of nearest neighbors for a disk is six.

Hence, the mathematical thrust of step one of the fictitious fluid approach is in showing that the overall rate of the energy dissipation $\widehat{W}_{\boldsymbol{\Pi}}$ of the
fictitious fluid  captures the singular behavior of  $\widehat{W}$, defined by \eqref{E:eff-visc}, as $\delta \to 0$.
More precisely, in Section \ref{S:dn} we prove the following theorem:
\begin{theorem}[\textit{Approximation by the fictitious fluid}] \label{T:main-thm1}
Suppose an array of inclusions satisfies the close packing
condition. Let $\widehat{W}$ be the overall viscous dissipation
rate defined by \eqref{E:V}-\eqref{E:eff-visc} and $\widehat{W}_{\boldsymbol{\Pi}}$
be the viscous dissipation rate of the fictitious fluid defined by
\eqref{E:W-fict}. Then the following asymptotic formula holds:
\begin{equation}   \label{E:lemma-ineq1}
\frac{|\widehat{W}-\widehat{W}_{\boldsymbol{\Pi}}|}{\widehat{W}_{\boldsymbol{\Pi}}}=O(\delta^{1/2}) \quad \mbox{as } \delta\rightarrow 0.
\end{equation}
\end{theorem}

In step two, we study asymptotics (blow up) of the overall viscous dissipation rate $\widehat{W}$ as $\delta\rightarrow 0$. In view
of step one, this is reduced to finding of asymptotics of $\widehat{W}_{\boldsymbol{\Pi}}$. The latter is done by a construction of a discrete network
approximation and introduction of a so-called discrete viscous dissipation rate $\mathcal{I}$. To show closeness
of the continuum and the discrete dissipation rates, $\widehat{W}_{\boldsymbol{\Pi}}$ and $\mathcal{I}$, respectively,
we employ the direct and dual variational techniques~\cite{bk,bn,bbp}.

Also note that the conditions (\ref{E:v_form}$d$,$e$) in the original problem led to significant technical difficulties in
variational analysis of the overall viscous dissipation rate of \cite{bbp}, which is why the analysis of \cite{bbp} is restricted to its
leading singular term. In contrast, analogs of these conditions in the fictitious fluid problem are satisfied {\it automatically} by
construction, which results in substantial simplification of the analysis and thus allows to capture all singular terms.

The approximation of the overall viscous dissipation rate $\widehat{W}$ by the discrete dissipation rate $\mathcal{I}$ is given by the following theorem.
\begin{theorem}[\textit{Approximation of the continuum dissipation rates by the discrete one}] \label{T:main-thm4}
Suppose $\Omega_F$ satisfies the close packing condition and
$\displaystyle \mathcal{I}=\min_{(\mathbb{U},\boldsymbol{\omega},\boldsymbol{\beta})\in \boldsymbol{\mathcal{R}}}Q(\mathbb{U},\boldsymbol{\omega},\boldsymbol{\beta})$,
where the positive definite quadratic form
\[
Q(\mathbb{U},\boldsymbol{\omega},\boldsymbol{\beta})=\mathcal{I}_1(\boldsymbol{\beta})\delta^{-5/2}+\mathcal{I}_2(\mathbb{U},\boldsymbol{\omega},\boldsymbol{\beta})\delta^{-3/2}+
\mathcal{I}_3(\mathbb{U},\boldsymbol{\omega},\boldsymbol{\beta})\delta^{-1/2}
\]
on the class of admissible discrete variables $(\mathbb{U},\boldsymbol{\omega},\boldsymbol{\beta})\in \boldsymbol{\mathcal{R}}$ defined in \eqref{E:constraints-new}-\eqref{E:Q_ij-bound}.
Then the following approximation to the viscous dissipation rate holds:
\begin{equation} \label{E:main-thm4}
\frac{|\widehat{W}-\mathcal{I}|}{\mathcal{I}}=O(\delta^{1/2}) \quad \mbox{as } \delta\rightarrow 0,
\end{equation}
\end{theorem}

\begin{remark}
$\mathcal{I}_i$ ($i=1,2,3$) are explicitly computable quadratic polynomials of $(\mathbb{U},\boldsymbol{\omega},\boldsymbol{\beta})\in \boldsymbol{\mathcal{R}}$
from equations \eqref{E:Q_ij}-\eqref{E:Q_ij-bound} below. These polynomials depend only on boundary data $\boldsymbol{f}$, viscosity $\mu$ and geometry of $\Omega_F$.
\end{remark}

The next subsection is devoted to the construction of these quadratic polynomials. We first
introduce a set of discrete variables $(\mathbb{U},\boldsymbol{\omega},\boldsymbol{\beta})$, define a new variable $\beta_{ij}$ in each neck $\Pi_{ij}$
and explain how the quadratic form $Q$ is obtained. We will also discuss the underlying structure of the flow in a neck
(microflow) and physical ramifications of the obtained asymptotics \eqref{E:main-thm4} and \eqref{E:main-thm3} below.

\subsection{Construction of the discrete network. Microflows.} \label{S:constr-micro}

$(a)$ We begin with the discretization of the boundary conditions.
Denote $\mathbb{U}=\{\boldsymbol{U}^i\}_{i\in \mathbb{I}\cup\mathbb{B}}$, $\boldsymbol{\omega}=\{\omega^i\}_{i\in \mathbb{I}\cup\mathbb{B}}$
where on the boundary $\partial\Omega$, that is, for $i\in\mathbb{B}$, these velocities are given by boundary conditions as follows:
\begin{equation} \label{E:boundary-condtns}
\begin{array}{l l l l }
& \boldsymbol{U}^i=
\begin{cases}
\begin{pmatrix} ax \\ cx \end{pmatrix}, & \boldsymbol{x}_i \in \partial\Omega_{lat} \\[10pt]
\begin{pmatrix} by \\ -ay \end{pmatrix}, & \boldsymbol{x}_i \in \partial\Omega^\pm
\end{cases} & \mbox{ and }
\omega^i=
\begin{cases}
c, & \boldsymbol{x}_i \in \partial \Omega^- \\
-c, & \boldsymbol{x}_i \in \partial \Omega^+ \\
b, & \boldsymbol{x}_i \in \partial \Omega^-_{lat} \\
-b, & \boldsymbol{x}_i \in \partial \Omega^+_{lat} \\
\end{cases}
\end{array}
\end{equation}
where $\boldsymbol{x}_i$ is the center of the quasidisk $B^i$ and
$\partial\Omega^+$ and $\partial\Omega^-$ are the upper and lower parts of $\partial\Omega$, respectively, and $\partial\Omega^+_{lat}$ and $\partial\Omega^-_{lat}$ are
the lateral (left and right) boundary.

\begin{figure}[!ht]
  \centering
  \includegraphics[scale=0.85]{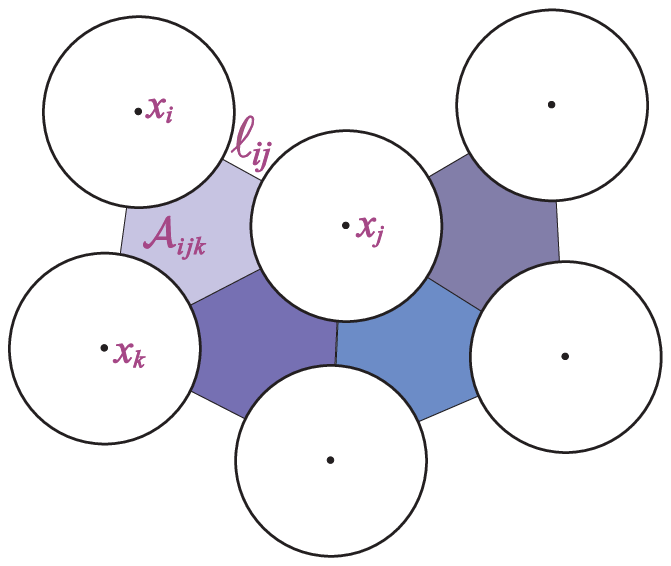} \hspace{1cm} \includegraphics[scale=1.35]{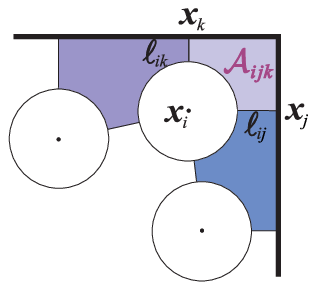}\\
  $(a)$ \hspace{5.5cm} $(b)$\\
  \caption{$(a)$ Decomposition of $\Omega_F$ into curvilinear hexagons $\mathcal{A}_{ijk}$ and line $\ell_{ij}$ connecting neighbors $B^i$ and $B^j$,
  $(b)$ Construction of $\mathcal{A}_{ijk}$ at the corner of the boundary $\partial\Omega$
}\label{F:six1}
\end{figure}

$(b)$ Discretization of the incompressibility condition is implemented as follows. Decompose the domain $\Omega_F$ into curvilinear hexagons $\mathcal{A}_{ijk}$
as in Fig. \ref{F:six1}:
\[
\Omega_F=\bigcup_{i\in \mathbb{I}, j,k\in \mathcal{N}_i}\mathcal{A}_{ijk}.
\]
Each $\mathcal{A}_{ijk}$ consists of the line segments $\ell_{ij}$, $\ell_{jk}$, $\ell_{ki}$ (Fig. \ref{F:six1}$a$) connecting disks $B^i$, $B^j$, $B^k$ and arcs
$a_i$, $a_j$, $a_k$ of the corresponding disk.\footnote{In the case when the disk $B^i$ has two quasidisk neighbors $B^{j}$ and $B^{k}$, that is,
when $B^k$ is in the corner as in Fig. \ref{F:six1}$b$,
then the domain $\mathcal{A}_{ijk}$ is actually a curvilinear pentagon. By slight abuse of terminology we still call it a ``curvilinear hexagon''.}
Then the weak incompressibility condition \eqref{E:weak-incompres} for the class $V_{\boldsymbol{\Pi}}$ \eqref{E:V_F} becomes
\begin{equation}  \label{E:zero-flux-six}
\int_{\mathcal{A}_{ijk}}\boldsymbol{u}\cdot\boldsymbol{n}ds=0, \quad \mbox{ for any } \mathcal{A}_{ijk}, \quad i\in \mathbb{I}, j,k \in \mathcal{N}_i.
\end{equation}

\begin{figure}[!ht]
  \centering
  \includegraphics[scale=0.95]{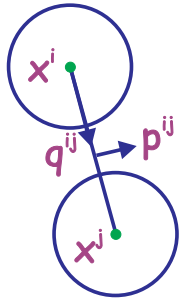}
  \caption{Unit vectors of the local coordinate system}\label{F:pq}
\end{figure}

In order to continue our analysis at this point we must introduce a new set of discrete variables. Here we define permeation constants:
\begin{equation}  \label{E:beta}
\begin{array}{l l l}
\displaystyle \beta_{ij}^*=\frac{1}{R}\int_{\ell_{ij}}\boldsymbol{u}\cdot\boldsymbol{n}ds, & i\in \mathbb{I}, \,\, j \in \mathcal{N}_i, \\[5pt]
\displaystyle \beta_{ij}^*=\frac{1}{R}\int_{\ell_{ij}}\boldsymbol{f}\cdot\boldsymbol{n}ds, & \ell_{ij}\subset\partial\Omega\,(i,j\in \mathbb{B}),
\end{array}
\end{equation}
where $\ell_{ij}$ is the line segment joining two neighbors $B^i$ and $B^j$, $\boldsymbol{u} \in V_{\boldsymbol{\Pi}}$ and $\boldsymbol{n}$ is an outer normal to $\mathcal{A}_{ijk}$.
Then \eqref{E:zero-flux-six} can be rewritten as:
\[
\beta_{ij}^*+\beta_{jk}^*+\beta_{ki}^*+\frac{1}{R}\int_{a_{i}}\boldsymbol{U}^i\cdot\boldsymbol{n}^ids+
\frac{1}{R}\int_{a_{j}}\boldsymbol{U}^j\cdot\boldsymbol{n}^jds+\frac{1}{R}\int_{a_{k}}\boldsymbol{U}^k\cdot\boldsymbol{n}^kds=0,
\]
for $i\in \mathbb{I}$, $j,k\in\mathcal{N}_i$ ($\boldsymbol{n}^i$ is a unit outer normal to $\partial B^i$), which can be further simplified as
\begin{equation} \label{E:disk-constraint}
\beta_{ij}^*+\beta_{jk}^*+\beta_{ki}^*+(\boldsymbol{U}^i+\boldsymbol{U}^j)\boldsymbol{p}^{ij}+(\boldsymbol{U}^j+\boldsymbol{U}^k)\boldsymbol{p}^{jk}+(\boldsymbol{U}^k+\boldsymbol{U}^i)\boldsymbol{p}^{ki}=0,
\end{equation}
where vectors $\boldsymbol{q}^{ij}$ and $\boldsymbol{p}^{ij}$ are the unit vectors of the local system coordinate of two neighboring
disks $B^i$ and $B^j$ as in Fig. \ref{F:pq}.
We call \eqref{E:disk-constraint} the \textit{weak incompressibility condition}.
This formula explains the scaling $\frac{1}{R}$ in \eqref{E:beta}. Indeed, from dimensional analysis $\int_{\ell_{ij}}\boldsymbol{u}\cdot\boldsymbol{n}ds$ must be divided by a
characteristic lengthscale which is $R$ in the curvilinear hexagonal $\mathcal{A}_{ijk}$.

$(c)$ We finally discretize the Stokes equations in necks. This discretization is based on lubrication theory \cite{ach,leal,batchelor}.
This theory describes thin film flows between two solid bodies sliding relative to each other.
It is well-known (see e.g. \cite{leal}) that in dense suspensions the dominant hydrodynamic contribution to the viscous dissipation rate occurs in necks
between closely spaced inclusions, where lubrication equations are relevant. Lubrication theory determines fluid motion in such necks as a result of relative
kinematic motion of the neighboring inclusions. This raises a question of the classification of \textit{microflows}, that is,
local flows in necks between two closely spaced neighbors.

Recall three classical types of microflows between two parallel plates: the shear, the squeeze, and the Poiseuille flow.
The last one, however, is not related to motions of two plates relative one to another and therefore it is not described by lubrication theory. In this work we show
that exactly these three types of microflows fully describe the motion of the fluid between two neighboring disks.
Asymptotic analysis of microflows between two parallel plates technically is much simpler than that for inclusions with curvilinear boundaries,
which is needed for suspensions.

As shown below, in classical 3D problem the Poiseuille microflow between two inclusions does not contribute to the singular behavior of the viscous dissipation rate.
However, in analogous 2D problem (e.g. thin films) all three microflows are present and, moreover, the Poiseuille flow results in anomalously strong singularity.
Note that in the previous studies of the overall properties of suspensions the Poiseuille type microflow was not taken into account.

Thus, it is necessary to analyze kinematics of a pair of neighboring inclusions when one moves relative to the other.
To this end for a pair of neighbors $B^i$ and $B^j$, centered at $\boldsymbol{x}_{i}$ and $\boldsymbol{x}_{j}$, we choose the
\textit{local} coordinate system where the origin is at $(\boldsymbol{x}_i+\boldsymbol{x}_j)/2$ and the $y$-axis is
directed along the vector connecting $\boldsymbol{x}_{i}$ and $\boldsymbol{x}_{j}$.

For clarity of presentation we consider two interior disks only (that is, $i\in \mathbb{I}$, $j\in \mathcal{N}_i\cap\mathbb{I}$).
An analogous construction in boundary necks ($i\in \mathbb{I}$, $j\in \mathcal{N}_i\cap\mathbb{B}$) is given in Appendix \ref{A:bound}.

\begin{figure}[!ht]
  \centering
  \includegraphics[scale=0.80]{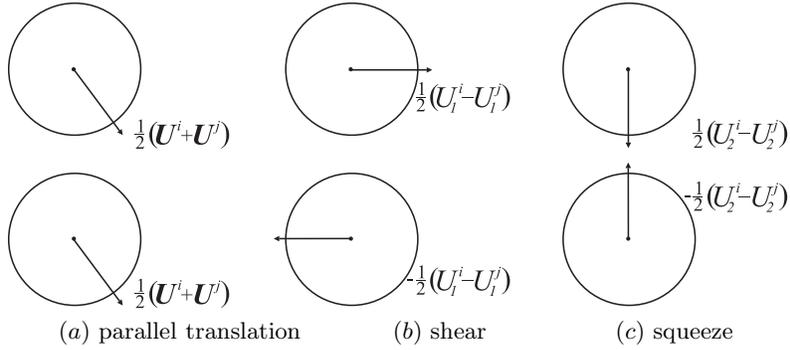}\\
  ($a$) parallel translation \hspace{1.cm} ($b$) shear \hspace{1.5cm} ($c$) squeeze
  \caption{Decomposition of translational velocities $\boldsymbol{U}^i$ and $\boldsymbol{U}^j$ into three elementary motions}\label{F:bc1}
\end{figure}

There are exactly five elementary kinematic motions of inclusions. To see this we consider boundary conditions on $\partial B^i$, $\partial B^j$ in
\eqref{E:V_F} and first assume that $\omega^i=\omega^j=0$. Then the conditions:
$\boldsymbol{u}=\boldsymbol{U}^i$ on $\partial B^{i}$ and $\boldsymbol{u}=\boldsymbol{U}^j$ on $\partial B^{j}$ can be rewritten as follows:
\[
\begin{array}{l l }
& \displaystyle \boldsymbol{u}|_{\partial B_i}=\frac{1}{2}(\boldsymbol{U}^i+\boldsymbol{U}^j)+\frac{1}{2}(U^i_1-U^j_1)\boldsymbol{e}_1+\frac{1}{2}(U^i_2-U^j_2)\boldsymbol{e}_2,\\[5pt]
\mbox{and  }& \displaystyle \boldsymbol{u}|_{\partial B^j}=\frac{1}{2}(\boldsymbol{U}^i+\boldsymbol{U}^j)-\frac{1}{2}(U^i_1-U^j_1)\boldsymbol{e}_1-\frac{1}{2}(U^i_2-U^j_2)\boldsymbol{e}_2.
\end{array}
\]
Hence, the translational velocities of disks are decomposed into three motions. First, when both inclusions and fluid move with the same velocity (see Fig. \ref{F:bc1}$a$),
and therefore this motion does not contribute to the singular behavior of the viscous dissipation rate. Second, when one inclusion moves
relative to the other producing the \textit{shear} type motion (Fig. \ref{F:bc1}$b$) and, finally, the \textit{squeeze} type motion (Fig. \ref{F:bc1}$c$) of the fluid.

\begin{figure}[!ht]
  \centering
  \includegraphics[scale=0.95]{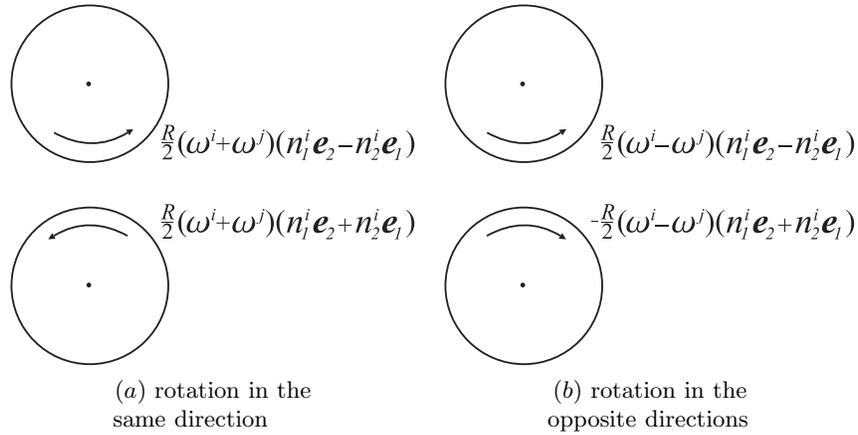}\\
  ($a$) rotation in the \hspace{3.0cm} ($b$) rotation in the\\
  same direction \hspace{3.5cm} opposite directions
  \caption{Decomposition of the angular velocities $\omega^i$ and $\omega^j$ into two elementary motions}\label{F:bc2}
\end{figure}

Similarly, assuming $\boldsymbol{U}^i=\boldsymbol{U}^j=\boldsymbol{0}$ in conditions on $\partial B^i$, $\partial B^j$ in
\eqref{E:V_F} we decompose them into two relative elementary motions:
rotations of the disks in the same direction (Fig. \ref{F:bc2}$a$), and rotations of the disks in the opposite directions (Fig. \ref{F:bc2}$b$). Thus,
\[
\begin{array}{l l }
& \displaystyle \boldsymbol{u}|_{\partial B_i}=\frac{R}{2}(\omega^i+\omega^j)(n_1^i \boldsymbol{e}_{2}-n_2^i \boldsymbol{e}_{1})+ \frac{R}{2}(\omega^i-\omega^j)(n_1^i \boldsymbol{e}_{2}-n_2^i \boldsymbol{e}_{1})\\[5pt]
\mbox{and  }& \displaystyle \boldsymbol{u}|_{\partial B^j}=\frac{R}{2}(\omega^i+\omega^j)(n_1^i\boldsymbol{e}_{2}+n_2^i \boldsymbol{e}_{1})+\frac{R}{2}(\omega^i-\omega^j)(-n_1^i \boldsymbol{e}_{2}-n_2^i\boldsymbol{e}_{1}),
\end{array}
\]
where $\boldsymbol{n}^{i}=(n_1^i,n_2^i)$ is the outer normal to $\partial B^i$.

\begin{figure}[!ht]
  \centering
  \includegraphics[scale=1.35]{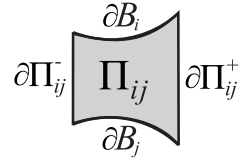}
  \caption{Boundary of the neck $\Pi_{ij}$}\label{F:neck4}
\end{figure}

Next consider microflows corresponding to the four kinematic motions contributing to the singular behavior of the dissipation rate described above.
We further choose the microflow $\boldsymbol{u}_{lub}$ which minimizes the viscous dissipation rate in the neck $\Pi_{ij}$ by imposing the
natural boundary conditions on the lateral boundaries $\partial\Pi^{\pm}_{ij}$ (Fig. \ref{F:neck4}). That is, the function $\boldsymbol{u}_{lub}$ minimizes
$W_{\Pi_{ij}}$ in the class:
\begin{equation} \label{E:V_ij-curl}
\begin{array}{r l l l }
\mathcal{V}_{ij}= & \displaystyle \left\{\boldsymbol{v}\in H^1(\Pi_{ij}):\,\,\bigtriangledown\cdot\boldsymbol{v}=0 \,\, \mbox{in}\,\,\Pi_{ij},\,\,
\boldsymbol{v}=\boldsymbol{g}_i  \,\, \mbox{on}\,\,\partial B^i, \,\,
\boldsymbol{v}=\boldsymbol{g}_j  \,\, \mbox{on}\,\,\partial B^j\right\}
\end{array}
\end{equation}
\begin{equation} \label{E:bc-g}
\mbox{with }\,\,\boldsymbol{g}_i=\boldsymbol{U}^{i}+R\omega^i(n_1^i\boldsymbol{e}_2-n_2^i\boldsymbol{e}_1), \quad \boldsymbol{g}_j=\boldsymbol{U}^{j}+R\omega^j(n_1^j\boldsymbol{e}_2-n_2^j\boldsymbol{e}_1).
\end{equation}

Due to linearity the minimizer $\boldsymbol{u}_{lub}$ is decomposed into five vector fields corresponding to the relative motions of inclusions
as in Fig. \ref{F:bc1}-\ref{F:bc2} as follows:
\begin{equation} \label{E:u-lub}
\begin{array}{r l l l }
\boldsymbol{u}_{lub} =  & \displaystyle \frac{1}{2}(\boldsymbol{U}^i+\boldsymbol{U}^j)+[(\boldsymbol{U}^i-\boldsymbol{U}^j)\cdot\boldsymbol{p}^{ij}]\boldsymbol{u}_{1}+
[(\boldsymbol{U}^i-\boldsymbol{U}^j)\cdot\boldsymbol{q}^{ij}]\boldsymbol{u}_{2}\\[5pt]
+  & \displaystyle R(\omega^i+\omega^j)\boldsymbol{u}_{3}+R(\omega^i-\omega^j)\boldsymbol{u}_{4},
\end{array}
\end{equation}
where functions $\boldsymbol{u}_k$, $k=1,\ldots,4$ are minimizers of $W_{\Pi_{ij}}$ in $\mathcal{V}_{ij}$ (equation \eqref{E:V_ij-curl}) where boundary conditions \eqref{E:bc-g} are replaced, respectively, by:\\
$1)$ the \textit{shear} motion of the fluid between two neighboring inclusions:
\begin{equation} \label{E:bc-1}
\boldsymbol{u}_1|_{\partial B^i}=\boldsymbol{g}_i=\frac{1}{2}\boldsymbol{e}_1, \quad \boldsymbol{u}_1|_{\partial B^j}=\boldsymbol{g}_j=-\frac{1}{2}\boldsymbol{e}_1,
\end{equation}
$2)$ the \textit{squeeze} motion of the fluid between neighbors:
\begin{equation} \label{E:bc-2}
\boldsymbol{u}_2|_{\partial B^i}=\boldsymbol{g}_i=\frac{1}{2}\boldsymbol{e}_2, \quad \boldsymbol{u}_2|_{\partial B^j}=\boldsymbol{g}_j=-\frac{1}{2}\boldsymbol{e}_2,
\end{equation}
$3)$ the \textit{rotation in the same directions}:
\begin{equation} \label{E:bc-3}
\boldsymbol{u}_3|_{\partial B^i}=\boldsymbol{g}_i=\frac{1}{2}(n_1^i\boldsymbol{e}_2-n_2^i\boldsymbol{e}_1), \quad
\boldsymbol{u}_3|_{\partial B^j}=\boldsymbol{g}_j=\frac{1}{2}(n_1^i\boldsymbol{e}_2+n_2^i\boldsymbol{e}_1),
\end{equation}
$4)$ the \textit{rotation in the opposite directions}:
\begin{equation} \label{E:bc-4}
\boldsymbol{u}_4|_{\partial B^i}=\boldsymbol{g}_i=\frac{1}{2}(n_1^i\boldsymbol{e}_2-n_2^i\boldsymbol{e}_1), \quad
\boldsymbol{u}_4|_{\partial B^j}=\boldsymbol{g}_j=-\frac{1}{2}(n_1^i\boldsymbol{e}_2+n_2^i\boldsymbol{e}_1).
\end{equation}

\begin{figure}[!ht]
  \centering
  \includegraphics[scale=0.95]{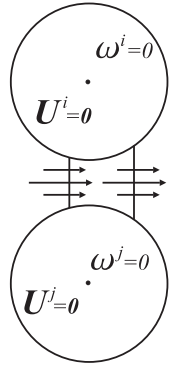}
  \caption{The Poiseuille's microflow between motionless inclusions}\label{F:bc3}
\end{figure}

Let $\boldsymbol{u}$ be the minimizer of \eqref{E:v_form-var} then
\begin{equation} \label{E:u-5}
\boldsymbol{u}_{p}:=\boldsymbol{u}-\boldsymbol{u}_{lub}
\end{equation}
minimizes $W_{\Pi_{ij}}$ in the class
\begin{equation}   \label{E:V_ij-p}
\begin{array}{r l l}
\mathcal{V}_{p}=& \displaystyle \left\{\boldsymbol{v} \in H^{1}(\Pi_{ij}):\,\,\nabla\cdot\boldsymbol{v}=0 \,\,\mbox{in}\,\,\Pi_{ij},\,\,
\boldsymbol{v}=\boldsymbol{0}\,\,\mbox{on}\,\, \partial B^{i},%\right.\\[7pt]
%& \displaystyle \left.
\,\,\boldsymbol{v}=\boldsymbol{0}\,\,\mbox{on}\,\, \partial B^{j}, \right. \\[7pt]
& \displaystyle \left.\,\,\frac{1}{R}\int_{\ell_{ij}}(\boldsymbol{v}+\boldsymbol{u}_{lub})\cdot\boldsymbol{n}ds=\beta_{ij}^* \right\}.
\end{array}
\end{equation}
The vector field $\boldsymbol{u}_p$ describes the flow between two motionless inclusions, hence, it is natural to call it the \textit{Poiseuille microflow}.

Denote by $\boldsymbol{\beta}^*=\{\beta_{ij}^*\}_{i\in \mathbb{I}\cup\mathbb{B},j\in \mathcal{N}_i}$ and introduce the following set of discrete variables
\begin{equation} \label{E:constraints}
\begin{array}{l l l}
\boldsymbol{\mathcal{R}}^*= & \displaystyle\left\{(\mathbb{U},\boldsymbol{\omega},\boldsymbol{\beta}^*):
\boldsymbol{U}^i,\omega^i \mbox{ satisfying } \eqref{E:boundary-condtns} \mbox{ for } i\in\mathbb{B}\right.\\[5pt]
& \displaystyle \left. \,\,\mbox{and}\,\, (\mathbb{U},\boldsymbol{\beta}^*) \mbox{ satisfying } \eqref{E:disk-constraint}\right\}.
\end{array}
\end{equation}
It is straightforward to show that (see Lemma \ref{L:two-min})
\begin{equation} \label{E:min-chain_2}
\widehat{W}_{\boldsymbol{\Pi}}=\min_{(\mathbb{U},\boldsymbol{\omega},\boldsymbol{\beta}^*)\in \boldsymbol{\mathcal{R}}^*}\sum_{i\in\mathbb{I}, j\in \mathcal{N}_i}\min_{V_{ij}} W_{\Pi_{ij}}(\cdot),
\end{equation}
where $V_{ij}$ is defined by
\begin{equation}   \label{E:V_ij}
\begin{array}{l l l}
V_{ij}=& \displaystyle \left\{\boldsymbol{v} \in H^{1}(\Pi_{ij}):\,\,\nabla\cdot\boldsymbol{v}=0
\,\,\mbox{in}\,\,\Pi_{ij},\,\,
\boldsymbol{v}=\!\boldsymbol{U}^{i}\!\!+R\omega^j(n_1^i\boldsymbol{e}_2\!-\!n_2^i\boldsymbol{e}_1) \,\,\mbox{on}\,\, \partial B^{i},\right.\\[7pt]
& \displaystyle
\left.\,\,\boldsymbol{v}\!=\!\boldsymbol{U}^{j}\!\!+\!\!R\omega^i(n_1^j\boldsymbol{e}_2-n_2^j\boldsymbol{e}_1)
\,\,\mbox{on}\,\, \partial B^{j}
\,\,(j\in\mathbb{I}) \mbox{ or } \boldsymbol{v}=\boldsymbol{f}\,\, \mbox{on}\,\, \partial B^j \,\,(j\in\mathbb{B})\right. \\[7pt]
& \displaystyle \left.\,\,
\frac{1}{R}\int_{\ell_{ij}}\boldsymbol{v}\cdot\boldsymbol{n}ds=\beta_{ij}^*\right\},
\end{array}
\end{equation}
with $\boldsymbol{U}^{i}$, $\omega^i$ and $\beta_{ij}^*$ to be components of some no longer arbitrary but fixed
$(\mathbb{U},\boldsymbol{\omega},\boldsymbol{\beta}^*)\in \boldsymbol{\mathcal{R}}^*$.
Functions from $V_{ij}$ are defined in a single neck $\Pi_{ij} \in \boldsymbol{\Pi}$.

By direct computations vector fields
\begin{equation} \label{E:various-u}
\begin{array}{r l}
\boldsymbol{u}_{t} & \displaystyle =\frac{1}{2}(\boldsymbol{U}^i+\boldsymbol{U}^j),\\
\boldsymbol{u}_{sh} & \displaystyle =[(\boldsymbol{U}^i-\boldsymbol{U}^j)\cdot\boldsymbol{p}^{ij}]\boldsymbol{u}_{1}+R(\omega^i+\omega^j)\boldsymbol{u}_{3},\\
\boldsymbol{u}_{sq} & \displaystyle =[(\boldsymbol{U}^i-\boldsymbol{U}^j)\cdot\boldsymbol{q}^{ij}]\boldsymbol{u}_{2},\\
\boldsymbol{u}_{per} & \displaystyle =R(\omega^i-\omega^j)\boldsymbol{u}_{4}+\boldsymbol{u}_{p},
\end{array}
\end{equation}
are orthogonal with respect to the scalar product induced by the dissipation functional $W_{\Pi_{ij}}$:
\begin{equation} \label{E:W-decomp}
\begin{array}{r l}
W_{\Pi_{ij}}(\boldsymbol{u})=& W_{\Pi_{ij}}(\boldsymbol{u}_{sh})+W_{\Pi_{ij}}(\boldsymbol{u}_{sq})+W_{\Pi_{ij}}(\boldsymbol{u}_{per}), \quad W_{\Pi_{ij}}(\boldsymbol{u}_{t})=0,\\[5pt]
\boldsymbol{u}=& \displaystyle \boldsymbol{u}_{t}+\boldsymbol{u}_{sh}+\boldsymbol{u}_{sq}+\boldsymbol{u}_{per}.
\end{array}
\end{equation}

Physically, the decomposition \eqref{E:W-decomp} corresponds to three well-known types of microflows. Namely,
($i$) the shear type arises when a pair of inclusions either rotates in the same direction or disks move into opposite directions (Fig. \ref{F:bc2}$a$, \ref{F:bc1}$a$),
($ii$) the squeeze type, when two inclusions in a pair move towards or away each other in thin gaps (Fig. \ref{F:bc1}$b$),
($iii$) and permeation of the fluid through the thin gaps between neighbors due to Poiseuille flow between motionless inclusions
or rotation of neighbors in the opposite directions (Fig. \ref{F:bc3}, \ref{F:bc2}$b$).

We are now ready to introduce a quadratic form which determines the discrete dissipation rate. To this end from now on instead of $\beta_{ij}^*$ as a discrete variable we will use:
\begin{equation} \label{E:beta0-new}
\beta_{ij}=\beta_{ij}^*-\frac{\delta_{ij}}{2R}\left[(\boldsymbol{U}^i+\boldsymbol{U}^j)\cdot\boldsymbol{p}^{ij}\right]
-(\omega^i-\omega^j)\frac{\delta_{ij}}{2R}\left[1+\frac{\delta_{ij}}{4R}\right].
\end{equation}
The reason for this replacement is that $\beta_{ij}$ is invariant with respect to Galilean transformation whereas $\beta_{ij}^*$ is not
(see remark in the end of Subsection \ref{SS:proof_ff-approx}).
For example, if a constant vector $\boldsymbol{U}^0$ is added to both $\boldsymbol{U}^i$ and $\boldsymbol{U}^j$
then the total flux $\beta_{ij}^*$ changes while $\beta_{ij}$ stays the same.
Also, $\beta_{ij}^*$ is the total flux through $\ell_{ij}$ of the entire flow (including the parallel translation, shear, squeeze and permeation) whereas
$\beta_{ij}$ is the flux due to the Poiseuille microflow solely. Finally, the use of $\beta_{ij}$ simplifies
the discrete dissipation form (equation \eqref{E:Q_ij}).

The use of $\beta_{ij}$ instead of $\beta_{ij}^*$ leads to the replacement of the class $\boldsymbol{\mathcal{R}}^*$ \eqref{E:constraints} by $\boldsymbol{\mathcal{R}}$ defined as follows:
\begin{equation} \label{E:constraints-new}
\begin{array}{l l l}
\boldsymbol{\mathcal{R}}\!=\!\! & \displaystyle\left\{\!(\mathbb{U},\boldsymbol{\omega},\boldsymbol{\beta})\!:
\boldsymbol{U}^i\!,\omega^i \mbox{ satisfying boundary condition on } \partial\Omega \,\, \eqref{E:boundary-condtns} \mbox{ for } i\in\mathbb{B},\right.\\[5pt]
& \displaystyle \left. \,\,\, (\mathbb{U},\boldsymbol{\beta}) \mbox{ satisfying weak incompresibility condition} \eqref{E:disk-constraint}, \,\eqref{E:beta0-new} \right\}.
\end{array}
\end{equation}

Introduce the \textit{overall discrete dissipation rate}:
\begin{equation} \label{E:I}
\mathcal{I}:=Q(\mathbb{\widehat{U}},\boldsymbol{\widehat{\omega}},\boldsymbol{\widehat{\beta}})=\min_{(\mathbb{U},\boldsymbol{\omega},\boldsymbol{\beta})\in
\boldsymbol{\mathcal{R}}}Q(\mathbb{U},\boldsymbol{\omega},\boldsymbol{\beta}),
\end{equation}
where
\begin{equation} \label{E:Q}
Q(\mathbb{U},\boldsymbol{\omega},\boldsymbol{\beta})=\sum_{i\in \mathbb{I}}\sum_{j\in\mathcal{N}_i} Q_{ij},
\end{equation}
\begin{equation} \label{E:Q_ij}
\begin{array}{r l l}
\displaystyle
Q_{ij}(\boldsymbol{U}^i,\boldsymbol{U}^j,\omega^{i},\omega^{j},\beta_{ij})
=&
\displaystyle \left[(\boldsymbol{U}^i-\boldsymbol{U}^j)\cdot\boldsymbol{p}^{ij}+ R\omega^i+R\omega^j\right]^2\boldsymbol{\mathcal{C}}_{1}^{ij}\delta^{-1/2}\\[7pt]
+& \displaystyle \left[(\boldsymbol{U}^i-\boldsymbol{U}^j)\cdot\boldsymbol{q}^{ij}\right]^2\left(\boldsymbol{\mathcal{C}}_{2}^{ij}\delta^{-3/2}+\boldsymbol{\mathcal{C}}_{3}^{ij}\delta^{-1/2}\right)\\[7pt]
+& \displaystyle \left.\beta_{ij}^2\left(\boldsymbol{\mathcal{C}}_{4}^{ij}\delta^{-5/2}+\boldsymbol{\mathcal{C}}_{5}^{ij}\delta^{-3/2}+\boldsymbol{\mathcal{C}}_{6}^{ij}\delta^{-1/2}\right)\right.\\[7pt]
+& \displaystyle R(\omega^i-\omega^j)\beta_{ij}\left(\boldsymbol{\mathcal{C}}_{7}^{ij}\delta^{-3/2}+\boldsymbol{\mathcal{C}}_{8}^{ij}\delta^{-1/2}\right)\\[7pt]
+& \displaystyle
R^2(\omega^i-\omega^j)^2\boldsymbol{\mathcal{C}}_{9}^{ij}\delta^{-1/2},
\quad \mbox{for } j \in \mathcal{N}_i \cap \mathbb{I},
\end{array}
\end{equation}
\begin{equation} \label{E:Q_ij-bound}
\begin{array}{r l l}
\displaystyle
Q_{ij}(\boldsymbol{U}^i,\omega^{i},\beta_{ij},\boldsymbol{f})
= & \displaystyle \beta_{ij}^2\left[\boldsymbol{\mathcal{B}}_{1}^{ij}\delta^{-5/2}+\boldsymbol{\mathcal{B}}_{2}^{ij}\delta^{-3/2}+\boldsymbol{\mathcal{B}}_{3}^{ij}\delta^{-1/2}\right] \\[7pt]
+ &  \displaystyle \left[(\boldsymbol{U}^i-\boldsymbol{U}^j)\cdot\boldsymbol{p}^{ij}+R\omega^i\right]^2\boldsymbol{\mathcal{B}}_{4}^{ij}\delta^{-1/2}\\[7pt]
+ &  \displaystyle R^2(\omega^i-\omega^j)^2\boldsymbol{\mathcal{B}}_{5}^{ij}\delta^{-1/2}\\[7pt]
+ &  \displaystyle \left[(\boldsymbol{U}^i-\boldsymbol{U}^j)\cdot\boldsymbol{q}^{ij}\right]^2\left(\boldsymbol{\mathcal{B}}_{6}^{ij}\delta^{-3/2}+\boldsymbol{\mathcal{B}}_{7}^{ij}\delta^{-1/2}\right)\\[7pt]
+ &  \displaystyle \beta_{ij} \left[(\boldsymbol{U}^i-\boldsymbol{U}^j)\cdot\boldsymbol{p}^{ij}+R\omega^i\right](\boldsymbol{\mathcal{B}}_{8}^{ij}\delta^{-3/2}+\boldsymbol{\mathcal{B}}_{9}^{ij}\delta^{-1/2})\\[7pt]
+ &  \displaystyle \beta_{ij} R(\omega^i-\omega^j)(\boldsymbol{\mathcal{B}}_{10}^{ij}\delta^{-3/2}+\boldsymbol{\mathcal{B}}_{11}^{ij}\delta^{-1/2})\\[7pt]
+ &  \displaystyle \beta_{ij} R\omega^i\boldsymbol{\mathcal{B}}_{12}^{ij}\delta^{-1/2}\\[7pt]
+ &  \displaystyle \left[(\boldsymbol{U}^i-\boldsymbol{U}^j)\cdot\boldsymbol{p}^{ij}+R\omega^i\right]R(\omega^i-\omega^j)\boldsymbol{\mathcal{B}}_{13}^{ij}\delta^{-1/2}\\[7pt]
+ & \displaystyle
\left[(\boldsymbol{U}^i-\boldsymbol{U}^j)\cdot\boldsymbol{q}^{ij}\right]Ra\boldsymbol{\mathcal{B}}_{14}^{ij}\delta^{-1/2},
\quad \mbox{for } j \in \mathcal{N}_i \cap \mathbb{B},
\end{array}
\end{equation}
called the \textit{discrete dissipation rates}, with coefficients
$\boldsymbol{\mathcal{C}}_{k}^{ij}$, $k=1,\ldots,9$,
$\boldsymbol{\mathcal{B}}_{m}^{ij}$, $m=1,\ldots,14$, which depend
on $\mu$, the ratio $\displaystyle \frac{R}{d_{ij}}$ and explicitly given by \eqref{E:coefficients} in Appendix \ref{A:coeffs}.

The solution of the discrete problem \eqref{E:I} is a set of discrete variables
$(\mathbb{\widehat{U}},\boldsymbol{\widehat{\omega}},\boldsymbol{\widehat{\beta}})\in
\mathcal{R}$, where $\mathbb{\widehat{U}}$ represents the set of
translational velocities of inclusions,
$\boldsymbol{\widehat{\omega}}$ the set of angular velocities and
$\boldsymbol{\widehat{\beta}}$ characterizes the Poiseuille
microflow in necks between neighboring inclusions.

\begin{remark} \label{R:pos-def}
$Q(\mathbb{U},\boldsymbol{\omega},\boldsymbol{\beta})$, defined by \eqref{E:Q}, is a positive definite quadratic form (see Appendix \ref{A:Q}).
\end{remark}

\begin{remark}  \label{R:coeff}
The agreement of the coefficients in \eqref{E:Q_ij} (explicitly given in \eqref{E:coefficients}) with
the previous results of \cite{bbp} is as follows. Only
coefficients $\boldsymbol{\mathcal{C}}_{1}^{ij}$,
$\boldsymbol{\mathcal{C}}_{2}^{ij}$ coincide with the
corresponding coefficients
$C^{ij}_{sh}$, $C^{ij}_{sp}$ in \cite{bbp}. This is because the coefficients in \cite{bbp} are
obtained by using the approximation of circular surfaces of
inclusions by parabolas whereas in this paper we use the true
circular surfaces. The main objective of \cite{bbp} was capturing
the strong blow up term of order $\delta^{-3/2}$ only, which
requires coefficients $C^{ij}_{sp}$
($\boldsymbol{\mathcal{C}}_{2}^{ij}$). The parabolic approximation
does not bring any discrepancy in $C^{ij}_{sp}$ whereas it may
bring a discrepancy in some other coefficients. Also since in
\cite{bbp} only the leading term was considered under special
boundary conditions there was no need to consider the Poiseuille
microflow.
\end{remark}

Both Theorems \ref{T:main-thm1} and \ref{T:main-thm4} are based on the following technical proposition.
\begin{proposition} %[\textit{Representation of the error term}]
\label{T:main-thm3}
Suppose $\Omega_F$ satisfies the close packing condition. Let the triple $(\mathbb{\widehat{U}},\boldsymbol{\widehat{\omega}},\boldsymbol{\widehat{\beta}}) \in \boldsymbol{\mathcal{R}}$ solve the discrete problem
\eqref{E:I}-\eqref{E:Q_ij-bound}. Then the following estimate holds:
\begin{equation} \label{E:main-thm3}
\begin{array}{r l l l }
|\widehat{W}-\mathcal{I}|\leq & \displaystyle
\mu\left(\sum_{i\in\mathbb{I}}\sum_{j\in\mathcal{N}_i}C_1\widehat{\beta}_{ij}^{2}+
C_2|\boldsymbol{\widehat{U}}^i-\boldsymbol{\widehat{U}}^j|^2+C_3\sum_{i\in\mathbb{I}\cup\mathbb{B}}R^2(\widehat{\omega}^i)^2\right).
\end{array}
\end{equation}
where $C_k$, $k=1,2,3$, are dimensionless constants.
\end{proposition}

The quadratic form $Q$ defined by \eqref{E:Q}-\eqref{E:Q_ij-bound} can be written in the following form:
\begin{equation} \label{E:form-decomp}
Q=Q^{in}_{sh}+Q^{in}_{sq}+Q^{in}_{per}+Q^{b}_{per}+Q^{b}_{sq},
\end{equation}
where
\[
Q^{in}_{sh}
= \sum_{i\in \mathbb{I}}\sum_{j\in\mathcal{N}_i\cap\mathbb{I}} \left[(\boldsymbol{U}^i-\boldsymbol{U}^j)\cdot\boldsymbol{p}^{ij}+ R\omega^i+R\omega^j\right]^2\boldsymbol{\mathcal{C}}_{1}^{ij}\delta^{-1/2},
\]
\[
Q^{in}_{sq}
= \sum_{i\in \mathbb{I}}\sum_{j\in\mathcal{N}_i\cap\mathbb{I}}\left[(\boldsymbol{U}^i-\boldsymbol{U}^j)\cdot\boldsymbol{q}^{ij}\right]^2\left(\boldsymbol{\mathcal{C}}_{2}^{ij}\delta^{-3/2}+\boldsymbol{\mathcal{C}}_{3}^{ij}\delta^{-1/2}\right),
\]
\[
\begin{array}{r l }
Q^{in}_{per}& \displaystyle = \sum_{i\in \mathbb{I}}\sum_{j\in\mathcal{N}_i\cap\mathbb{I}}\left\{\beta_{ij}^2\left(\boldsymbol{\mathcal{C}}_{4}^{ij}\delta^{-5/2}+\boldsymbol{\mathcal{C}}_{5}^{ij}\delta^{-3/2}+\boldsymbol{\mathcal{C}}_{6}^{ij}\delta^{-1/2}\right)\right.\\[10pt]
& \displaystyle + R(\omega^i-\omega^j)\beta_{ij}\left(\boldsymbol{\mathcal{C}}_{7}^{ij}\delta^{-3/2}+\boldsymbol{\mathcal{C}}_{8}^{ij}\delta^{-1/2}\right)\\[10pt]
& \displaystyle +\left.R^2(\omega^i-\omega^j)^2\boldsymbol{\mathcal{C}}_{9}^{ij}\delta^{-1/2}\right\}
\end{array}
\]
\[
\begin{array}{r l }
Q_{per}^{b} & \displaystyle = \sum_{i\in \mathbb{I}}\sum_{j\in\mathcal{N}_i\cap\mathbb{B}}
\left\{ \beta_{ij}^2\left[\boldsymbol{\mathcal{B}}_{1}^{ij}\delta^{-5/2}+\boldsymbol{\mathcal{B}}_{2}^{ij}\delta^{-3/2}+\boldsymbol{\mathcal{B}}_{3}^{ij}\delta^{-1/2}\right] \right.\\[7pt]
& +  \displaystyle \boldsymbol{\mathcal{B}}_{8}^{ij}\beta_{ij}\left[(\boldsymbol{U}^i-\boldsymbol{U}^j)\cdot\boldsymbol{p}^{ij}+R\omega^i\right]\delta^{-3/2}+\boldsymbol{\mathcal{B}}_{10}^{ij}\beta_{ij} R(\omega^i-\omega^j)\delta^{-3/2}\\[7pt]
& +  \displaystyle \boldsymbol{\mathcal{B}}_{9}^{ij}\beta_{ij}\left[(\boldsymbol{U}^i-\boldsymbol{U}^j)\cdot\boldsymbol{p}^{ij}+R\omega^i\right]\delta^{-1/2}+\boldsymbol{\mathcal{B}}_{11}^{ij}\beta_{ij} R(\omega^i-\omega^j)\delta^{-1/2}\\[7pt]
& +  \displaystyle \boldsymbol{\mathcal{B}}_{12}^{ij}\beta_{ij}R\omega^i\delta^{-1/2}+\boldsymbol{\mathcal{B}}_{4}^{ij}\left[(\boldsymbol{U}^i-\boldsymbol{U}^j)\cdot\boldsymbol{p}^{ij}+R\omega^i\right]^2\delta^{-1/2}\\[7pt]
& +  \displaystyle \boldsymbol{\mathcal{B}}_{13}^{ij}\left[(\boldsymbol{U}^i-\boldsymbol{U}^j)\cdot\boldsymbol{p}^{ij}+R\omega^i\right]R(\omega^i-\omega^j)\delta^{-1/2}\\[7pt]
& \left.+  \displaystyle \boldsymbol{\mathcal{B}}_{5}^{ij}R^2(\omega^i-\omega^j)^2\delta^{-1/2}\right\},
\end{array}
\]
\begin{equation} \label{E:various-Q}
\begin{array}{r l }
Q_{sq}^{b} & \displaystyle = \sum_{i\in \mathbb{I}}\sum_{j\in\mathcal{N}_i\cap\mathbb{B}} \left\{
\left[(\boldsymbol{U}^i-\boldsymbol{U}^j)\cdot\boldsymbol{q}^{ij}\right]^2\left(\boldsymbol{\mathcal{B}}_{6}^{ij}\delta^{-3/2}+\boldsymbol{\mathcal{B}}_{7}^{ij}\delta^{-1/2}\right)\right.\\[7pt]
& \displaystyle \left. + \boldsymbol{\mathcal{B}}_{14}^{ij}\left[(\boldsymbol{U}^i-\boldsymbol{U}^j)\cdot\boldsymbol{q}^{ij}\right]Ra\delta^{-1/2}\right\}.
\end{array}
\end{equation}

\begin{remark}
Both $Q_{per}^{in}$ and $Q_{per}^{b}$ are quadratic forms (see Appendix I).
\end{remark}

\begin{remark}  \label{R:Q}
The discrete dissipation rate $\mathcal{I}$ makes physics transparent. We can see how microflows enter this dissipation rate.
Indeed, the \textit{discrete dissipation form} $Q$ is presented as a sum of three motions corresponding to the decomposition \eqref{E:W-decomp}.
The decomposition \eqref{E:form-decomp} reflects the physics of the problem.
Namely, for the interior necks the first term $Q_{sh}^{in}$ corresponds to the dissipation rate $W_{\Pi_{ij}}(\boldsymbol{u}_{sh})$ in \eqref{E:W-decomp} and describes the shear flow
between inclusions. This type of flow is produced by two motions: rotation in the same direction (Fig. \ref{F:bc2}$a$) and relative
shear (Fig. \ref{F:bc1}$a$). The second term $Q_{sq}^{in}$ corresponds to the dissipation rate $W_{\Pi_{ij}}(\boldsymbol{u}_{sq})$ in \eqref{E:W-decomp} and describes the
local flow due to the squeeze motion (Fig. \ref{F:bc1}$b$). The third term $Q_{per}^{in}$ corresponds to the dissipation rate $W_{\Pi_{ij}}(\boldsymbol{u}_{per})$
in \eqref{E:W-decomp} and describes the ``permeation" type flow due to the Poiseuille microflow between motionless inclusions (Fig.
\ref{F:bc3}) and the rotation into the opposite directions (Fig. \ref{F:bc2}$b$). Finally, the last term $Q_b$ describes the local flow in the thin gap between a
disk and the external boundary. Here the decomposition into the parts corresponding to permeation $Q_{per}^{b}$ and squeeze $Q_{sq}^{b}$ in boundary necks is similar
to one in interior necks.
\end{remark}

\begin{remark}  \label{R:some}
The asymptotics \eqref{E:main-thm3} is given in terms of the \textit{discrete
dissipation rates} $Q_{ij}$. These dissipation rates are quadratic forms of the key physical parameters $\mathbb{U}, \boldsymbol{\omega}, \boldsymbol{\beta}$.
In particular, the overall discrete dissipation rate reveals the functional dependance of the overall properties of suspensions (dissipation rate,
effective viscosity, effective permeability) on the microflows (the shear, squeeze and permeation between neighbors).
\end{remark}

\Section{Anomalous Rate of Blow up of the Dissipation Rate, Qualitative Conclusions, Discussions and Open Problems} \label{S:discussions}

\begin{figure}[!ht]
  \centering
  \includegraphics[scale=1.05]{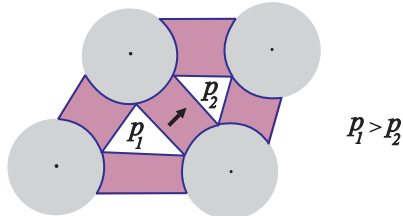}
  \caption{Local flux due to pressure drop}\label{F:poiseuille}
\end{figure}

In this section we show that in 3D the Poiseuille microflow does not occur (the fluid simply flows around the neck).
In contrast, in 2D incompressible fluid may need to permeate through necks. Necks separate fluid domain into
\textit{disconnected} regions, triangles, which may have different pressures $p_1,\,p_2$ (Fig. \ref{F:poiseuille}).
In this sense, the 2D problem becomes more complicated than an analogous 3D problem.
Moreover, this type of flow in experimental 2D settings may lead to an observable physical effect.

The techniques used in the previous study \cite{bbp} of 2D flow
were capable of capturing only the strong blow up term of order
$O(\delta^{-3/2})$ of the asymptotics of the effective viscosity.
The subsequent study of \cite{bp} reveals the significance of the
weak blow up of order $O(\delta^{-1/2})$ in 2D. It was shown that it
becomes the leading term in the asymptotics of the shear effective
viscosity. The objectives of \cite{bp} was to study the degeneracy
of the strong blow up term and evaluate the order of the magnitude
of the next term which was shown to exhibit the weak blow up. The
qualitative conclusion of the analysis of \cite{bp} was that the
shear viscosity exhibits the weak blow up in both 2D and 3D. While
this study highlighted the significance of the weak blow up term,
the techniques of \cite{bbp} did not allow to calculate it. As
shown in Theorem \ref{T:main-thm3} above the fictitious fluid
approach suggested in this paper allows to calculate
\textit{\textbf{any}} singular term, and, in particular, the weak
blow up. While a generalization of techniques of \cite{bbp} for 3D
suspensions leads to increasing technical difficulties, the
fictitious fluid approach provides an appropriate tool to attack
this problem. We anticipate that this approach would also be
useful in a variety of similar physical problems (e.g. rigid
inclusions in an elastic medium in both 2D and 3D). Moreover, the
analysis by this approach reveals the significance of the
Poiseuille microflows. In this section we present an example which
illustrates the following. For the suspensions of free inclusions
the Poiseuille microflows contribute to the weak blow up. However,
if an external field, which ``clamps'' inclusions, is imposed on inclusions then the Poiseuille
microflow may result in a new type singular behavior of viscous
dissipation rate (superstrong blow up). We also present an example of one disk in the fluid that can be clamped by the fluid flow with no external field. However, it
is not clear whether this example can be generalized to an ensemble of inclusions. Our example may suggest that on a suspension of free inclusions
the superstrong blow up may occur only due to boundary layer effects and therefore for the large number of inclusions becomes negligible.

\subsection{Suspensions in a Pinning Field} \label{SS:superstrong-blow_up}

\begin{example} \label{EX:disc-1}
The techniques described above can be applied to the problems of suspensions of non-neutrally buoyant rigid inclusions defined by
\begin{equation}   \label{E:ex_form}
\left\{
\begin{array}{r l l}
(a) & \displaystyle \mu\triangle \boldsymbol{u}=\nabla p, & \boldsymbol{x} \in \Omega_F \\[5pt]
(b) & \displaystyle \nabla\cdot \boldsymbol{u}=0, &  \boldsymbol{x} \in \Omega_F\\[5pt]
(c) & \displaystyle \boldsymbol{u}=\boldsymbol{U}^{i}+R\omega^i(n_1^i\boldsymbol{e}_2-n_2^i\boldsymbol{e}_1), & \boldsymbol{x} \in \partial B^{i}, \quad i=1\ldots N \\[5pt]
(d) & \displaystyle \int_{\partial B^{i}}\boldsymbol{\sigma}(\boldsymbol{u})\boldsymbol{n}^i ds=m\boldsymbol{g}, & i=1\ldots N \\[10pt]
(e) & \displaystyle \int_{\partial B^{i}}\boldsymbol{n}^i\times\boldsymbol{\sigma}(\boldsymbol{u})\boldsymbol{n}^ids=0, & i=1\ldots N \\[5pt]
(f) & \displaystyle \boldsymbol{u}=\boldsymbol{f}, &  \boldsymbol{x} \in \partial\Omega
\end{array}
\right.
\end{equation}
where the domain $\Omega_F$ is depicted in Fig. \ref{F:per_ex}. This problem corresponds to minimization of $\displaystyle W_{\Omega_F}(\boldsymbol{u})+ \sum_{i=1}^N
m\boldsymbol{g}\cdot \boldsymbol{U}^i$. The additional term does not change the analysis.

\begin{figure}[!ht]
  \centering
  \includegraphics[scale=1.05]{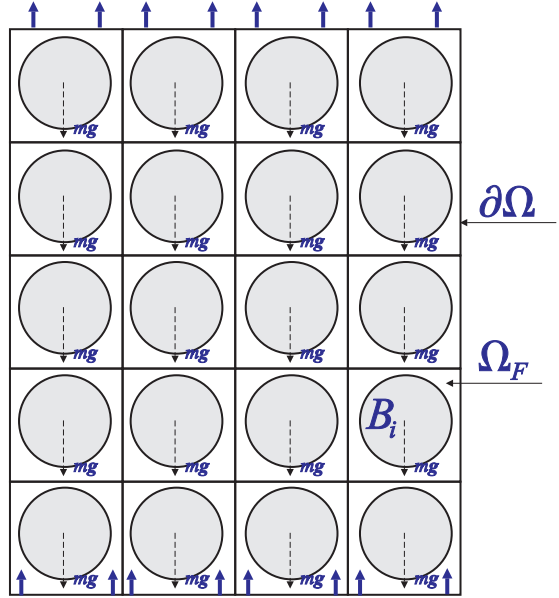}
  \caption{Example of a domain occupied by a suspension in a presence of a pinning field (gravity)}\label{F:per_ex}
\end{figure}

Here we suppose that the density of the solid inclusions is $\rho_s$ and
the fluid density is $\rho_f$. Then the force on the inclusion in the left hand side of
(\ref{E:ex_form}$d$) counteracts the external gravitational field
$-m\boldsymbol{g}$, where $\boldsymbol{g}=(0,g)$ is acceleration
due to the gravity and $m=\pi R^2(\rho_s-\rho_f)$ is the excess
mass. We choose the external boundary condition $\boldsymbol{f}$ and the force exerted by the heavy disks on the fluid so that the inclusions do not move and fluid is forced
to permeate through the thin gaps between motionless inclusions.
The force exerted by the disks is equal to their weight
$\pi R^2(\rho_s-\rho_f)g$, where we choose inclusions to be superheavy:
\begin{equation}   \label{E:rho_s}
\rho_s=C\delta^{-5/2}, \mbox{ for some } C=C(\mu,R)>0,
\end{equation}
and the applied boundary data $\boldsymbol{f}$ is chosen so that
\begin{equation}   \label{E:boundary-f}
\boldsymbol{f}=\begin{pmatrix} 0 \\ f_{2}(x,y) \end{pmatrix} \in H^{1/2}(\partial\Omega),
\end{equation}
 where  $f_2$  is some periodic function of  $(x,y) \in \partial\Omega$.
Such a boundary data $\boldsymbol{f}$ and the asymptotics of $\rho_s$ are selected so that inclusions do not move (the gravity balances the viscous force).
This balance can be found by solving an auxiliary problem similar to the one considered in \cite{has} (for details see Appendix \ref{A:hasimoto}).

Then the following proposition holds.

\begin{proposition}[\textit{Superstrong Blow Up due to a Pinning Field}]   \label{P:discussion}
Let $\widehat{W}$ be the overall viscous dissipation rate of the problem \eqref{E:ex_form}.
There exist $\rho_s$ and $\boldsymbol{f}$ of the form \eqref{E:rho_s} and \eqref{E:boundary-f}, respectively, such that the following
asymptotic representation holds:
\[
\begin{array}{l l}
\widehat{W} & \displaystyle = N \left(\mathcal{C}_1\delta^{-5/2}+\mathcal{C}_2\delta^{-3/2}+\mathcal{C}_3\delta^{-1/2}\right)+
O(1), \quad \mbox{as}\quad \delta\rightarrow 0,
\end{array}
\]
\[
\begin{array}{l l l l}
\mbox{\textit{where} } & \displaystyle
\mathcal{C}_1=\frac{9}{4}\pi\mu R^{5/2}, & \displaystyle
\mathcal{C}_2=\frac{99}{160}\pi\mu R^{3/2}, & \displaystyle
\mathcal{C}_2=\frac{29241}{17920}\pi\mu R^{1/2},
\end{array}
\]
and $N$ is the total number of inclusions.
\end{proposition}
The proof of this proposition is also given in Appendix \ref{A:hasimoto}).

\end{example}

\subsection{No Singularity of the Dissipation Rate due to the Poiseuille Microflow in 3D}
\label{SS:3D_poiseuille}

\begin{example} \label{EX:disc-3}

In order to explain why the the Poiseuille microflow does not contribute to
the singular behavior of the dissipation rate in 3D let us consider what happens with the fluid between two neighboring inclusions. Let $K=[-L,L]^3$ be a cube. The
parts of two neighboring inclusions are modeled by the hemispheres attached to the top and bottom faces of the cube as shown in Fig. \ref{F:3d-pois}.
Consider a 3D analog of the problem \eqref{E:v_form} with boundary condition
$\boldsymbol{f}$ to be a given constant vector $\boldsymbol{V}$ on
two opposite sides of the lateral boundary and zero vector on the rest of the
boundary (Fig. \ref{F:3d-pois}). Since the effective viscous dissipation rate $\widehat{W}$ is bounded from above by the dissipation rate
$W(\boldsymbol{w})$ for any test function $\boldsymbol{w}$ it suffices to find $\boldsymbol{w}$ such that $W(\boldsymbol{w})=O(1)$.
\begin{figure}[!ht]
  \centering
  \includegraphics[scale=.75]{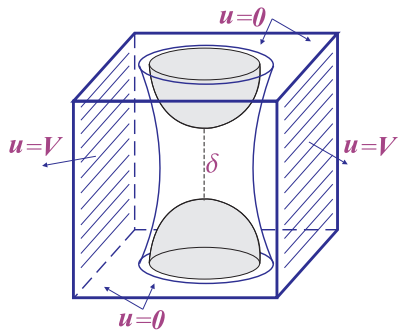}
  \caption{Poiseuille microflow in 3D}\label{F:3d-pois}
\end{figure}

Consider a ``hourglass domain'' $\Psi$ inside of the box containing two hemispheres as in Fig. \ref{F:3d-pois}.
We choose the trial function $\boldsymbol{w}$ to be:
\begin{equation} \label{E:trial-3dpois}
\boldsymbol{w}=
\begin{cases}
\boldsymbol{0}, & \mbox{when }
\boldsymbol{x} \in \Psi\\
\boldsymbol{W}, & \mbox{when } \boldsymbol{x} \in
K\setminus\Psi
\end{cases}
\end{equation}
where $\boldsymbol{W}$ solves the Stokes problem with Dirichlet boundary conditions: $\boldsymbol{W}=\boldsymbol{0}$ on the
boundary of $\Psi$ and $\boldsymbol{W}=\boldsymbol{f}$ on $\partial K$ (such a solution exists and is unique).

Evaluating the dissipation rate on such a trial function we
obtain:
\[
W(\boldsymbol{w})\leq \|\boldsymbol{w}\|^2_{H^{1}(K)}\leq C\|\boldsymbol{V}\|^2_{H^{1/2}(\partial K)}=O(1).
\]
\end{example}

\subsection{Boundary Layer Effects Leading to Superstrong Blow Up} \label{SS:boundary_layer}

\begin{example} \label{EX:boundary_layer}

Consider the domain $\Omega=(-1,1)^2$ which contains only one inclusion $B$.
Decompose the domain outside $B$ into necks $\displaystyle \bigcup_{i=1}^{4}\Pi_i$ and squares $\displaystyle \bigcup_{i=1}^{4}\Box_{i}$ as in Fig. \ref{F:boundary_layer}$(a)$.
Choose
\begin{equation*}
\boldsymbol{f}=
\begin{cases}
\begin{pmatrix} -1 \\ 1 \end{pmatrix}\zeta_{1}, & \partial \Box_1 \cap \partial \Omega=:\Gamma_1\\
\begin{pmatrix} -1 \\ -1 \end{pmatrix}\zeta_{2}, & \partial \Box_2 \cap \partial \Omega=:\Gamma_2\\
\begin{pmatrix} 1 \\ -1 \end{pmatrix}\zeta_{3}, & \partial \Box_3 \cap \partial \Omega=:\Gamma_3\\
\begin{pmatrix} 1 \\ 1 \end{pmatrix}\zeta_{4}, & \partial \Box_4 \cap \partial \Omega=:\Gamma_4\\
\boldsymbol{0}, & \mbox{elsewhere}
\end{cases}
\end{equation*}
where $\zeta_{i}$, $i=1,\ldots,4$ are smooth function having a compact support outside of $\Gamma_{i}$ such that $\zeta_i=1$ in $\Gamma_{i}$. Moreover, they satisfy the following symmetry condition:
\[
\zeta_{1}(x,y)=\zeta_{2}(-x,y)=\zeta_{3}(-x,-y)=\zeta_{4}(x,-y) \quad \mbox{and}\quad \zeta_{1}(x,y)=\zeta_{1}(y,x).
\]

\begin{figure}[!ht]
  \centering
  \includegraphics[scale=0.55]{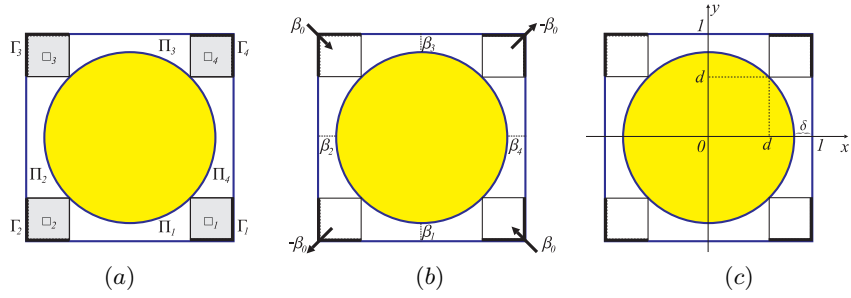}\\
  $(a)$ \hspace{3.5cm} $(b)$ \hspace{3.5cm} $(c)$
  \caption{One inclusion example: Boundary layer lead to superstrong blow up of viscous dissipation rate}\label{F:boundary_layer}
\end{figure}

The fluxes through the parts of the boundary $\Gamma_1$ and $\Gamma_3$ which are equal and we denote them as
\[
\beta_0=\frac{2}{R}(-1+d), \quad d=\delta+R(1-\frac{1}{\sqrt{2}}).
\]
Then the fluxes through $\Gamma_2$ and $\Gamma_4$ are $-\beta_0$.

Due to symmetry of the problem the inclusion does not rotate, that is, $\omega=0$ and $\beta_1=\beta_3=-\beta_2=-\beta_4$ (see Fig. \ref{F:boundary_layer}$(b)$). Then
\begin{equation*}
\begin{array}{r l}
\mathcal{I} & \displaystyle =\min_{\boldsymbol{U},\beta_i}\sum_{i=1}^{4}A\delta^{-3/2}[\boldsymbol{U}\cdot\boldsymbol{q}^i]^2+B\delta^{-1/2}[\boldsymbol{U}\cdot\boldsymbol{p}^i]^2+CR^2\delta^{-5/2}\beta_i^2\\[7pt]
& \displaystyle -p_1(\beta_0+\beta_1-\beta_4+U_1-U_2)-p_2(-\beta_0+\beta_2-\beta_1-U_1-U_2)\\[7pt]
& \displaystyle -p_3(\beta_0+\beta_3-\beta_2-U_1+U_2)-p_4(-\beta_0+\beta_4-\beta_3+U_1+U_2)\\[7pt]
& \displaystyle =\min_{\boldsymbol{U},\beta_i}\,\{2A\delta^{-3/2}[U_1^2+U_2^2]+2B\delta^{-1/2}[U_1^2+U_2^2]+4CR^2\delta^{-5/2}\beta_1^2\\[7pt]
& \displaystyle -p_1(2\beta_0+4\beta_1)-p_2(-2\beta_0-4\beta_1)\},
\end{array}
\end{equation*}
where $p_i$, $i=1,\ldots,4$ are the Lagrange multipliers corresponding to the weak incompressibility condition \eqref{E:disk-constraint}  and
$p_1=p_3$, $p_2=p_4$.

Solving the Euler-Lagrange equations for this minimization problem we obtain
\[
U_1=U_2=0, \quad \beta_1=-\frac{\beta_0}{2}=\frac{1}{R}(1-d),
\]
which provides the following asymptotics:
\[
\mathcal{I}=C(1-d)^2\delta^{-5/2}.
\]

Thus, we see that the superstrong blow up can occur due to the boundary layer where Poiseuille flow is significant.

\end{example}

\subsection{Free Suspensions (no external field)} \label{SS:free_suspensions}

In the above example \ref{EX:disc-1} we demonstrated that the Poiseuille microflow dominates the asymptotics of the overall
viscous dissipation rate and may cause the superstrong blow up of order $O(\delta^{-5/2})$. The key ingredient of this example is the presence of a strong external
pinning field which ``clamps" inclusions or alternatively the presence of the boundary layer as in Example \ref{EX:boundary_layer}.
In typical suspensions, with no external field, the inclusions are free to move. Then we expect that $\beta$ is asymptotically
small, as $\delta\rightarrow 0$, so that the Poiseuille microflow does not contribute to the superstrong blow up by Theorem
\ref{T:main-thm3}. This observation is also supported by the analysis of the periodicity
cell problem with five inclusions in \cite{bgmp}.
Below we present an example of a suspension with a hexagonal periodic array of inclusions
and prove that for the extensional external boundary conditions the viscous dissipation rate exhibits the strong blow up $O(\delta^{-3/2})$ since all $\beta_{ij}=0$.
Such a hexagonal array in 2D is a representative of a densely packed array of disks.

\begin{figure}[!ht]
  \centering
  \includegraphics[scale=1.05]{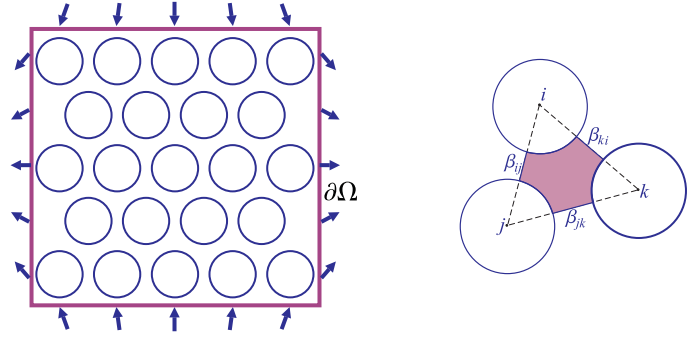}
  \caption{Periodic domain occupied by a suspension under the extensional boundary conditions that exhibits the strong blow up}\label{F:hex-ex}
\end{figure}
\begin{example} \label{EX:disc-4}

Consider a square domain $\Omega=(-M,M)^2$, $M>0$ with hexagonal array of disks $B^i$ centered at $(x_i,y_i)\in\Omega$. We consider
a periodic hexagonal array of such inclusions, as depicted on Figure \ref{F:hex-ex}).

The boundary value problem \eqref{E:v_form} is supplemented with boundary conditions
\begin{equation}   \label{E:bc-hex}
\boldsymbol{f}=\begin{pmatrix} x \\ -y \end{pmatrix} \quad \mbox{on} \quad\partial\Omega.
\end{equation}
We first prove that the overall dissipation rate $\widehat{W}_{hex}=O(\delta^{-3/2})$ because we construct a trial vector field $\boldsymbol{v}$ such that $W_{\Omega_F}(\boldsymbol{v})=O(\delta^{-3/2})$.
We thus obtain an upper bound for $\widehat{W}_{hex}$.

The construction is as follows. For each inclusion $B^i$, centered at $\boldsymbol{x}_i = (x_i, y_i)$ we prescribe the translational velocity to be exactly
\begin{equation}   \label{E:bc-ex}
\boldsymbol{U}^i  = \begin{pmatrix} x_i\\ -y_i \end{pmatrix},
\end{equation}
and the rotational velocity  $\omega^i=0$.
For such velocities the zero flux constraint \eqref{E:disk-constraint} for every fluid region $\mathcal{A}_{ijk}$ (see the right part of Fig. \ref{F:hex-ex}) takes the form
\begin{equation}   \label{E:beta-sum}
\beta_{ij}+\beta_{jk}+\beta_{ki}=0,
\end{equation}
because elementary computations show that for every $i,j,k$
\[
(\boldsymbol{U}^i+\boldsymbol{U}^j)\boldsymbol{p}^{ij}+(\boldsymbol{U}^j+\boldsymbol{U}^k)\boldsymbol{p}^{jk}+(\boldsymbol{U}^k+\boldsymbol{U}^i)\boldsymbol{p}^{ki}=0.
\]
Hence, we can choose all $\beta_{ij}$ to be identically zero. Then
Theorem \ref{T:main-thm3} implies that there exists a trial vector
field $\boldsymbol{v}$ such that
\[
W_{\Omega_F}(\boldsymbol{v})=O(\delta^{-3/2}),
\]
and, therefore, $\widehat{W}_{hex}=O(\delta^{-3/2})$.

To show that for this array of inclusions $\widehat{W}_{hex}>C\delta^{-3/2}$ we consider a chain of disks that connects upper and lower boundaries of the domain $\Omega$: $\partial\Omega^+=\{(x,y): y=M\}$ and $\partial\Omega^-\{(x,y): y=-M\}$, respectively
(see Fig. \ref{F:hex_ex_1}).
\begin{figure}[h!]
  \centering
  \includegraphics[scale=0.95]{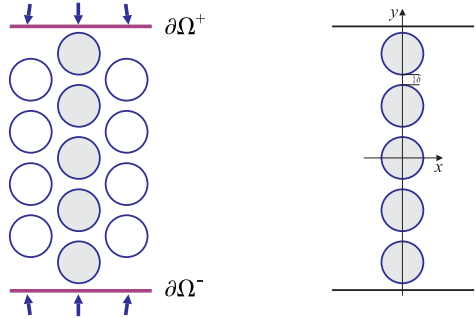}
  \caption{Chain of disks connecting the upper and lower boundaries of $\Omega$} \label{F:hex_ex_1}
\end{figure}
We choose this chain so that the $y$-axis of the coordinate system passes through the centers of disks in this chain.
Then
\[
\widehat{W}=\min_{\boldsymbol{\mathcal{R}}}Q(\mathbb{U},\boldsymbol{\omega},\boldsymbol{\beta})\geq (\mathcal{C}_2\delta^{-3/2}+\mathcal{C}_3\delta^{-1/2})\min_{\mathbb{U}}\sum_{chain}(U^i_2-U^j_2)^2.
\]
where the last minimum is taken over the disks in the chain (shadowed disks in Fig. \ref{F:hex_ex_1}).
Then
\[
\min_{\mathbb{U}}\sum_{chain}(U^i_2-U^j_2)^2=A>0,
\]
therefore $\widehat{W}_{hex}>C\delta^{-3/2}$.
This leads to a conclusion that the overall viscous dissipation rate $\widehat{W}_{per}$ exhibits the strong blow up of order $\delta^{-3/2}$.

\end{example}

It is known that for local movements of pairs of inclusions the squeeze type motion (Fig. \ref{F:bc1}$b$) provides the strongest
singularity (of order $\delta^{-3/2}$) whereas all other types of motions: rotations and shear provide a weaker singularity (of
order $\delta^{-1/2}$). Thus, it is natural to expect that in a suspension where inclusions are free to move with the extensional
boundary conditions the viscous dissipation exhibits the singularity of order $\delta^{-3/2}$. The above example shows that this
is indeed the case, that is, the superstrong blow up does not occur and the anomalous rate of $O(\delta^{-5/2})$ can be achieved if an external field is applied.

There is also one more case where the superstrong blow up can be obtained. Namely, Example
\ref{EX:boundary_layer} shows that it occurs in the boundary layers. Hence the above $O(\delta^{-3/2})$ conclusion applies to bulk effective properties
of free suspensions of large number of inclusions, when the boundary effects are negligible.

\Section{The Fictitious Fluid Problem} \label{SS:FFP}

The use of the fictitious fluid approach immediately gives a lower estimate on the viscous dissipation rate $\widehat{W}$ as follows.
\begin{lemma}   \label{L:main-ineq}
For $\widehat{W}_{\boldsymbol{\Pi}}$ defined by \eqref{E:W-fict} and $\widehat{W}$
defined by \eqref{E:W} the following inequality holds:
\begin{equation}   \label{E:main-ineq}
\widehat{W}_{\boldsymbol{\Pi}}\leq\widehat{W}.
\end{equation}
\end{lemma}

\textit{Proof.} The idea of the proof is to observe that for any incompressible vector field $\boldsymbol{u}$:
$W_{\boldsymbol{\Pi}}(\boldsymbol{u})\leq W_{\Omega_F}(\boldsymbol{u})$.
Therefore, consider the minimizer $\boldsymbol{u}$ of \eqref{E:v_form-var} and
\begin{equation}   \label{E:lemma1-proof}
\widehat{W}=\mu\int_{\Omega_F}\boldsymbol{D}(\boldsymbol{u}):\boldsymbol{D}(\boldsymbol{u})~d\boldsymbol{x}\geq
\mu\int_{\boldsymbol{\Pi}}\boldsymbol{D}(\boldsymbol{u}):\boldsymbol{D}(\boldsymbol{u})~d\boldsymbol{x}.
\end{equation}
since we disregard a part of the domain $\Omega_F$. Hence, $W_{\boldsymbol{\Pi}}(\boldsymbol{u})\leq\widehat{W}$.

We only need to show now that the restriction of the minimizer $\boldsymbol{u}$ on the set $\boldsymbol{\Pi}$ is an admissible trial field for minimization problem of the fictitious fluid \eqref{E:W-fict}, that is, $\boldsymbol{u}|_{\boldsymbol{\Pi}} \in V_{\boldsymbol{\Pi}}$. Indeed, $\boldsymbol{u}|_{\boldsymbol{\Pi}}$ belongs to $\boldsymbol{H}^{1}(\boldsymbol{\Pi})$, is divergence-free in
$\boldsymbol{\Pi}$ and satisfies
\[
0=\int_{\triangle_{ijk}}\nabla\cdot\boldsymbol{u}d\boldsymbol{x}=
\int_{\partial\triangle_{ijk}}\boldsymbol{u}\cdot\boldsymbol{n}ds.
\]
Thus,
\begin{equation*}
\min_{\boldsymbol{v}\in V_{\boldsymbol{\Pi}}}W_{\boldsymbol{\Pi}}(\boldsymbol{v})\leq W_{\boldsymbol{\Pi}}(\boldsymbol{u})\leq\widehat{W}= W_{\Omega_F}(\boldsymbol{u}).
\end{equation*}
$\Box$

Another significant advantage of the fictitious fluid problem is that the global minimization problem \eqref{E:W-fict}
can be split into two consecutive problems: one of them is on a single neck $\Pi_{ij}$, and the other one is a minimization problem on discrete variables $(\mathbb{U}, \boldsymbol{\omega},\boldsymbol{\beta}^*)\in\boldsymbol{\mathcal{R}}^*$.

\begin{lemma} \label{L:two-min} (Iterative minimization lemma).
Suppose $\widehat{W}_{\boldsymbol{\Pi}}$ is defined by \eqref{E:W-fict}. Then
\begin{equation} \label{E:min-chain_3}
\widehat{W}_{\boldsymbol{\Pi}}=\min_{(\mathbb{U},\boldsymbol{\omega},\boldsymbol{\beta}^*)\in \boldsymbol{\mathcal{R}}^*}\sum_{i\in\mathbb{I}, j\in \mathcal{N}_i}\min_{V_{ij}} W_{\Pi_{ij}}(\cdot).
\end{equation}
Moreover, the minimizer of $W_{\Pi_{ij}}$ over $V_{ij}$ satisfies the following Euler-Lagrange equations:
\begin{equation}   \label{E:EL-lower-beta}
\left\{
\begin{array}{r l l}
(a) & \displaystyle \mu\triangle \boldsymbol{u}=\nabla p, & \boldsymbol{x} \in \Pi_{ij}, \\
(b) & \displaystyle \nabla\cdot \boldsymbol{u}=0, & \boldsymbol{x} \in \Pi_{ij},\\
(c') & \displaystyle \boldsymbol{u}=\boldsymbol{U}^{i}+R\omega^i(n_1^i\boldsymbol{e}_2-n_2^i\boldsymbol{e}_1), & \boldsymbol{x} \in \partial B^{i}, \\
(c'') & \displaystyle \boldsymbol{u}=\boldsymbol{U}^{j}+R\omega^j(n_1^j\boldsymbol{e}_2-n_2^j\boldsymbol{e}_1), & \boldsymbol{x} \in \partial B^{j}, \\
(d) & \displaystyle \frac{1}{R}\int_{\ell_{ij}}\boldsymbol{u}\cdot\boldsymbol{n}ds=\beta_{ij}^*,\\
(e) & \displaystyle \boldsymbol{\sigma}(\boldsymbol{u})\boldsymbol{n}=-p^{\pm}_{ij}\boldsymbol{n}, & \boldsymbol{x} \in \partial\Pi^{\pm}_{ij},\\
(f) & \displaystyle \boldsymbol{u}=\boldsymbol{f}, & \boldsymbol{x} \in \partial\Pi_{ij}\cap \partial\Omega,
\end{array}
\right.
\end{equation}
where $p_{ij}^\pm$ are the Lagrange multipliers for the weak incompressibility condition \eqref{E:disk-constraint}.
\end{lemma}

\begin{figure}[!ht]
  \centering
  \includegraphics[scale=.55]{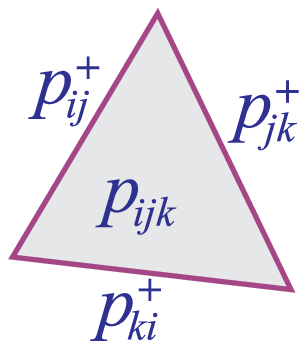} \hspace{1.5cm} \includegraphics[scale=.55]{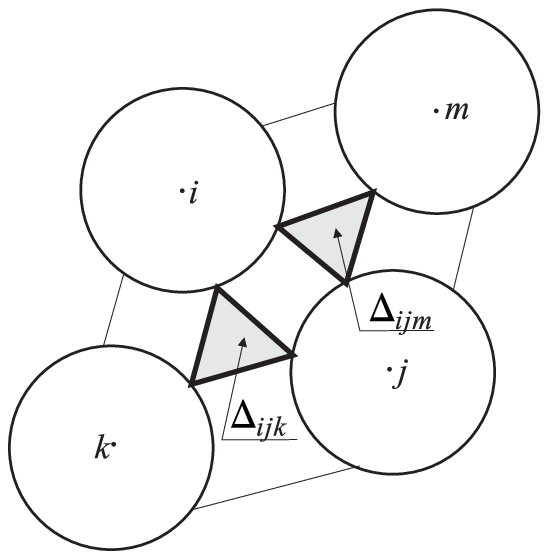}\\
  $(a)$ \hspace{3.5cm} $(b)$
  \caption{$(a)$ Pressures on the boundary of the triangle $\triangle_{ijk}$; $(b)$ Two triangles adjacent to the neck $\Pi_{ij}$}\label{F:pressures}
\end{figure}

\begin{proof}
Minimizing $W_{\boldsymbol{\Pi}}(\boldsymbol{u})$ over $V_{\boldsymbol{\Pi}}$ leads to the
Euler-Lagrange equations
\begin{equation}  \label{E:EL-boldneck}
\left\{
\begin{array}{r l l}
(a) & \displaystyle \mu\triangle \boldsymbol{u}=\nabla p, & \boldsymbol{x} \in \boldsymbol{\Pi}\\[5pt]
(b) & \displaystyle \nabla\cdot \boldsymbol{u}=0, &  \boldsymbol{x} \in \boldsymbol{\Pi}\\[5pt]
(c) & \displaystyle \boldsymbol{u}=\boldsymbol{U}^{i}+R\omega^i(n_1^i\boldsymbol{e}_2-n_2^i\boldsymbol{e}_1), & \boldsymbol{x} \in \partial B^{i}, \quad i=1\ldots N \\[5pt]
(d) & \displaystyle \int_{\partial B^{i}}\boldsymbol{\sigma}(\boldsymbol{u})\boldsymbol{n}^i ds=\boldsymbol{0} & i=1\ldots N \\[10pt]
(e) & \displaystyle \int_{\partial B^{i}}\boldsymbol{n}^i\times\boldsymbol{\sigma}(\boldsymbol{u})\boldsymbol{n}^ids=\boldsymbol{0}, & i=1\ldots N \\[7pt]
(f) & \displaystyle \int_{\partial\triangle_{ijk}}\boldsymbol{u}\cdot\boldsymbol{n}ds=0, &  i \in \mathbb{I}, \,\, j,k \in \mathcal{N}_i   \\[10pt]
(g) & \displaystyle \boldsymbol{\sigma}(\boldsymbol{u})\boldsymbol{n}=p_{ijk}\boldsymbol{n}, &  \boldsymbol{x} \in \partial  \triangle_{ijk}  \\[7pt]
(h) & \displaystyle \boldsymbol{u}=\boldsymbol{f}, & \boldsymbol{x} \in \partial\Omega,
\end{array}
\right.
\end{equation}
where the pressure $p(\boldsymbol{x})$ arises from the incompressibility condition  $\nabla\cdot \boldsymbol{u}=0$, $\boldsymbol{x} \in \boldsymbol{\Pi}$ and ``pressure constants" $p_{ijk}$  arise
from weak incompressibility condition
\[
\int_{\partial \triangle_{ijk}} \boldsymbol{v}\cdot\boldsymbol{n}ds=0.
\]
Given the boundary data  $f$ in \eqref{E:EL-boldneck}  we uniquely determine (see Appendix \ref{A:ex-un-2}) unknowns
\begin{equation}\label{set-1}
\boldsymbol{u},~p,~\boldsymbol{U}^{i},~\omega^i,p_{ijk},~i\in \mathbb{I},j,k\in \mathcal{N}_i.
\end{equation}
Fix a neck $\Pi_{ij}$ and consider the problem \eqref{E:EL-lower-beta} on it. For this pair of
indices $i,j$ take
\begin{equation}\label{set-2}
f, \boldsymbol{U}^{i}, \omega^i, \boldsymbol{U}^{j}, \omega^j, \beta_{ij}^*=\frac{1}{R}\int_{\ell_{ij}}\boldsymbol{u}\cdot\boldsymbol{n}ds
\end{equation}
 found from \eqref{set-1}. Using \eqref{set-2} as a given data, boundary value problem  \eqref{E:EL-lower-beta} can be solved uniquely (see Appendix \ref{A:ex-un-3}). Due to the unique solvability of
 both  \eqref{E:EL-boldneck} and \eqref{E:EL-lower-beta} the pair
 $(\boldsymbol{u},p)$ in \eqref{set-1} must solve  \eqref{E:EL-lower-beta} and
 \footnote{
 For notational convenience we identify
 $p_{ij}^+=p_{jk}^+=p_{ki}^+=p_{ijk}$, and $p_{ij}^+=p_{ji}^-$.
 }
 \[
 p_{ij}^+=p_{ijk},  p_{ij}^-=p_{ijm},
 \]
where triangles $\triangle_{ijk}$ and $\triangle_{ijm}$ are adjacent to the neck $\Pi_{ij}$
(see Fig.\ref{F:pressures}).
Hence  \eqref{E:EL-boldneck}  reduces to \eqref{E:EL-lower-beta}.

At this point the completion of the proof would be trivial if we did not have $\beta_{ij}^*$. Indeed,
in the entire domain $\Omega_F$ a result analogous to \eqref{E:min-chain_3} is simply
\[
\min_{(\mathbb{U},\boldsymbol{\omega}, \boldsymbol{u})}
W_{\Omega_F}(\boldsymbol{u})=\min_{(\mathbb{U},\boldsymbol{\omega})}\left(
\min_{ \boldsymbol{u}, \hbox{ when } (\mathbb{U},\boldsymbol{\omega}) \hbox{ fixed }}W_{\Omega_F}(\boldsymbol{u})\right).
\]
Hence it only remains to show that for any given  $(\mathbb{U},\boldsymbol{\omega})$ we can find
at least one set $\boldsymbol{\beta}^*$ satisfying the weak incompressibility condition \eqref{E:disk-constraint}. To this end we fix $(\mathbb{U},\boldsymbol{\omega})$ and  let $\boldsymbol{u}$ be the solution of the Stokes equation $\mu\triangle\boldsymbol{u}=\nabla p$ and $\nabla\cdot\boldsymbol{u}=0$ in the domain $\Omega_F$
with the Dirichlet data $\boldsymbol{u}|_{\partial B^{i}}=\boldsymbol{U}^{i}+R\omega^i(n_1^i\boldsymbol{e}_2-n_2^i\boldsymbol{e}_1)$,
$\boldsymbol{u}|_{\partial\Omega}=\boldsymbol{f}$. Set $\displaystyle \beta_{ij}^*=\frac{1}{R}\int_{\ell_{ij}}\boldsymbol{u}\cdot\boldsymbol{n}ds$.
Hence, we obtain $\boldsymbol{\beta}^*=\{\beta_{ij}^*\}$ such that $(\mathbb{U},\boldsymbol{\omega},\boldsymbol{\beta}^*)\in \boldsymbol{\mathcal{R}}^*$.
This completes the proof of lemma \ref{L:two-min}.
\end{proof}

\begin{remark}
For a given $(\mathbb{U},\boldsymbol{\omega})$ permeation constants $\boldsymbol{\beta}^*$ may not be unique. In fact, the choice of $\beta_{ij}^*$ has
$N$ degrees of freedom where $N$ is the number of inclusions. Indeed, $\boldsymbol{\beta}^*$ is found from solving a linear system $A\boldsymbol{\beta}^*=\boldsymbol{b}$
where the number of unknowns equals the number $P$ of interior necks and the number of equations equals the number of triangles, but there are only $P-N$ linearly independent
ones. Hence, the number of free parameters is equal to the number of inclusions.
\end{remark}

\begin{remark}
The unknowns of the problem \eqref{E:EL-boldneck} are the velocity field $\boldsymbol{u}(\boldsymbol{x})$, the pressure $p(\boldsymbol{x})$,
the constant translational $\boldsymbol{U}^i$ and angular $\omega^i$ velocities of the disk $B^i$, $i=1,\ldots,N$ and
constants $p_{ijk}$. Formally these constants appear as the Lagrange multipliers for the constraints \eqref{E:weak-incompres}.
Similar to how the pressure $p(\boldsymbol{x})$ appears as the Lagrange multiplier in the variational formulation corresponding to the Stokes equation, the
weak incompressibility condition for the fictitious fluid, inherited from the original fluid,
gives rise to the constant Lagrange multipliers, that one can regard as a constant pressure on the boundary of the fictitious fluid domain.
This also motivates the notations $p_{ijk}$. Thus, the fictitious fluid may be also interpreted as follows:
an incompressible fluid occupies necks while triangular pockets are filled with ``gas'' of the constant pressure $p_{ijk}$.
Naturally, the unknowns of the problem \eqref{E:EL-lower-beta} are the functions $\boldsymbol{u}(\boldsymbol{x})$, $p(\boldsymbol{x})$ and
constants $p_{ij}^{\pm}$, representing the velocity field, the pressure in the neck $\Pi_{ij}$ and the constant pressures on the lateral boundary $\partial\Pi_{ij}^{\pm}$, respectively.
\end{remark}

\begin{remark}
A major difficulty in applying previous `one-step' discretization techniques \cite{bbp,bn,bk} to
vectorial problems is the presence of integral constraints in the dual variational formulation.
In our `two-step' discretization approach, due to the Iterative Minimization Lemma, the inner minimization
problem has Dirichlet boundary conditions on inclusions and therefore neither this problem, nor its
dual have any integral conditions. On the other hand, due to \eqref{E:min-chain_3} second minimization
implies that these integral conditions are automatically satisfied.
\end{remark}

\begin{lemma} \label{L:upper}
Suppose $\Omega_F$ satisfies the close packing condition. Then for $\widehat{W}_{\boldsymbol{\Pi}}$ defined by \eqref{E:W-fict} and $\widehat{W}$
defined by \eqref{E:W} the following inequality holds:
\begin{equation}   \label{E:lemma-main}
\widehat{W}\leq\widehat{W}_{\boldsymbol{\Pi}} + \mu\left(\sum_{i\in\mathbb{I}}\sum_{j\in\mathcal{N}_i}C_1R^2 \widehat{\beta}_{ij}^{2}+
C_2|\boldsymbol{\widehat{U}}^i-\boldsymbol{\widehat{U}}^j|^2+C_3\sum_{i\in\mathbb{I}\cup\mathbb{B}}R^2(\widehat{\omega}^i)^2\right),
\end{equation}
where $(\mathbb{\widehat{U}},\boldsymbol{\widehat{\omega}},\boldsymbol{\widehat{\beta}^*})$ minimizes \eqref{E:min-chain_3}.
\end{lemma}
The proof of this lemma relies on the technical construction that appears in Subsection \ref{SS:proof_ff-approx} and therefore is put in Subsection \ref{SS:proof_upper}.

\begin{remark}
Lemma \ref{L:upper} is the only place where the close packing condition is necessary to obtain the desired estimate \eqref{E:lemma-main} because we needed uniform Lipschitz regularity of triangles $\triangle_{ijk}$.
\end{remark}

As a corollary of Lemmas \ref{L:main-ineq} and \ref{L:upper} we have the
main result of this section: the accuracy of approximation of the overall viscous dissipation rate by the dissipation rate of the fictitious fluid
given in the following proposition.
\begin{proposition}  \label{P:prop-main}
Suppose $\Omega_F$ satisfies the close packing condition. Then
\begin{equation}   \label{E:prop-main}
|\widehat{W}-\widehat{W}_{\boldsymbol{\Pi}}|\leq \mu\left(\sum_{i\in\mathbb{I}}\sum_{j\in\mathcal{N}_i}C_1 R^2 \widehat{\beta}_{ij}^{2}+
C_2|\boldsymbol{\widehat{U}}^i-\boldsymbol{\widehat{U}}^j|^2+C_3\sum_{i\in\mathbb{I}\cup\mathbb{B}}R^2(\widehat{\omega}^i)^2\right),
\end{equation}
where $(\mathbb{\widehat{U}},\boldsymbol{\widehat{\omega}},\boldsymbol{\widehat{\beta}^*})$ minimizes \eqref{E:min-chain_3}.
\end{proposition}

\Section{Discrete Network} \label{S:dn}

In the previous chapter we described a discrete network that arises from the fictitious fluid approach. Indeed, the equation \eqref{E:min-chain_3} in view of
\eqref{E:EL-lower-beta} is
\begin{equation}\label{uh-1}
\widehat{W}_{\boldsymbol{\Pi}}=\min_{(\mathbb{U},\boldsymbol{\omega},\boldsymbol{\beta}^*)\in \boldsymbol{\mathcal{R}}^*}\sum_{i\in\mathbb{I}, j\in \mathcal{N}_i} W_{\Pi_{ij}}(\boldsymbol{u})
=:\min_{(\mathbb{U},\boldsymbol{\omega},\boldsymbol{\beta}^*)\in \boldsymbol{\mathcal{R}}^*}\mathcal{W}(\mathbb{U},\boldsymbol{\omega},\boldsymbol{\beta}^*),
\end{equation}
where $\boldsymbol{u}$ is the solution of \eqref{E:EL-lower-beta} and $\mathcal{W}$ is a positive definite quadratic form of $(\mathbb{U},\boldsymbol{\omega},\boldsymbol{\beta}^*)$.

Our next objective is to find coefficients of $\mathcal{W}$ asymptotically as characteristic distance $\delta\rightarrow 0$. We have the following result.
\begin{lemma} \label{L:ff-approx}
Suppose $\widehat{W}_{\boldsymbol{\Pi}}$ is defined by \eqref{E:W-fict},
$\mathcal{I}$ is defined by \eqref{E:I},
 $\mathcal{W}$  is defined by \eqref{uh-1},
and $Q$ is defined by  \eqref{E:constraints-new}-\eqref{E:Q_ij-bound}.
Then
\begin{equation}   \label{E:ff-approx}
|\mathcal{W}(\mathbb{U},\boldsymbol{\omega},\boldsymbol{\beta}^*)-Q(\mathbb{U},\boldsymbol{\omega},\boldsymbol{\beta})| \leq
\mu\left(\sum_{i\in\mathbb{I}}\sum_{j\in\mathcal{N}_i}C_1 R^2 \beta_{ij}^{2}+
C_2|\boldsymbol{U}^i-\boldsymbol{U}^j|^2+C_3\sum_{i\in\mathbb{I}\cup\mathbb{B}}R^2(\omega^i)^2\right)
\end{equation}
for any $(\mathbb{U},\boldsymbol{\omega},\boldsymbol{\beta^*})\in \boldsymbol{\mathcal{R}}^*$
and $\beta_{ij}$ related to $\beta_{ij}^*$ through  \eqref{E:beta0-new}.
 In particular,
\begin{equation}   \label{E:ff-approx-min}
|\widehat{W}_{\boldsymbol{\Pi}}-\mathcal{I}| \leq
\mu\left(\sum_{i\in\mathbb{I}}\sum_{j\in\mathcal{N}_i}C_1R^2 \widehat{\beta}_{ij}^{2}+
C_2|\widehat{\boldsymbol{U}}^i-\widehat{\boldsymbol{U}}^j|^2+C_3\sum_{i\in\mathbb{I}\cup\mathbb{B}}R^2(\widehat{\omega}^i)^2\right).
\end{equation}
where $(\mathbb{\widehat{U}},\boldsymbol{\widehat{\omega}},\boldsymbol{\widehat{\beta}})$ is the minimizer of $Q$.
\end{lemma}

This Lemma shows that coefficients of $\mathcal{W}(\mathbb{U},\boldsymbol{\omega},\boldsymbol{\beta}^*)$ tend to infinity as $\delta\rightarrow 0$
because the corresponding coefficients of $Q(\mathbb{U},\boldsymbol{\omega},\boldsymbol{\beta})$ are asymptotically large and given in \eqref{E:Q_ij}-\eqref{E:Q_ij-bound}. The proof of this Lemma is given in Subsection \ref{SS:proof_ff-approx}.

Combining Proposition \ref{P:prop-main} and Lemma \ref{L:ff-approx}
we obtain the claim of Proposition \ref{T:main-thm3}.

In order to prove Theorems \ref{T:main-thm1} and \ref{T:main-thm4} it remains  to show that the error term of the right hand side of \eqref{E:ff-approx-min} becomes relatively small compared to the overall discrete dissipation
rate $\mathcal{I}$. In order to show that we prove in the next lemma that $\mathcal{I}\rightarrow\infty$ as $\delta\rightarrow 0$.
More specifically we have the following result.

\begin{lemma} \label{L:thm3-l2}
Suppose $\Omega_F$ satisfies the close packing condition. Then there exists a constant $C>0$ such that for every $\{\boldsymbol{U}^i\}$:
\begin{equation} \label{E:thm3-l2}
\begin{array}{l l }
\mathcal{I}& \displaystyle \geq \mu\sum_{i\in\mathbb{I}} \sum_{j\in \mathcal{N}_i}C_1\delta^{-3/2}|\boldsymbol{U}^i-\boldsymbol{U}^j|^2+
C_2\delta^{-1/2}R^2(\omega^i+\omega^j)^2\\[7pt]
& \displaystyle + C_3\delta^{-1/2}R^2(\omega^i-\omega^j)^2+C_4\delta^{-5/2}\beta_{ij}^2, \quad \mbox{as} \quad \delta\rightarrow 0.
\end{array}
\end{equation}
\end{lemma}

\begin{proof}
We write
\[
\begin{array}{l l }
\mathcal{I}& \displaystyle = \mu\sum_{i\in\mathbb{I}} \sum_{j\in \mathcal{N}_i} C_1\delta^{-3/2}[(\boldsymbol{U}^i-\boldsymbol{U}^j)\cdot\boldsymbol{q}^{ij}]^2 \\
& \displaystyle +C_2\delta^{-1/2}[(\boldsymbol{U}^i-\boldsymbol{U}^j)\cdot\boldsymbol{p}^{ij}+R\omega^i+R\omega^j]^2+C_3\delta^{-1/2}R^2(\omega^i-\omega^j)^2\\
& \displaystyle + C_4\delta^{-5/2}R^2\beta_{ij}^2.
\end{array}
\]
Using $(a+b)^2\geq 3/4a^2-3b^2$ we have
\[
\begin{array}{l l }
\mathcal{I}& \displaystyle \geq \mu\sum_{i\in\mathbb{I}} \sum_{j\in \mathcal{N}_i} C_1\delta^{-3/2}[(\boldsymbol{U}^i-\boldsymbol{U}^j)\cdot\boldsymbol{q}^{ij}]^2-3C_2\delta^{-1/2}|\boldsymbol{U}^i-\boldsymbol{U}^j|^2 \\
& \displaystyle +\frac{3}{4}C_2\delta^{-1/2}R^2(\omega^i+\omega^j)^2+C_3\delta^{-1/2}R^2(\omega^i-\omega^j)^2+ C_4\delta^{-5/2}R^2\beta_{ij}^2.
\end{array}
\]

By discrete Korn's inequality (Lemma \ref{L:korn}, Appendix \ref{A:korn}) when $\delta\rightarrow 0$:
\[
C_1\delta^{-3/2}[(\boldsymbol{U}^i-\boldsymbol{U}^j)\cdot\boldsymbol{q}^{ij}]^2-3C_2\delta^{-1/2}|\boldsymbol{U}^i-\boldsymbol{U}^j|^2
\geq C\delta^{-3/2}|\boldsymbol{U}^i-\boldsymbol{U}^j|^2.
\]
Hence, \eqref{E:thm3-l2} holds.

\end{proof}

From Lemma \ref{L:thm3-l2} and Proposition \ref{T:main-thm3} we have Theorem \ref{T:main-thm4}, and Theorem \ref{T:main-thm1} follows from Theorem \ref{T:main-thm4}
and Proposition \ref{P:prop-main}.

\Section{Proofs of Lemmas \ref{L:ff-approx} and \ref{L:upper}} \label{S:proof}

\subsection{Proof of Lemma \ref{L:ff-approx}} \label{SS:proof_ff-approx}

\begin{figure}[!ht]
  \centering
  \includegraphics[scale=1.05]{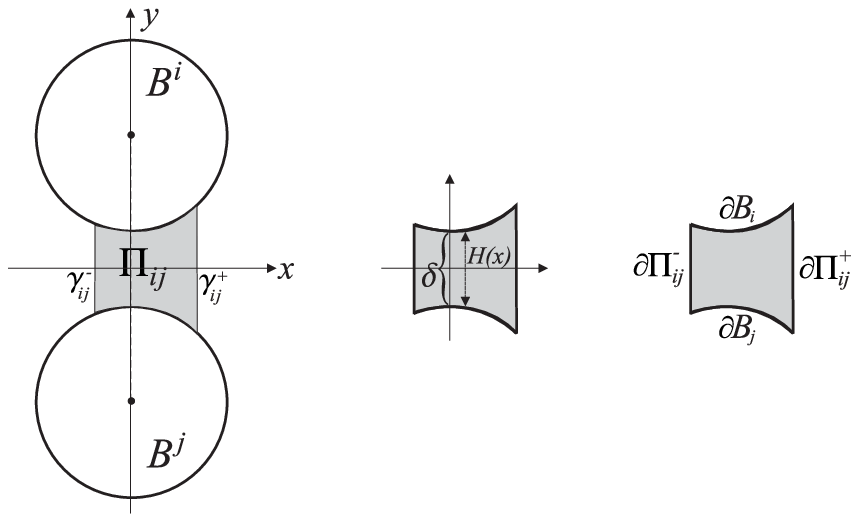}\\
  ($a$) \hspace{2.5cm}($b$)\hspace{2.5cm}($c$)
  \caption{($a$) The neck $\Pi_{ij}$ connecting two neighbors $B^i$ and $B^j$, ($b$) the distance $H(x)$ between two neighbors,
  ($c$) the boundary $\partial\Pi_{ij}$ of the neck}\label{F:neck00}
\end{figure}

{\bf Boundary value problem on a neck.} In the local coordinate system  the neck  $\Pi_{ij}$ (see Fig.\ref{F:neck00}) is given by
\begin{equation}   \label{E:neck_ij}
\Pi_{ij}=\left\{\boldsymbol{x}=(x,y) \in \mathbb{R}^2: \, -\gamma_{ij}^-<x < \gamma_{ij}^+, \,\,-H_{ij}(x)/2<y<H_{ij}(x)/2\right\},
\end{equation}
where $|\gamma_{ij}^+-\gamma_{ij}^-|$ is the neck width,
and $H_{ij}$ is the the distance between $\partial B_{i}$ and
$\partial B_{j}$ defined by
\begin{equation}   \label{E:distance}
H_{ij}(x)=
\displaystyle \delta_{ij}+2R - 2\sqrt{R^2-x^2}.
\end{equation}

Fix $(\mathbb{U}, \boldsymbol{\omega},\boldsymbol{\beta}^*) \in \boldsymbol{\mathcal{R}}^*$. In the chosen coordinate system the boundary conditions on $\partial B^i$ and $\partial B^j$ for functions defined in \eqref{E:various-u} are written as:
\begin{equation}  \label{E:G1}
\begin{array}{ r l l }
\boldsymbol{u}_{sh}|_{\partial B^i}=\boldsymbol{G}_{1}^{i}= & \displaystyle \frac{1}{2}(U^i_1-U^j_1)\boldsymbol{e}_{1}+R(\omega^i+\omega^j)\begin{pmatrix} \frac{1}{2}\sqrt{1-\frac{x^2}{R^2}} \\ \frac{x}{2R} \end{pmatrix}, \\
\boldsymbol{u}_{sh}|_{\partial B^j}=\boldsymbol{G}_{1}^{j}= & \displaystyle -\frac{1}{2}(U^i_1-U^j_1)\boldsymbol{e}_{1}+R(\omega^i+\omega^j)\begin{pmatrix} -\frac{1}{2}\sqrt{1-\frac{x^2}{R^2}} \\ \frac{x}{2R} \end{pmatrix},
\end{array}
\end{equation}
\begin{equation}  \label{E:G2}
\begin{array}{ r l l }
\boldsymbol{u}_{sq}|_{\partial B^i}=\boldsymbol{G}_{2}^{i}=  & \displaystyle  \frac{1}{2}(U^i_1-U^j_1)\boldsymbol{e}_{2},\\
\boldsymbol{u}_{sq}|_{\partial B^j}=\boldsymbol{G}_{2}^{j}= & \displaystyle \displaystyle -\frac{1}{2}(U^i_1-U^j_1)\boldsymbol{e}_{2},
\end{array}
\end{equation}
\begin{equation}  \label{E:G3}
\begin{array}{ r l l }
\boldsymbol{u}_{per}|_{\partial B^j}=\boldsymbol{G}_{3}^{i}= & \displaystyle R(\omega^i-\omega^j)\begin{pmatrix} \frac{1}{2}\sqrt{1-\frac{x^2}{R^2}} \\ \frac{x}{2R} \end{pmatrix}, \\
\boldsymbol{u}_{per}|_{\partial B^j}=\boldsymbol{G}_{3}^{j}= & \displaystyle R(\omega^i-\omega^j)\begin{pmatrix} \frac{1}{2}\sqrt{1-\frac{x^2}{R^2}} \\ -\frac{x}{2R} \end{pmatrix},
\end{array}
\end{equation}
and $\boldsymbol{u}_{sh}$, $\boldsymbol{u}_{sq}$, $\boldsymbol{u}_{per}$ are minimizers of $W_{\Pi_{ij}}$ over the following sets, respectively:
\begin{equation} \label{E:V-sh}
\begin{array}{l l}
V_{sh}=\left\{\boldsymbol{v} \in \boldsymbol{H}^1(\Pi_{ij}): \nabla\cdot\boldsymbol{v}=0\,\mbox{ in }\,\Pi_{ij}, \,\,\boldsymbol{v}= \boldsymbol{G}_{1}^{i}\,\mbox{ on }\,\partial B^i,\,\,\boldsymbol{v}=\boldsymbol{G}_{1}^{j}\,\mbox{ on }\,\partial B^j \right\}\\[7pt]
\end{array}
\end{equation}
\begin{equation} \label{E:V-sq}
\begin{array}{l l}
V_{sq}=\left\{\boldsymbol{v} \in \boldsymbol{H}^1(\Pi_{ij}): \nabla\cdot\boldsymbol{v}=0\,\mbox{ in }\,\Pi_{ij}, \,\, \boldsymbol{v}=\boldsymbol{G}_{2}^{i}\,\mbox{ on }\,\partial B^i,\,\, \,\, \boldsymbol{v}=\boldsymbol{G}_{2}^{i}\,\mbox{ on }\,\partial B^j \right\}\\[5pt]
\end{array}
\end{equation}
\begin{equation} \label{E:V-per}
\begin{array}{l l l}
V_{per}=
\left\{\boldsymbol{v} \in \boldsymbol{H}^1(\Pi_{ij}):\right.& \displaystyle\nabla\cdot\boldsymbol{v}=0\,\mbox{ in }\,\Pi_{ij},\,\, , \boldsymbol{v}=\boldsymbol{G}_{3}^{i}\,\mbox{ on }\,\partial B^i,\,\,\boldsymbol{v}=\boldsymbol{G}_{3}^{j}\,\mbox{ on }\,\partial B^j \\
& \displaystyle \left. \,\, \frac{1}{R}\int_{\ell_{ij}}\boldsymbol{v}\cdot\boldsymbol{n}ds=\beta_{ij}^*-\frac{\delta_{ij}}{2R}\left(U^i_1+U^j_1\right) \right\},
\end{array}
\end{equation}

{\bf Variational Duality.} For any $(\mathbb{U}, \boldsymbol{\omega},\boldsymbol{\beta}^*) \in \boldsymbol{\mathcal{R}}^*$ and any $\Pi_{ij} \in \boldsymbol{\Pi}$ we let
\begin{equation} \label{E:W-curl}
\mathcal{W}_{ij}(\boldsymbol{U}^i,\boldsymbol{U}^j,\omega^{i},\omega^{j},\beta_{ij}^*):=\min_{V_{ij}}W_{\Pi_{ij}}(\cdot).
\end{equation}
By orthogonality \eqref{E:W-decomp}
\[
\mathcal{W}_{ij}(\boldsymbol{U}^i,\boldsymbol{U}^j,\omega^{i},\omega^{j},\beta_{ij}^*)=W_{\Pi_{ij}}(\boldsymbol{u}_{sh})+W_{\Pi_{ij}}(\boldsymbol{u}_{sq})+W_{\Pi_{ij}}(\boldsymbol{u}_{per}).
\]
Using variational duality (see e.g. \cite{bk,bgn})\footnote{Since the construction is done in the neck $\Pi_{ij}$ hereafter we drop such a subscript in the functional $W_{\Pi_{ij}}=W$.
We also drop the indices $i,j$ in $H_{ij}(x)=H(x)$.}
\begin{equation} \label{E:max-sh}
W(\boldsymbol{u}_{sh})=\max_{\boldsymbol{\mathcal{S}}\in F }W^*_1(\boldsymbol{\mathcal{S}}),
\end{equation}
\begin{equation} \label{E:max-sq}
W(\boldsymbol{u}_{sq})=\max_{\boldsymbol{\mathcal{S}}\in F }W^*_2(\boldsymbol{\mathcal{S}}),
\end{equation}
\begin{equation} \label{E:max-per}
W(\boldsymbol{u}_{per})=\max_{\boldsymbol{\mathcal{S}}\in F_{per} }W^*_3(\boldsymbol{\mathcal{S}}),
\end{equation}
where the dual functionals are defined by ($k=1,2$)
\begin{equation} \label{E:dual-fnl-1}
\begin{array}{l l}
W^*_k(\boldsymbol{\mathcal{S}})& \displaystyle
=\int_{\partial B^i}\boldsymbol{G}^i_k\cdot \boldsymbol{\mathcal{S}}\boldsymbol{n}ds+ \int_{\partial B^j}\boldsymbol{G}^j_k\cdot \boldsymbol{\mathcal{S}}\boldsymbol{n}ds
-\frac{1}{4\mu}\int_{\Pi_{ij}}\left[\boldsymbol{\mathcal{S}}:\boldsymbol{\mathcal{S}}-\frac{(\mbox{tr}\,\boldsymbol{\mathcal{S}})^2}{2}\right]d\boldsymbol{x},
\end{array}
\end{equation}
\begin{equation} \label{E:dual-fnl-2}
\begin{array}{l l}
W^*_3(\boldsymbol{\mathcal{S}})& \displaystyle
=\int_{\partial B^i}\boldsymbol{G}^i_3\cdot \boldsymbol{\mathcal{S}}\boldsymbol{n}ds+ \int_{\partial B^j}\boldsymbol{G}^j_3\cdot \boldsymbol{\mathcal{S}}\boldsymbol{n}ds \\[10pt]
& \displaystyle +\frac{R\beta_{ij}^*}{H(\gamma_{ij}^+)}\int_{\partial \Pi_{ij}^+}\boldsymbol{n}\cdot \boldsymbol{\mathcal{S}}\boldsymbol{n}ds-
\frac{R\beta_{ij}^*}{H(\gamma_{ij}^-)}\int_{\partial \Pi_{ij}^-}\boldsymbol{n}\cdot \boldsymbol{\mathcal{S}}\boldsymbol{n}ds  \\[10pt]
& \displaystyle -\frac{1}{4\mu}\int_{\Pi_{ij}}
\left[\boldsymbol{\mathcal{S}}:\boldsymbol{\mathcal{S}}-
\frac{(\mbox{tr}\,\boldsymbol{\mathcal{S}})^2}{2}\right]
d\boldsymbol{x}
\end{array}
\end{equation}
over the sets
\begin{equation} \label{E:dual-set-1}
F=\left\{\boldsymbol{S} \in \mathbb{R}^{2\times 2}, \boldsymbol{S}=\boldsymbol{S}^{T},\,\,
\boldsymbol{S}_{ij} \in \boldsymbol{L}^{2}(\Pi_{ij}):\,\nabla\cdot\boldsymbol{S}=\boldsymbol{0} \mbox{ in } \Pi_{ij},
\boldsymbol{S}\boldsymbol{n}=\boldsymbol{0} \mbox{ on } \partial \Pi_{ij}^\pm \right\}.
\end{equation}
\begin{equation} \label{E:dual-set-2}
\begin{array}{l l}
F_{per}=\left\{\boldsymbol{S} \in \mathbb{R}^{2\times 2},\,\,\boldsymbol{S}=\boldsymbol{S}^{T},\,\,\boldsymbol{S}_{ij} \in \boldsymbol{L}^{2}(\Pi_{ij}):\right. & \displaystyle \,\bigtriangledown\cdot\boldsymbol{S}=\boldsymbol{0} \mbox{ in }\Pi_{ij},\\[5pt]
& \displaystyle \left.\,\,\boldsymbol{S}\boldsymbol{n}=\chi^\pm_{ij}\boldsymbol{n}\mbox{ on }\partial \Pi_{ij}^\pm\right\},
\end{array}
\end{equation}
and $\chi^+_{ij}$ and $\chi^-_{ij}$ are arbitrary constants. In \eqref{E:dual-fnl-2} the numbers $H(\gamma_{ij}^\pm)$ are the lengths of the lateral boundaries, that is, $H(\gamma_{ij}^+)=|\partial\Pi_{ij}^+|$,
$H(\gamma_{ij}^-)=|\partial\Pi_{ij}^-|$. It is straightforward to check that the maximizer of problem \eqref{E:max-sh}, \eqref{E:dual-fnl-1}, \eqref{E:dual-set-1} is the stress tensor
$\boldsymbol{\sigma}(\boldsymbol{u}_{sh})$, the maximizer of \eqref{E:max-sq}, \eqref{E:dual-fnl-1}, \eqref{E:dual-set-1} is $\boldsymbol{\sigma}(\boldsymbol{u}_{sq})$ (see e.g. \cite{bbp}) and
the maximizer of \eqref{E:max-per}, \eqref{E:dual-fnl-2}, \eqref{E:dual-set-2} is $\boldsymbol{\sigma}(\boldsymbol{u}_{per})$  (see Appendix \ref{A:dual-fnl}).

{\bf Proof.} By direct and dual variational principles
for any $\boldsymbol{\mathcal{S}}_1 \in F$, $\boldsymbol{\mathcal{S}}_2 \in F$, $\boldsymbol{\mathcal{S}}_3 \in F_{per}$ and
$\boldsymbol{v}_1 \in V_{sh}$, $\boldsymbol{v}_2 \in V_{sq}$, $\boldsymbol{v}_3 \in V_{per}$
we have the following bounds:
\begin{equation}   \label{E:bound}
W^*_1(\boldsymbol{\mathcal{S}}_1)+W^*_2(\boldsymbol{\mathcal{S}}_2)+W^*_3(\boldsymbol{\mathcal{S}}_3)\leq
W(\boldsymbol{u})\leq W(\boldsymbol{v}_1)+W(\boldsymbol{v}_2)+W(\boldsymbol{v}_3).
\end{equation}
Therefore, if we find trial fields $\boldsymbol{\mathcal{S}}_i$ and $\boldsymbol{v}_i$ such
that the following three inequalities hold (with universal constants $\mathcal{C}_i$, $i=1,...,5$)
\begin{equation}   \label{E:error-in-1}
\begin{array}{l l}
\left|W(\boldsymbol{v}_1)-W^*_1(\boldsymbol{\mathcal{S}}_{1})\right| & \displaystyle \leq\mu \left(\mathcal{C}_1[(\boldsymbol{U}^i-\boldsymbol{U}^j)\cdot\boldsymbol{p}^{ij}+R\omega^i+R\omega^j]^2
+ \mathcal{C}_2R^2(\omega^i+\omega^j)^2\right)\\[5pt]
& \displaystyle  =:E_1(\boldsymbol{U}^i,\boldsymbol{U}^j,\omega^i,\omega^j),
\end{array}
\end{equation}
\begin{equation}   \label{E:error-in-2}
\left|W(\boldsymbol{v}_{2})-W^*_2(\boldsymbol{\mathcal{S}}_{2})\right|\leq \mu \mathcal{C}_3[(\boldsymbol{U}^i-\boldsymbol{U}^j)\cdot\boldsymbol{q}^{ij}]^2=: E_2(\boldsymbol{U}^i,\boldsymbol{U}^j),
\end{equation}
\begin{equation}  \label{E:error-in-3}
\left|W(\boldsymbol{v}_{3})-W^*_3(\boldsymbol{\mathcal{S}}_{3})\right|\leq \mu \left(\mathcal{C}_4R^2(\omega^i-\omega^j)^2+\mathcal{C}_5R^2\beta_{ij}^2\right)=: E_3(\omega^i,\omega^j,\beta_{ij}),
\end{equation}
then for $\boldsymbol{v}=\boldsymbol{u}_t+\boldsymbol{v}_1+\boldsymbol{v}_2+\boldsymbol{v}_3$
we have
\[
\left|W(\boldsymbol{u})-W(\boldsymbol{v})\right|\leq
\left|W(\boldsymbol{v}_1)-W^*_1(\boldsymbol{\mathcal{S}}_{1})\right|+
\left|W(\boldsymbol{v}_2)-W^*_1(\boldsymbol{\mathcal{S}}_{2})\right|+
\left|W(\boldsymbol{v}_3)-W^*_1(\boldsymbol{\mathcal{S}}_{3})\right|
\]
\[
\leq \mu C_1\beta_{ij}^{2}+C_2|\boldsymbol{U}^i-\boldsymbol{U}^j|^2+C_3R^2(\omega^i)^2+C_4R^2(\omega^i)^2,
\]
and  Lemma \ref{L:ff-approx} will follow from summation of the above inequality over all necks and
explicit evaluation of $W(\boldsymbol{v}_1)$, $W(\boldsymbol{v}_2)$ and $W(\boldsymbol{v}_3)$.

The construction of $\boldsymbol{\mathcal{S}}_i$ and $\boldsymbol{v}_i$  for the interior necks
is done, for simplicity, for symmetric necks ($\gamma_{ij}^+=\gamma_{ij}^-=:\gamma_{ij}$)
and it is as follows.
\begin{equation}  \label{E:v1}
\begin{array}{l l}
\boldsymbol{v}_1 & \displaystyle = (U^i_1-U^j_1+R\omega^i+R\omega^j)\begin{pmatrix} \displaystyle yG(x)\\ \displaystyle F(x)-\frac{y^2}{2}G'(x) \end{pmatrix}\\[12pt]
& \displaystyle +R(\omega^i+R\omega^j)\begin{pmatrix} \displaystyle yK(x)\\ \displaystyle M(x)-\frac{y^2}{2}K'(x) \end{pmatrix} \quad \mbox{in} \quad \Pi_{ij},
\end{array}
\end{equation}
where functions $F$, $G$, $K$ and $M$ are found so that $\boldsymbol{v}_1$ satisfies \eqref{E:G1} on $\partial B^i, \partial B^j$, that is,
\begin{equation}  \label{E:F,G,K,M-1}
\begin{array}{l l l l}
& \displaystyle G(x)=\frac{1}{H(x)}, & & \displaystyle  F(x)=-\frac{H'(x)}{8}, \\[5pt]
& \displaystyle K(x)=\frac{\sqrt{1-\frac{x^2}{R^2}}-1}{H(x)}, & & \displaystyle M(x)=-\frac{x\delta}{8R\sqrt{R^2-x^2}}.
\end{array}
\end{equation}
\begin{equation}   \label{E:S1(1)}
\boldsymbol{\mathcal{S}}_{1} =
\mu \begin{pmatrix} \displaystyle 0 & \displaystyle G(x)-C \\ \displaystyle G(x)-C &  \displaystyle -yG'(x) \end{pmatrix} \quad \mbox{in} \quad \Pi_{ij},
\end{equation}
where the constant $C$ is chosen so that $\boldsymbol{\mathcal{S}}_{1}\in F$, namely $\displaystyle C=\frac{1}{H(\gamma_{ij})}$.
\begin{equation}  \label{E:v2}
\boldsymbol{v}_2(x,y)=(U^i_2-U^j_2)\begin{pmatrix} \displaystyle G(x)+3y^2F(x) \\
\displaystyle -yG'(x)-y^3F'(x) \end{pmatrix} \quad \mbox{in} \quad \Pi_{ij},
\end{equation}
where $G(x)$ and $F(x)$ are chosen so that $\boldsymbol{v}_2$ to satisfy \eqref{E:G2} on $\partial B^{i}, \, \partial B^{j}$, that is,
\begin{equation}  \label{E:G&F-v2}
G(x)=-\frac{3x}{2H(x)}, \quad \quad F(x)=\frac{2x}{H^3(x)}.
\end{equation}
\begin{equation}   \label{E:S2(1)}
\begin{array}{r l l }
& \displaystyle (\boldsymbol{\mathcal{S}}_{2})_{11}  = \mu (U^i_2-U^j_2)\left(G'+9y^2F'-6\int_{\gamma_0}^x F-C_1\frac{x^2}{2}-\frac{3}{2}C_2x^2y^2+C_3+C_4y^2\right)\\[5pt]
& \displaystyle (\boldsymbol{\mathcal{S}}_{2})_{12}  = (\boldsymbol{\mathcal{S}}_{2})_{21} = \mu(U^i_2-U^j_2)\left(\displaystyle 6yF-yG''-3y^3F''+C_1xy+C_2xy^3\right)\\[5pt]
& \displaystyle (\boldsymbol{\mathcal{S}}_{2})_{22}  = \mu(U^i_2-U^j_2)\left(-3G'-3y^2F'-6\int_{\gamma_0}^x F+\frac{y^2}{2}G'''+\frac{3}{4}y^4F'''-C_1\frac{y^2}{2}-C_2\frac{y^4}{4}\right)
\end{array}
\end{equation}
where constants $C_i$ ($i=1..4$) are chosen so that $\boldsymbol{\mathcal{S}}_{2}\boldsymbol{n}=\boldsymbol{0}$ on $\partial \Pi_{ij}^\pm$.
\begin{equation}  \label{E:v3}
\begin{array}{l l}
\boldsymbol{v}_3 & \displaystyle = R(\omega^i-\omega^j)\begin{pmatrix} \displaystyle P(x)+3y^2Q(x) \\ \displaystyle -yP'(x)-y^3Q'(x) \end{pmatrix}
+ \beta_{ij}R\begin{pmatrix} \displaystyle K(x)+3y^2M(x)\\ \displaystyle -yK(x)-y^3M'(x) \end{pmatrix},
\end{array}
\end{equation}
where functions $P$, $Q$, $K$ and $M$ are found so that $\boldsymbol{v}_3$ satisfies \eqref{E:G3} on $\partial B^i, \partial B^j$, that is,
\begin{equation}  \label{E:P,Q,K,M-3}
\begin{array}{l l l l}
& \displaystyle P(x)=-\frac{3x^2}{4RH(x)}+\frac{3H(x)}{16R}+\frac{1}{2}\sqrt{1-\frac{x^2}{R^2}}, & & \displaystyle  Q(x)=\frac{x^2}{RH^3(x)}-\frac{1}{4RH(x)}, \\[10pt]
& \displaystyle K(x)=\frac{3}{2H(x)}, & & \displaystyle M(x)=-\frac{2}{H^3(x)}.
\end{array}
\end{equation}
\begin{equation}   \label{E:S3(1)}
\begin{array}{r l}
(\boldsymbol{\mathcal{S}}_{3})_{11} = & \displaystyle \mu\left(P'+9y^2Q'-6\int Q -C_1\frac{x^2}{2}-\frac{3}{2}C_2x^2y^2+C_3+C_4y^2\right) ,\\[5pt]
(\boldsymbol{\mathcal{S}}_{3})_{12} = & \displaystyle (\boldsymbol{\mathcal{S}}_{3})_{21} =\mu\left(6yQ-yP''-3y^3Q''+C_1xy+C_2xy^3\right),\\[5pt]
(\boldsymbol{\mathcal{S}}_{3})_{22} = & \displaystyle\mu\left(-3P'-3y^2Q'-6\int Q +\frac{y^2}{2}P'''+\frac{3}{4}y^4Q'''-C_1\frac{y^2}{2}-C_2\frac{y^4}{4}\right),
\end{array}
\end{equation}
where
\begin{equation}  \label{E:G,F-3}
G(x)=R(\omega^i-\omega^j)P(x)+\beta_{ij}RK(x), \quad  F(x)=R(\omega^i-\omega^j)Q(x)+\beta_{ij}RM(x),
\end{equation}
and constants $C_i$ ($i=1..4$) are chosen so that $\boldsymbol{\mathcal{S}}_{3}\boldsymbol{n}=\chi^\pm\boldsymbol{n}$ on $\partial \Pi_{ij}^\pm$ for some constants $\chi^\pm$.

In Appendix \ref{A:constr-in} we verify that the chosen trial fields indeed satisfy inequalities  \eqref{E:error-in-1}, \eqref{E:error-in-2}, \eqref{E:error-in-3}, and also show that
\begin{equation} \label{E:W(1)}
\begin{array}{l l}
W(\boldsymbol{v}_1) & = \displaystyle \frac{1}{2}\pi\mu R^{1/2}(U_1^i-U_1^j+R\omega^i+R\omega^j)^2\delta_{ij}^{-1/2}\\[5pt]
& + \displaystyle \mathbf{C}_1\mu(U_1^i-U_1^j+R\omega^i+R\omega^j)^2+\mathbf{C}_2\mu R^2(\omega^i+\omega^j)^2\\[5pt]
& + \displaystyle \mathbf{C}_3\mu(U_1^i-U_1^j+R\omega^i+R\omega^j)R^2(\omega^i+\omega^j),
\end{array}
\end{equation}
\begin{equation}  \label{E:W(2)}
\begin{array}{l l l}
W(\boldsymbol{v}_2) & \displaystyle =(U^i_2-U^j_2)^2\left[\frac{3}{8}\pi\mu R^{3/2}\delta_{ij}^{-3/2}+\frac{207}{320}\pi\mu R^{-1/2}\delta_{ij}^{-1/2}\right]+\mathbf{C}_4\mu(U^i_2-U^j_2)^2,
\end{array}
\end{equation}
and
\begin{equation} \label{E:W(3)}
\begin{array}{l l}
W(\boldsymbol{v}_3) & \displaystyle
= \beta_{ij}^2\left(\frac{9}{4}\pi\mu R^{5/2}\delta_{ij}^{-5/2} +\frac{99}{160}\pi\mu R^{3/2}\delta_{ij}^{-3/2} +\frac{29241}{17920}\pi\mu R^{1/2}\delta_{ij}^{-1/2}\right)\\[7pt]
& \displaystyle + R^2(\omega^i-\omega^j)\beta_{ij}\left(-3\pi\mu R^{3/2}\delta_{ij}^{-3/2}+\frac{9}{40}\pi\mu R^{1/2}\delta_{ij}^{-1/2}\right)\\[7pt]
& \displaystyle + R^2(\omega^i-\omega^j)^2 \frac{3}{2}\pi\mu R^{1/2}\delta_{ij}^{-1/2}\\[7pt]
& \displaystyle +\mathbf{C}_5R^2(\omega^i-\omega^j)^2+\mathbf{C}_6R^2(\omega^i-\omega^j)\beta_{ij}+\mathbf{C}_7R^2\beta_{ij}^2.
\end{array}
\end{equation}

Coefficients  of the quadratic forms $Q_{ij}$ \eqref{E:Q_ij} can be read off from the formulas
 \eqref{E:W(1)}- \eqref{E:W(3)}.
Construction for the boundary necks and computation of the coefficients of $Q_{ij}$ \eqref{E:Q_ij-bound} is identical and it  is given in Appendix \ref{A:bound}.

\begin{remark}
Note that the first term of the decomposition \eqref{E:v3} indicates the contribution from the fluid motion due to the rotation of neighboring disks in opposite
directions while the second term describes the velocity of the fluid between motionless inclusions that is the Poiseuille microflow in a neck. Constants $\beta_{ij}$,
called \textit{permeation constants}, characterizing this microflow, are the local fluxes through the line segment $\ell_{ij}$
connecting two motionless inclusions and given by \eqref{E:beta0-new}.
\end{remark}

\subsection{Proof of Lemma \ref{L:upper}} \label{SS:proof_upper}

In the previous subsection we showed that for any $(\mathbb{U},\boldsymbol{\omega},\boldsymbol{\beta^*})\in \boldsymbol{\mathcal{R}}^*$ there exists a trial vector field
$\boldsymbol{v}_0$ on the fictitious fluid domain such that
\[
\boldsymbol{v}_0|_{\partial B^i}=\boldsymbol{U}^{i}+
R\omega^i(n_1^i\boldsymbol{e}_2-n_2^i\boldsymbol{e}_1), \quad
\boldsymbol{v}_0|_{\partial \Omega}=\boldsymbol{f},
\]
and for any neck $\Pi_{ij}$
\begin{equation}  \label{E:necks-estimate}
0\leq W_{\Pi_{ij}}(\boldsymbol{v}_0)-
W_{\Pi_{ij}}(\boldsymbol{u})\leq
C_1|\boldsymbol{U}^i-\boldsymbol{U}^j|^2+C_2R^2(\omega^i)^2+C_3R^2(\omega^j)^2+C_4R^2\beta_{ij}^2,
\end{equation}
where $\boldsymbol{u}$ is the minimizer of \eqref{E:W-fict}.

Hence, in order to prove Lemma \ref{L:ff-approx} it suffices to show that we can extend $\boldsymbol{v}_0$ in every triangle
$\triangle_{ijk}$ continuously so that this extension, called $\boldsymbol{w}$, satisfies
\[
\begin{array}{l l}
W_{\triangle_{ijk}}(\boldsymbol{w})  & \displaystyle \leq C_1|\boldsymbol{U}^i-\boldsymbol{U}^j|^2+C_2|\boldsymbol{U}^j-\boldsymbol{U}^k|^2+C_3|\boldsymbol{U}^k-\boldsymbol{U}^i|^2\\[5pt]
 & \displaystyle + C_4R^2(\omega^i)^2+C_5R^2(\omega^j)^2+C_6R^2(\omega^k)^2\\[5pt]
 & \displaystyle +C_7R^2\beta_{ij}^2+C_8R^2\beta_{jk}^2+C_9R^2\beta_{ki}^2,
\end{array}
\]
in every triangle $\triangle_{ijk}$. We construct this extension by assuming that $\boldsymbol{w}$ solves the Stokes equation on
$\triangle_{ijk}$\footnote{Here we discuss interior triangles, estimates for the triangles at the boundary are identical.}:
\begin{equation}  \label{E:problem-triangle}
\left\{
\begin{array}{r l l l}
(a) & \mu\triangle\boldsymbol{w}=\nabla p, & \mbox{in } \Delta_{ijk}\\
(b) & \nabla\cdot\boldsymbol{w}=0, & \mbox{in } \Delta_{ijk}\\
(c) & \boldsymbol{w}=\boldsymbol{v}_0, & \mbox{on }
\partial\Delta_{ijk}
\end{array}
\right.
\end{equation}

For \eqref{E:problem-triangle} by \cite{galdi} we have
\begin{equation}  \label{E:galdi}
W_{\triangle_{ijk}}\leq C\|\boldsymbol{v}_0-\boldsymbol{U}\|^2_{H^{1/2}(\partial\triangle_{ijk})},
\end{equation}
where $\boldsymbol{U}$ is any constant vector. The constant $C$ may in principle depend on the Lipschitz constant of
$\triangle_{ijk}$, but, since the close packing condition is assumed, angles in all $\triangle_{ijk}$ are bounded uniformly
from $0$ and $\pi$, and therefore, $C$ in \eqref{E:galdi} is universal.

\begin{figure}[!ht]
  \centering
  \includegraphics[scale=.55]{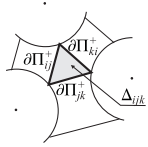}
  \caption{Boundary of the triangle $\triangle_{ijk}$}\label{F:bound-tr}
\end{figure}

Let $\boldsymbol{U}=1/3(\boldsymbol{U}^i+\boldsymbol{U}^j+\boldsymbol{U}^k)$.
Then $(\boldsymbol{v}_0-\boldsymbol{U})$ on $\partial\triangle_{ijk}$ is explicitly given. Namely,
\[
\begin{array}{r l}
(\boldsymbol{v}_{0}-\boldsymbol{U})|_{\partial\Pi_{ij}^+} = &
\displaystyle \frac{1}{6}(\boldsymbol{U}^i-\boldsymbol{U}^k)+
\frac{1}{6}(\boldsymbol{U}^j-\boldsymbol{U}^k)+
\left[(\boldsymbol{U}^i-\boldsymbol{U}^j)\cdot\boldsymbol{p}^{ij}+R\omega^i+R\omega^j\right]\boldsymbol{v}_{1}^{1} \\[5pt]
+ & \displaystyle R(\omega^i+\omega^j)\boldsymbol{v}_{1}^{2}+\left[(\boldsymbol{U}^i-\boldsymbol{U}^j)\cdot\boldsymbol{q}^{ij}\right]\boldsymbol{v}_2
+R(\omega^i-\omega^j)\boldsymbol{v}_{3}^{1}+\beta_{ij}\boldsymbol{v}_{3}^{2},
\end{array}
\]
where
$\partial\triangle_{ijk}=\partial\Pi^+_{ij}\cup\partial\Pi^+_{jk}\cup\partial\Pi^+_{ki}$
(see Fig. \ref{F:bound-tr}) and $\boldsymbol{v}_{1}^{1}$, $\boldsymbol{v}_{1}^{2}$, $\boldsymbol{v}_2$, $\boldsymbol{v}_{3}^{1}$,
$\boldsymbol{v}_{3}^{2}$ are polynomial vector fields independent of $(\mathbb{U},\boldsymbol{\omega},\boldsymbol{\beta}^*)$ and
$\delta$, and given in \eqref{E:v1}, \eqref{E:v2}, \eqref{E:v3}, and $\beta_{ij}$ is defined by \eqref{E:beta0}.
Therefore,
\begin{equation}  \label{E:estimate-triangle}
\begin{array}{l l}
\|\boldsymbol{v}_0-\boldsymbol{U}\|^2_{H^{1/2}(\partial\triangle_{ijk})} & \displaystyle \leq C_1|\boldsymbol{U}^i-\boldsymbol{U}^j|^2+C_2|\boldsymbol{U}^j-\boldsymbol{U}^k|^2+C_3|\boldsymbol{U}^k-\boldsymbol{U}^i|^2\\[5pt]
 & \displaystyle + C_4R^2(\omega^i)^2+C_5R^2(\omega^j)^2+C_6R^2(\omega^k)^2\\[5pt]
 & \displaystyle +C_7R^2\beta_{ij}^2+C_8R^2\beta_{jk}^2+C_9R^2\beta_{ki}^2.
\end{array}
\end{equation}

Combining \eqref{E:estimate-triangle} with \eqref{E:galdi} we obtain
\[
\begin{array}{r l}
\displaystyle \sum_{i \in \mathbb{I}} \sum_{j \in \mathcal{N}_i}
W_{\triangle_{ijk}}(\boldsymbol{w})\leq & \displaystyle \sum_{i
\in \mathbb{I}} \sum_{j \in \mathcal{N}_i}
C_1|\boldsymbol{U}^i-\boldsymbol{U}^j|^2+C_2|\boldsymbol{U}^j-\boldsymbol{U}^k|^2+C_3|\boldsymbol{U}^k-\boldsymbol{U}^i|^2\\[5pt]
 & \displaystyle + C_4R^2(\omega^i)^2+C_5R^2(\omega^j)^2+C_6R^2(\omega^k)^2\\[5pt]
 & \displaystyle +C_7R^2\beta_{ij}^2+C_8R^2\beta_{jk}^2+C_9R^2\beta_{ki}^2.
\end{array}
\]
Finally, the above inequality and \eqref{E:necks-estimate} give the desired result \eqref{E:lemma-main}.

\Section{Conclusions} \label{S:conclusions}

In this paper the asymptotic formula for the overall viscous
dissipation rate is obtained, where \emph{all} singular terms are
derived and justified. This is done by developing a new technical
tool - the two-step fictitious fluid approach. Such approach is
expected to be helpful in evaluation of effective properties of
various highly packed particulate composites. The obtained
asymptotics provides for a complete picture of microflows, while
previous studies, mentioned in Introduction, gave only partial
analysis of microflows and singular terms. In particular, a new
term due to the Poiseuille microflow, which was not taken into
account previously, is obtained. It is shown that this Poiseuille
microflow does not contribute to singular behavior of viscous
dissipation rate in 3D. In contrast, in 2D it may result in an
anomalous rate of blow-up (of order $\delta^{-5/2}$). Indeed, such
a rate of blow-up is obtained in the presence of external field
(e.g. gravity) or due to the boundary layer effects.
Our analysis suggests that in  absence of an external field the
anomalous blow-up does not occur.

The obtained asymptotics expresses the {\it continuum} dissipation
rate in terms of a {\it discrete} dissipation rate, and the latter
reveals dependence on the key physical parameters.

We remark that this work leads to a somewhat surprising
observation that suspensions are actually harder to analyze in 2D
than in 3D. As we mentioned above, the Poiseuille type microflow
is significant in 2D and it is negligible in 3D. The key reason
here is topological: in 2D thin gaps between closely spaced
inclusions partition the fluid domain into disconnected regions,
which is not the case in 3D. Hence, in 2D permeation of fluid
between two inclusions contributes into the singular behavior of
the overall viscous dissipation rate.

Finally, we note that 2D mathematical models were often used to analyze qualitative behavior of 3D problems in order to reduce the analytical and
computational complexity of the problem. Our work clearly shows
limits of validity of such modeling.

\appendix
\Section{Appendix}

\subsection{Coefficients of the quadratic form $Q$} \label{A:coeffs}

The coefficients $\boldsymbol{\mathcal{C}}_{k}^{ij}$, $k=1,\ldots,9$, and $\boldsymbol{\mathcal{B}}_{m}^{ij}$, $m=1,\ldots,14$ appearing in \eqref{E:Q_ij} and \eqref{E:Q_ij-bound}, respectively are given by
\begin{equation} \label{E:coefficients}
\begin{array}{l l l l }
& \displaystyle \boldsymbol{\mathcal{C}}_{1}^{ij} =
\frac{1}{2}\pi\mu \left(\frac{R}{d_{ij}}\right)^{1/2}, &
\displaystyle \boldsymbol{\mathcal{C}}_{2}^{ij} =
\frac{3}{4}\pi\mu \left(\frac{R}{d_{ij}}\right)^{3/2},
& \displaystyle \boldsymbol{\mathcal{C}}_{3}^{ij} = \frac{207}{320}\pi\mu \left(\frac{R}{d_{ij}}\right)^{1/2},\\[7pt]
& \displaystyle \boldsymbol{\mathcal{C}}_{4}^{ij} =
\frac{9}{4}\pi\mu \left(\frac{R}{d_{ij}}\right)^{5/2}, &
\displaystyle \boldsymbol{\mathcal{C}}_{5}^{ij} =
\frac{99}{160}\pi\mu \left(\frac{R}{d_{ij}}\right)^{3/2},
& \displaystyle \boldsymbol{\mathcal{C}}_{6}^{ij} = \frac{29241}{17920}\pi\mu \left(\frac{R}{d_{ij}}\right)^{1/2},\\[7pt]
& \displaystyle \boldsymbol{\mathcal{C}}_{7}^{ij} = -3\pi\mu
\left(\frac{R}{d_{ij}}\right)^{3/2}, & \displaystyle
\boldsymbol{\mathcal{C}}_{8}^{ij} = \frac{9}{40}\pi\mu
\left(\frac{R}{d_{ij}}\right)^{1/2},
& \displaystyle \boldsymbol{\mathcal{C}}_{9}^{ij} = \frac{3}{2}\pi\mu \left(\frac{R}{d_{ij}}\right)^{1/2},\\[7pt]
& \displaystyle \boldsymbol{\mathcal{B}}_{1}^{ij} = 18\pi\mu
\left(\frac{R}{d_{ij}}\right)^{5/2}, & \displaystyle
\boldsymbol{\mathcal{B}}_{2}^{ij} = \frac{51}{20}\pi\mu
\left(\frac{R}{d_{ij}}\right)^{3/2},
& \displaystyle \boldsymbol{\mathcal{B}}_{3}^{ij} = \frac{20889}{2240}\pi\mu \left(\frac{R}{d_{ij}}\right)^{1/2},\\[7pt]
& \displaystyle \boldsymbol{\mathcal{B}}_{4}^{ij} = 4\pi\mu
\left(\frac{R}{d_{ij}}\right)^{1/2}, & \displaystyle
\boldsymbol{\mathcal{B}}_{5}^{ij} = \frac{9}{2}\pi\mu
\left(\frac{R}{d_{ij}}\right)^{1/2},
& \displaystyle \boldsymbol{\mathcal{B}}_{6}^{ij} = 6\pi\mu \left(\frac{R}{d_{ij}}\right)^{3/2},\\[7pt]
& \displaystyle \boldsymbol{\mathcal{B}}_{7}^{ij} =
\frac{63}{20}\pi\mu \left(\frac{R}{d_{ij}}\right)^{1/2}, &
\displaystyle \boldsymbol{\mathcal{B}}_{8}^{ij} = 6\pi\mu
\left(\frac{R}{d_{ij}}\right)^{3/2},
& \displaystyle \boldsymbol{\mathcal{B}}_{9}^{ij} = \frac{19}{20}\pi\mu \left(\frac{R}{d_{ij}}\right)^{1/2},\\[7pt]
& \displaystyle \boldsymbol{\mathcal{B}}_{10}^{ij} = -3\pi\mu
\left(\frac{R}{d_{ij}}\right)^{3/2}, & \displaystyle
\boldsymbol{\mathcal{B}}_{11}^{ij} = -\frac{3}{8}\pi\mu
\left(\frac{R}{d_{ij}}\right)^{1/2},
& \displaystyle \boldsymbol{\mathcal{B}}_{12}^{ij} = -3\pi\mu \left(\frac{R}{d_{ij}}\right)^{1/2},\\[7pt]
& \displaystyle \boldsymbol{\mathcal{B}}_{13}^{ij} = -3\pi\mu
\left(\frac{R}{d_{ij}}\right)^{1/2}, & \displaystyle
\boldsymbol{\mathcal{B}}_{14}^{ij} = 6\pi\mu
\left(\frac{R}{d_{ij}}\right)^{1/2}.
\end{array}
\end{equation}

\subsection{Discrete Korn Inequality} \label{A:korn}

\begin{lemma}[Discrete Korn's Inequality] \label{L:korn}
Suppose $\Omega_F$ satisfies the close packing condition. Given $\{\boldsymbol{U}^i,i\in\mathbb{B}\}$ there exists a constant $C>0$ such that for every $\{\boldsymbol{U}^i,i\in\mathbb{I}\}$ the following inequality holds:
\begin{equation} \label{E:korn}
\sum_{i\in\mathbb{I}} \sum_{j\in \mathcal{N}_i}[(\boldsymbol{U}^i-\boldsymbol{U}^j)\cdot\boldsymbol{q}^{ij}]^2\geq
C\sum_{i\in\mathbb{I}} \sum_{j\in \mathcal{N}_i}(\boldsymbol{U}^i-\boldsymbol{U}^j)^2.
\end{equation}
\end{lemma}

\begin{proof}
By contradiction we assume that there exists $\{\boldsymbol{U}^i_{k}\}$ such that
\[
\sum_{i\in\mathbb{I}} \sum_{j\in \mathcal{N}_i}[(\boldsymbol{U}^i_{k}-\boldsymbol{U}^j_{k})\cdot\boldsymbol{q}^{ij}]^2<
\frac{1}{k}\sum_{i\in\mathbb{I}} \sum_{j\in \mathcal{N}_i}(\boldsymbol{U}^i_{k}-\boldsymbol{U}^j_{k})^2.
\]

Then we have two options here.

\underline{Case 1}. There exists a constant $M>0$ such that for all $i\in\mathbb{I}$ and $k\in\mathbb{N}$: $|\boldsymbol{U}^i_{k}|<M$.

\underline{Case 2}. There exists $i_0\in \mathbb{I}$ and a subsequence $\{k_n\}$ such that $\displaystyle \boldsymbol{U}^{i_0}_{k_n}\rightarrow \infty$ as $k_n\rightarrow \infty$.
Without lost of generality we set: $\displaystyle |\boldsymbol{U}^{i_0}_{k_n}|=\max_{i\in\mathbb{I}}|\boldsymbol{U}^{i}_{k_n}|$.

First we prove the Case 1. By Bolzano-Weierstrass lemma there exists a subsequence $\{k_n\}$ and a set of vectors $\{\boldsymbol{V}^{i},i\in\mathbb{I}\}$ such that
$\boldsymbol{U}^{i}_{k_n}\rightarrow\boldsymbol{U}^{i}$ as $k_n\rightarrow\infty$ for $i\in\mathbb{I}$ and for $\boldsymbol{V}^{i}$ we have:

$(a_1)$ $(\boldsymbol{V}^{i}-\boldsymbol{V}^{j})\cdot\boldsymbol{q}^{ij}=0$, $i,j\in\mathbb{I}\cup\mathbb{B}$;

$(b_1)$ $\boldsymbol{V}^{i}=\begin{pmatrix} a & b\\ c & -a \end{pmatrix}\begin{pmatrix} x \\ y \end{pmatrix}$, $i\in\mathbb{B}$.

Recall that in this paper we focus on the quasi-hexagonal array of inclusions which triangulates the domain $\Omega$ and in which
every disk typically has six neighbors. For such an array the condition $(a_1)$ implies that all disks, including quasi-disks,
are moving as a rigid body, that is, translate and rotate as a whole. Hence, there exist numbers $\alpha$, $\beta$, $\varphi$ such that
\begin{equation} \label{E:rigid-body}
\boldsymbol{V}^{i}=\begin{pmatrix} \alpha \\ \beta \end{pmatrix}+\begin{pmatrix} \cos \varphi & \sin \varphi \\ -\sin \varphi & \cos \varphi \end{pmatrix}\begin{pmatrix} x \\ y \end{pmatrix}.
\end{equation}
However, $(b_1)$ implies that $\alpha=\beta=0$, $\cos \varphi=a=0$, $b=-s=\sin \varphi$ which are excluded by our boundary conditions \eqref{E:ext-bc}.

For the Case 2 denote by $\displaystyle a_{k_n}:=|\boldsymbol{U}^{i_0}_{k_n}|=\max_{i\in\mathbb{I}}|\boldsymbol{U}^{i}_{k_n}| \rightarrow\infty$
as $k_n\rightarrow\infty$.

Consider $\boldsymbol{V}^{i}_{k_n}=\frac{\boldsymbol{U}^{i}_{k_n}}{a_{k_n}}$. We know that $|\boldsymbol{V}^{i}_{k_n}|\leq 1$, $|\boldsymbol{V}^{i_0}_{k_n}|= 1$. Hence,
there exists a subsequence $\{k_n\}$ and a set of vectors $\{\boldsymbol{V}^{i},i\in\mathbb{I}\}$  such that

$(a_2)$ $\boldsymbol{V}^{i}_{k_n}\rightarrow\boldsymbol{V}^{i}$ as $k_n\rightarrow\infty$;

$(b_2)$ $|\boldsymbol{V}^{i_0}|=1$;

$(c_2)$ $\boldsymbol{V}^{i}=\boldsymbol{0}$, for $i\in\mathbb{B}$;

$(d_2)$ $(\boldsymbol{V}^i-\boldsymbol{V}^j)\cdot\boldsymbol{q}^{ij}=0$, $i,j\in\mathbb{I}\cup\mathbb{B}$.

As before, the condition $(d_2)$ implies that all disks move as a rigid body, that is, \eqref{E:rigid-body} holds. Using $(c_2)$ we have that $\alpha=\beta=\varphi=0$, hence
$\boldsymbol{V}^{i}=\boldsymbol{0}$ for $i\in\mathbb{I}\cup\mathbb{B}$ which contradicts $(b_2)$.

\end{proof}

\subsection{Construction of the network} \label{A:network}

Here we describe how the network $\mathcal{G}$ is constructed.

\begin{figure}[h!]
  \centering
  \includegraphics[scale=0.95]{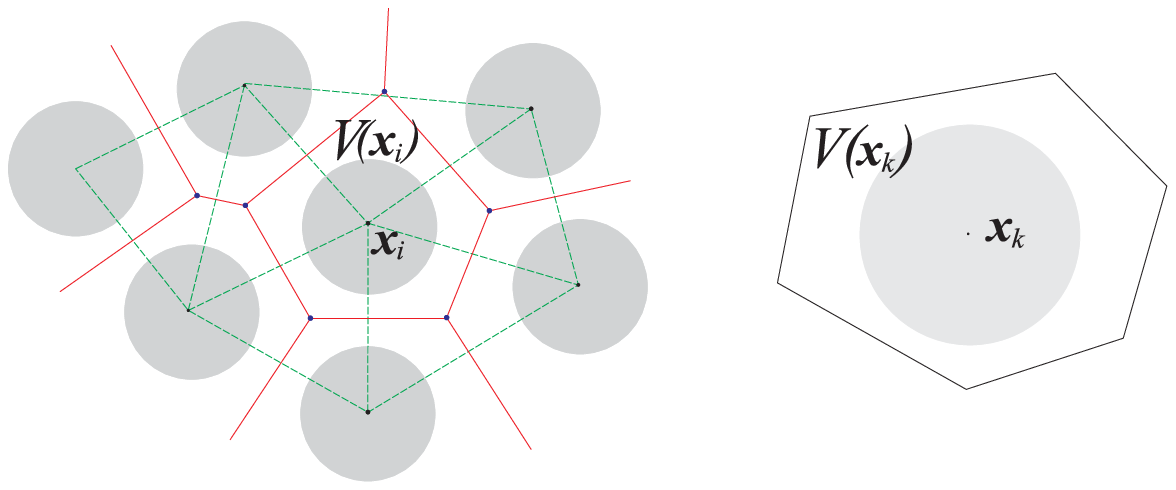}
  \caption{Voronoi diagram and Voronoi cell of the site $\boldsymbol{x}_{k}$} \label{F:voronoi}
\end{figure}

\begin{figure}[h!]
  \centering
  \includegraphics[scale=0.65]{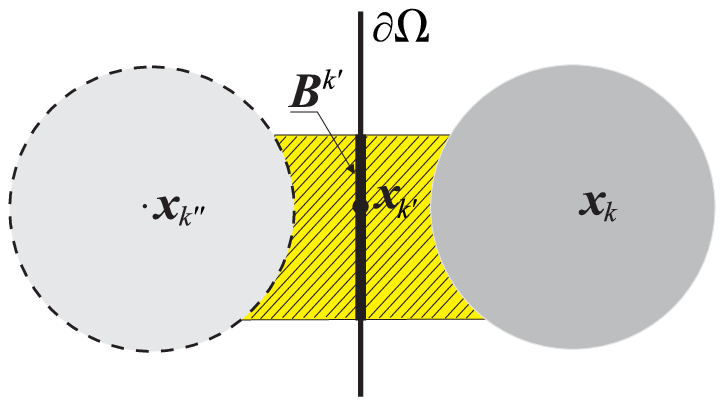}
  \caption{The quasidisk $B^{k'}$ centered at $\boldsymbol{x}_{k'}$ on the boundary of $\Omega$} \label{F:auxiliary}
\end{figure}

We follow \cite{bk} where the notion of neighboring inclusions were defined via the \textit{Voronoi tessellation (diagram)} \cite{fortune} of the domain $\Omega_F$
(Fig.\ref{F:voronoi}). For an arbitrary distribution of the centers $\boldsymbol{x}_{i}$ of inclusions $B^{i}$, called \textit{sites} or
\textit{vertices}:
\[
X=\{\boldsymbol{x}_{i} \in \Omega,\, i=1,\ldots, N\},
\]
the Voronoi diagram is the partition of the domain $\Omega$ into non-overlapping \textit{Voronoi cells} $V(\boldsymbol{x}_{i})$.
Each $V(\boldsymbol{x}_{i})$ is the set of all points in $\Omega$ that are closer to $\boldsymbol{x}_{i}$ than to any other site from $X$.
The situation when site $\boldsymbol{x}_{k}$ is located near one of boundaries of $\Omega$ is treated similarly to the one in \cite{bk,bn}.
Namely, one can consider the reflection $\boldsymbol{x}_{k'}$ of $\boldsymbol{x}_{k}$ along this boundary and construct the
corresponding Voronoi cell using the obtained auxiliary cite $\boldsymbol{x}_{k''}$ outside of $\Omega$ (see Fig.\ref{F:auxiliary}, where the hatched region
is a neck), introducing a new site $\boldsymbol{x}_{k'}$ on the boundary, which is the center of a line segment $B^{k'}$ in $\partial \Omega$,
called a quasidisk (see also \cite{bk,bn}).
Two disks are called \textit{neighbors} if their Voronoi cells share a common edge.

\begin{figure}[h!]
  \centering
  \includegraphics[scale=0.65]{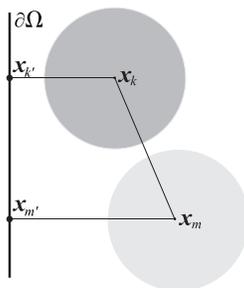}
  \caption{A ``triangle" near the boundary of $\Omega$} \label{F:auxiliary1}
\end{figure}

Now we can define a network $\mathcal{G}$ as introduced in Subsection \ref{SS:FFA}.

This graph is simply the Delaunay triangulation of $X$. Recall that the geometric dual of the Voronoi diagram in two
dimensions is a Delaunay triangulation \cite{fortune}. The basic element (or cell) of the Delaunay triangulation is a triangle that
may degenerate into quadrilateral or $n$-gon, when four or more disk centers lie on the same circle.

The \textit{neck-triangle partition} is a decomposition of the domain $\Omega_F$ into simple geometric objects:
necks and triangles, based on the \textit{central projection}.
\begin{figure}[ht]
  \begin{center}
 \includegraphics[scale=0.99]{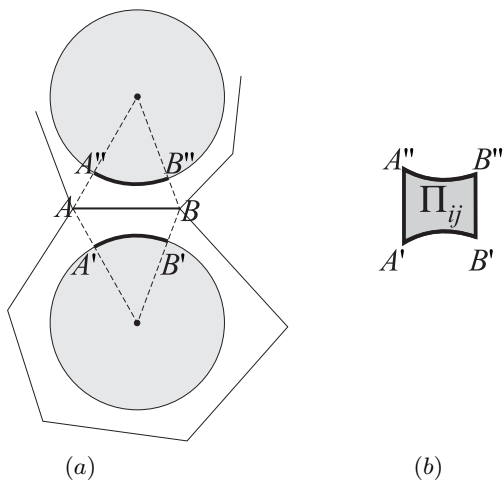}\\
 ($a$) \hspace{4cm} ($b$)
  \end{center}
  \caption{($a$) Central projection projects the Voronoi edge $AB$ onto arcs $A'B'$ and $A''B''$, ($b$) the neck $\Pi_{ij}$} \label{F:neck}
\end{figure}

The central projection partitions the boundary of a disk using the corresponding Voronoi cells as follows $\pi_{i}:\partial
V(\boldsymbol{x}_{i})\rightarrow \partial B^{i}$, $i \in \mathbb{I}$, where
$\pi_{i}(\boldsymbol{x})=R(\boldsymbol{x}-\boldsymbol{x}_{i})/|\boldsymbol{x}-\boldsymbol{x}_{i}|$.

As shown on Fig. \ref{F:neck}, the \textit{neck} $\Pi_{ij}$ connects arcs $A'B'$ with the arc $A''B''$, obtained by the
central projection of the same edge $AB$. We remark here that with this construction the edges of a neck are always parallel. Also a neck is necessarily symmetric
with respect to the line connecting centers of the neighboring disks. The set of all necks in $\Omega_F$ is $\boldsymbol{\Pi}$
(the shadowed region on Fig. \ref{F:f_domain}), called the fictitious fluid domain.

Once necks are defined the complementary part of the domain is the union of triangles (see e.g. $\triangle_{ijk}$ on Fig. \ref{F:f_domain}) and possibly
trapezoids near the boundary $\partial\Omega$. As mentioned above, by slight abuse of terminology we call such trapezoids ``triangles''. The set of all triangles
is $\displaystyle \boldsymbol{\Delta}=\bigcup_{i\in \mathbb{I}, j,k\in \mathcal{N}_i}\triangle_{ijk}$.

Thus, we obtain the \textit{neck-triangle partition} of the fluid domain: $\Omega_F=\boldsymbol{\Delta}\cup\boldsymbol{\Pi}$.
Such a partition of the fluid domain $\Omega_F$ allows us to decompose the functional $W_{\Omega_F}$ as follows:
\[\widehat{W}=W_{\Omega_F}(\cdot)=W_{\boldsymbol{\Pi}}(\cdot)+W_{\boldsymbol{\Delta}}(\cdot)=:\widehat{W}_{\boldsymbol{\Pi}}+\widehat{W}_\Delta.\]

\begin{figure}[!ht]
  \centering
  \includegraphics[scale=1.05]{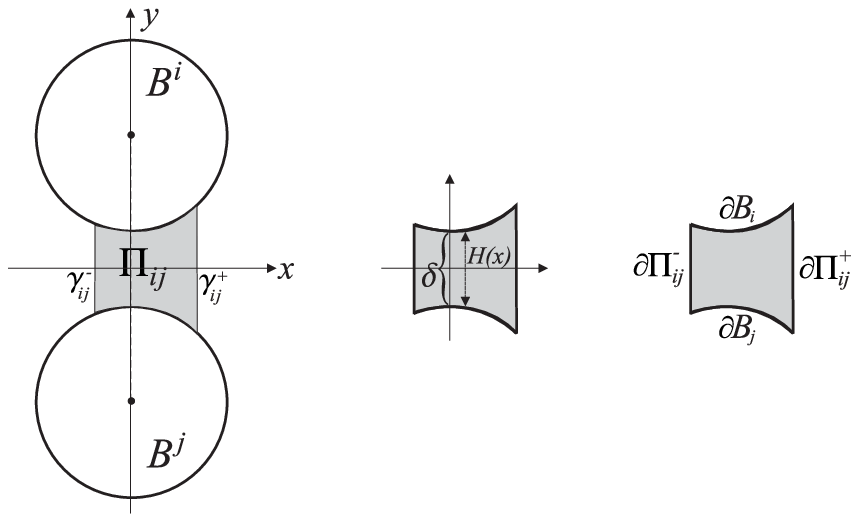}\\
  ($a$) \hspace{2.5cm}($b$)\hspace{2.5cm}($c$)
  \caption{($a$) The neck $\Pi_{ij}$ connecting two neighbors $B^i$ and $B^j$, ($b$) the distance $H(x)$ between two neighbors,
  ($c$) the boundary $\partial\Pi_{ij}$ of the neck}\label{F:neck0}
\end{figure}

\subsection{Auxiliary Problem for \eqref{E:ex_form}} \label{A:hasimoto}

Here we prove Proposition \ref{P:discussion}. First we show that it is possible to choose such a density $\rho_s$ and boundary data $\boldsymbol{f}$ that disks do not move.

\begin{theorem}  \label{P:hasimoto}
There exist $\rho_s$ and $\boldsymbol{f}$ such that $\boldsymbol{U}^i=0$ and $\omega^i=0$, $i=1,\ldots,N$.
\end{theorem}

\begin{figure}[!ht]
  \centering
  \includegraphics[scale=1.05]{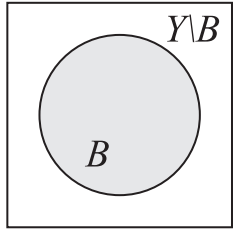}
  \caption{Periodicity cell $Y$}\label{F:hasimoto}
\end{figure}

\textit{Proof}. On the periodicity cell $Y=\left(-\frac{1}{2},\frac{1}{2}\right)\times\left(-\frac{1}{2},\frac{1}{2}\right)$ (Fig.\ref{F:hasimoto}) consider:
\begin{equation}   \label{E:hasimoto}
\left\{
\begin{array}{r l l}
(a) & \displaystyle \mu\triangle \boldsymbol{v}=\nabla p, & \boldsymbol{x} \in Y\backslash B \\[5pt]
(b) & \displaystyle \nabla\cdot \boldsymbol{v}=0, &  \boldsymbol{x} \in Y\backslash B \\[5pt]
(c) & \displaystyle \boldsymbol{v}=\boldsymbol{0}, & \boldsymbol{x} \in \partial B\\[5pt]
(d) & \displaystyle
\langle\boldsymbol{v}\rangle=\frac{1}{|Y|}\int_{Y}\boldsymbol{u}d\boldsymbol{x}=
\begin{pmatrix} 0 \\ 1 \end{pmatrix}\\[10pt]
(e) & \displaystyle \boldsymbol{v}\mbox{ is } Y-\mbox{periodic } \\[5pt]
\end{array}
\right.
\end{equation}
This problem is equivalent to the one considered in \cite{has}.

The existence of the unique solution to \eqref{E:hasimoto} can be verified by the classical energy methods. Observe
that the solution $\boldsymbol{v}=(v_1,v_2)$ of \eqref{E:hasimoto} satisfies the equations $(a)$, $(b)$ of the problem
\eqref{E:ex_form}, and the condition (\ref{E:ex_form}$c$) (with $\boldsymbol{U}^i=\boldsymbol{0}$, $\omega^i=0$) for $\Omega=Y$ (i.e. $N=1$).
Using symmetry of the problem one can show that $\boldsymbol{v}$ also satisfies the balance of the angular momentum (\ref{E:ex_form}$e$) with $\Omega=Y$:
\[
\int_{\partial B}\boldsymbol{\sigma}(\boldsymbol{v})\boldsymbol{n}\times\boldsymbol{n}ds\left(= \int_{\partial B}\left([\boldsymbol{\sigma}
(\boldsymbol{v})\boldsymbol{n}]_1n_2 -[\boldsymbol{\sigma}(\boldsymbol{v})\boldsymbol{n}]_2n_1 \right)ds\right)=0,
\]
\[
\int_{\partial B}\boldsymbol{\sigma}(\boldsymbol{v})\boldsymbol{n} ds=\begin{pmatrix} 0 \\ a \end{pmatrix},
\]
for some $a \in \mathbb{R}$ which is not, in general, zero (this integral is the force exerted by the fluid on the inclusion $B$).

We estimate $a$ in terms of $\delta$. For this we compare problems \eqref{E:ex_form} and
\eqref{E:hasimoto} for $\boldsymbol{u}$ and $\boldsymbol{v}$,
respectively. If we choose $\boldsymbol{f}$ in
(\ref{E:ex_form}$e$) to be a restriction of $\boldsymbol{v}$ on
$\partial\Omega$ and $a=mg=\pi R^2(\rho_s-\rho_f)g$ then
$\boldsymbol{v}$ solves \eqref{E:ex_form}. Hence,
$\boldsymbol{u}=\boldsymbol{v}$ in $\Omega=Y$, where, in
particular, the translational and angular velocities are zeros:
$\boldsymbol{U}=\boldsymbol{0}$, $\omega=0$.

Consider now the problem \eqref{E:ex_form} for $\boldsymbol{u}$
obtained by a periodic extension of $\boldsymbol{v}$ defined in $Y$ to $N$ cells. If we now choose $\rho_s$ and $\boldsymbol{f}$ as above, namely,
\begin{equation} \label{E:rho_s,f}
\rho_s=-\frac{a}{N \pi R^2 g}+\rho_f, \quad
\boldsymbol{f}=\boldsymbol{v}|_{\partial\Omega},
\end{equation}
then uniqueness implies that all inclusions are motionless (due to (\ref{E:hasimoto}$c$)). \hspace{5.5cm} $\Box$

Theorem \ref{P:hasimoto} shows that if we choose the density
$\rho_s$ and the external boundary conditions $\boldsymbol{f}$ in
\eqref{E:ex_form} so that inclusions are motionless then the only
microflow between neighboring inclusions is permeation (or the
Poiseuille microflow). We now show that the overall viscous
dissipation rate of \eqref{E:ex_form}
exhibits a superstrong blow up. Hereinafter we use the equivalence
of \eqref{E:ex_form} and \eqref{E:hasimoto} extended by
periodicity to $\Omega$. Before we obtain the asymptotics of the
dissipation rate $\widehat{W}$ for this problem we find an asymptotics of $\rho_s$ as $\delta\rightarrow 0$ and the form of $\boldsymbol{f}$.

\begin{lemma}
The following asymptotics holds for the force exerted by the fluid on the disks $B$ defined by the right hand side of (\ref{E:hasimoto}$d$):
\[
a=C\delta^{-5/2}+O(\delta^{-3/2}) \mbox{ with some constant } C=C(\mu, R, N)>0 \mbox{ as } \delta\rightarrow 0.
\]
\end{lemma}

\textit{Proof.} Consider an auxiliary function
$\boldsymbol{w}=\boldsymbol{v}-\langle\boldsymbol{v}\rangle$ that
solves:
\begin{equation}   \label{E:hasimoto-aux}
\left\{
\begin{array}{r l l}
(a) & \displaystyle \mu\triangle \boldsymbol{w}=\nabla p, & \boldsymbol{x} \in Y\backslash B \\[5pt]
(b) & \displaystyle \nabla\cdot \boldsymbol{w}=0, &  \boldsymbol{x} \in Y\backslash B \\[5pt]
(c) & \displaystyle \boldsymbol{w}=\begin{pmatrix} 0 \\ -1 \end{pmatrix}, & \boldsymbol{x} \in \partial B\\[5pt]
(d) & \displaystyle \langle\boldsymbol{w}\rangle=\boldsymbol{0}\\[10pt]
(e) & \displaystyle \boldsymbol{w}\mbox{ is periodic } \\[5pt]
\end{array}
\right.
\end{equation}
Then if we multiply (\ref{E:hasimoto}$a$) by $\boldsymbol{w}$ and
integrate by parts over $Y\backslash B$ we obtain:
\[
0=-\int_{Y\backslash
B}\boldsymbol{D}(\boldsymbol{w}):\boldsymbol{\sigma}(\boldsymbol{v})d\boldsymbol{x}+
\int_{\partial
B}\boldsymbol{w}\cdot\boldsymbol{\sigma}(\boldsymbol{v})\boldsymbol{n}ds+
\int_{\partial
Y}\boldsymbol{w}\cdot\boldsymbol{\sigma}(\boldsymbol{v})\boldsymbol{n}ds.
\]
Taking into account that
$\boldsymbol{D}(\boldsymbol{w})=\boldsymbol{D}(\boldsymbol{v})$
and (\ref{E:hasimoto-aux}$c$) we continue
\[
\begin{array}{l l l}
0=& \displaystyle -\int_{Y\backslash
B}\boldsymbol{D}(\boldsymbol{v}):\boldsymbol{\sigma}(\boldsymbol{v})d\boldsymbol{x}+
\int_{\partial B}[\boldsymbol{\sigma}(\boldsymbol{v})\boldsymbol{n}]_2ds+
\int_{\partial Y}\boldsymbol{w}\cdot\boldsymbol{\sigma}(\boldsymbol{v})\boldsymbol{n}ds\\[7pt]
= & \displaystyle -W_{Y\backslash B}(\boldsymbol{v})+a+
\int_{\partial
Y}\boldsymbol{w}\cdot\boldsymbol{\sigma}(\boldsymbol{v})\boldsymbol{n}ds.
\end{array}
\]
Due to periodicity of the function $\boldsymbol{w}$, indicated
above properties of the components of the vector
$\boldsymbol{\sigma}(\boldsymbol{v})\boldsymbol{n}$ and the
symmetry of the domain $Y$ we obtain that $\displaystyle
\int_{\partial
Y}\boldsymbol{w}\cdot\boldsymbol{\sigma}(\boldsymbol{v})\boldsymbol{n}ds=0$.
Hence, $\displaystyle a=W_{Y\backslash B}(\boldsymbol{v})$.
Therefore, in order to obtain an estimate for $a$ with respect to
$\delta$ we have to study a problem of the asymptotics of the
overall viscous dissipation rate of the suspension given by
\eqref{E:ex_form}.
Since inclusions in problem \eqref{E:ex_form} are motionless then the fluid permeates through thin gaps between them.
 So in view of formula \eqref{E:form-decomp} of Theorem \ref{T:main-thm3}  the viscous dissipation rate of this suspension
 exhibits the superstrong blow up of order $O(\delta^{-5/2})$ if we  prove that $\beta_{ij}=O(1)$ for at least one pair $(i,j)$: $i\in \mathbb{I}, j\in \mathcal{N}_i$.

Recall that $\boldsymbol{f}$ now is given by \eqref{E:rho_s,f} and
$v_1(-\frac{1}{2},y)=v_1(\frac{1}{2},y)$ due to periodicity of
$\boldsymbol{v}$. But $v_1$ is odd in $x$, hence, $v_1$ does not
depend on $x$. Similarly,
$v_1(x,-\frac{1}{2})=v_1(x,\frac{1}{2})=0$ due of oddness of $v_1$
in $y$. Thus, $v_1(x,y)=0$. Moreover, due to the periodicity of
$v_2(x,y)$ we have that the second component $f_2(x,y)$ of
$\boldsymbol{f}$ is even in both $x$ and $y$. Hence,
$\boldsymbol{f}=\begin{pmatrix} 0 \\ f_{2}(x,y)\end{pmatrix}$
where $f_2$ is an even function of both $x$ and $y$ on
$\partial\Omega$.

\begin{figure}[!ht]
  \centering
  \includegraphics[scale=0.85]{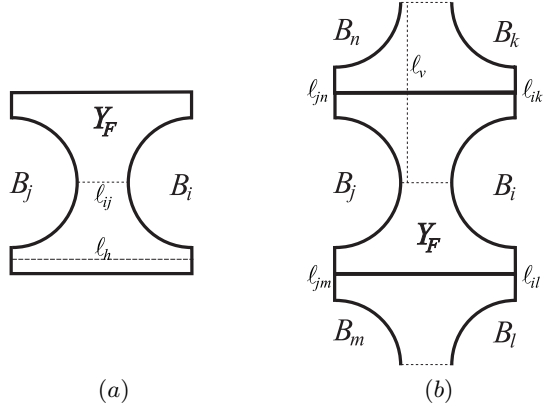}\\
  $(a)$ \hspace{3.7cm} $(b)$
  \caption{Periodicity cell $Y_F$}\label{F:per_ex1}
\end{figure}

Since for $\rho_s$ and $\boldsymbol{f}$ given by \eqref{E:rho_s,f}
we have $\boldsymbol{U}^i=\boldsymbol{0}$, and $\omega^i=0$ for
all $i$, the problem \eqref{E:ex_form} becomes the problem of
minimization of $W_{\Omega_F}(\boldsymbol{v})$ introduced in
\eqref{E:W} over the class $V$ defined by \eqref{E:v_form-var}, where
$\boldsymbol{U}^i=\boldsymbol{0}$, and $\omega^i=0$ for all $i$,
considered above. As mentioned above $a=\widehat{W}$ asymptotics of which is given by Propostion \ref{E:main-thm3}. To find an asymptotics of $\widehat{W}$
we need to know the flux $\beta_{ij}$ through the line $\ell_{ij}$ joining two neighboring disks $B^i$ and $B^j$.
In order to find this parameter consider a periodicity cell $Y_F$ as
in Fig. \ref{F:per_ex1}($a$). We note that due to divergence-free
property and periodicity of $\boldsymbol{u}$ one can obtain that
$\beta_{ij}$ is equal to the flux through any horizontal line
$\ell_{h}$ connecting the lateral boundaries of $Y_F$:
\[
\beta_{ij}=\frac{1}{R}\int_{\ell_{ij}}\boldsymbol{u}\cdot\boldsymbol{n}ds=\frac{1}{R}\int_{-\delta/2}^{\delta/2}u_2(x,0)dx=\frac{1}{R}
\int_{\ell_{h}}u_2(x,y)dx,
 \]
 for any fixed  $y\in(-1/2,1/2)$.
Using this fact, $\boldsymbol{u}=\boldsymbol{0}$ on $\partial B$, $\langle\boldsymbol{u}\rangle=\begin{pmatrix} 0 \\
1 \end{pmatrix}$ and Fubini's theorem we have
\[
1=\frac{1}{|Y|}\int_{Y}u_2(x,y)d\boldsymbol{x}=
\int_{-\frac{1}{2}}^{\frac{1}{2}}dy\int_{-\frac{1}{2}}^{\frac{1}{2}}u_2(x,y)dx=
\int_{-\frac{1}{2}}^{\frac{1}{2}}dy\int_{\ell_{h}}u_2(x,y)dx=R\beta_{ij}.
\]
In Fig. \ref{F:per_ex1}$a$ $\ell_{ij}$ is shown as a horizontal line.
There are also vertical half-middle lines connecting $B^i$ and
$B^j$ with their neighbors in the vertical direction ($B^k$, $B^l$
and $B^n$, $B^m$, respectively, in Fig. \ref{F:per_ex1}$b$). As
before using the periodicity and the divergence-free properties of
$\boldsymbol{u}$ we obtain that:
\[
\beta_{ik}=\beta_{jn}=\frac{1}{R}\int_{\ell_v}u_1(x,y)dy \mbox{ for any fixed
} x\in\left(-\frac{1}{2},\frac{1}{2}\right).
\]
where $\ell_{v}$ is any vertical line in $Y_F$
(Fig. \ref{F:per_ex1}$b$). And finally using the average of
$\boldsymbol{u}$ and that $\boldsymbol{u}$ is zero on the
inclusions we obtain:
\[
0=\frac{1}{|Y|}\int_{Y}u_1(x,y)d\boldsymbol{x}=
\int_{-\frac{1}{2}}^{\frac{1}{2}}dx\int_{\ell_{v}}u_1(x,y)dy=R\beta_{ik}=R\beta_{jn}.
\]
Proceeding similarly one can obtain that $\beta_{jm}=\beta_{jn}=0$.

Hence, we obtain that the flux through any horizontal line
connecting the neighbors equals $1$, while the flux through any
vertical line is zero. Also we note that the number of pairs
$(i,j)$ for which $\beta_{ij}$ is not zero is equal to the total
number $N$ of the inclusions in the suspension.

The fact that all inclusions are motionless and that there are $N$
fluxes though the middle lines connecting neighbors of our network
together with Proposition \ref{T:main-thm3} results in $a=N\widehat{W}=\frac{9}{4}\pi\mu N R^{5/2}\delta^{-5/2}+O(\delta^{-3/2})$.

\Section{Appendix}

\subsection{Details of construction of the trial fields} \label{A:constr-in}

\subsubsection{Verification of \eqref{E:error-in-1}}

Consider the minimizer $\displaystyle \boldsymbol{u}_{sh}$ of $W$ over $V_{sh}$.
The test function $\displaystyle \boldsymbol{v}_{1} \in V_{sh}$ must satisfy the boundary conditions \eqref{E:G1} on $\partial B^i$ and $\partial B^j$, which
can be rewritten as
\begin{equation}   \label{E:G1-new}
\begin{array}{l l}
\boldsymbol{G}^i_1 & \displaystyle = (U^i_1-U^j_1+R\omega^i+R\omega^j)\left(\frac{1}{2}\boldsymbol{e}_1\right)+R(\omega^i+\omega^j)\begin{pmatrix} \frac{1}{2}\left[\sqrt{1-\frac{x^2}{R^2}}-1\right]\\ \frac{x}{2R}\end{pmatrix} \\[15pt]
& \displaystyle =:
(U^i_1-U^j_1+R\omega^i+R\omega^j)\boldsymbol{g}^i_1+R(\omega^i+\omega^j)\boldsymbol{g}^i_2\,,\\[7pt]
\boldsymbol{G}^j_1 & \displaystyle = (U^i_1-U^j_1+R\omega^i+R\omega^j)\left(-\frac{1}{2}\boldsymbol{e}_1\right)+R(\omega^i+\omega^j)\begin{pmatrix} -\frac{1}{2}\left[\sqrt{1-\frac{x^2}{R^2}}-1\right]\\ \frac{x}{2R}\end{pmatrix}\\[15pt]
& \displaystyle =: (U^i_1-U^j_1+R\omega^i+R\omega^j)\boldsymbol{g}^j_1+R(\omega^i+\omega^j)\boldsymbol{g}^j_2\,.
\end{array}
\end{equation}

Hence, $\boldsymbol{v}_{1}$ can be written as
\begin{equation}   \label{E:v1-decomp}
\boldsymbol{v}_{1}=(U^i_1-U^j_1+R\omega^i+R\omega^j)\boldsymbol{v}^1_1+R(\omega^i+\omega^j)\boldsymbol{v}^2_1,
\end{equation}
where
\[
\left.\boldsymbol{v}^1_1\right|_{\partial B^i}=\boldsymbol{g}^i_1, \quad \left.\boldsymbol{v}^2_1\right|_{\partial B^i}=\boldsymbol{g}^i_2, \quad\quad\quad
\left.\boldsymbol{v}^1_1\right|_{\partial B^j}=\boldsymbol{g}^j_1, \quad \left.\boldsymbol{v}^2_1\right|_{\partial B^j}=\boldsymbol{g}^j_2.
\]
Using the decomposition \eqref{E:v1-decomp} we define $\boldsymbol{v}_{1}$ by \eqref{E:v1}-\eqref{E:F,G,K,M-1}.

Next we evaluate the difference between $W(\boldsymbol{v}_1)$ and $W(\boldsymbol{u}_{sh})$. Here denote
$\xi=U_1^i-U_1^j+R\omega^i+R\omega^j$ and $\zeta=R(\omega^i+\omega^j)$ then
\[
W(\boldsymbol{u}_{sh})=W\left(\xi\boldsymbol{u}_{sh}^1+\zeta\boldsymbol{u}_{sh}^2\right),
\]
where $\boldsymbol{u}_{sh}^1$ is the minimizer of $W$ satisfying $\boldsymbol{g}^i_1$, $\boldsymbol{g}^j_1$ and
$\boldsymbol{u}_{sh}^2$ is the minimizer of $W$ satisfying $\boldsymbol{g}^i_2$, $\boldsymbol{g}^j_2$ on $\partial B^i$, $\partial B^j$, respectively.

Compute the functional $W(\boldsymbol{v}_1)$ as $\delta\rightarrow 0$:
\begin{equation} \label{E:W(v1)-new}
\begin{array}{l l}
W(\boldsymbol{v}_1) & = \displaystyle \xi^2W(\boldsymbol{v}_1^1)+\zeta^2W(\boldsymbol{v}_1^2)+2\xi\zeta W(\boldsymbol{v}_1^1,\boldsymbol{v}_1^2)\\[5pt]
& = \displaystyle \frac{1}{2}(U_1^i-U_1^j+R\omega^i+R\omega^j)^2 \pi\mu R^{1/2} \delta_{ij}^{-1/2}+\zeta^2W(\boldsymbol{v}_1^2)+2\xi\zeta W(\boldsymbol{v}_1^1,\boldsymbol{v}_1^2)\\[5pt]
& =: \displaystyle (U_1^i-U_1^j+R\omega^i+R\omega^j)^2 \boldsymbol{\mathcal{C}}_1^{ij}\delta^{-1/2}+\zeta^2W(\boldsymbol{v}_1^2)+2\xi\zeta W(\boldsymbol{v}_1^1,\boldsymbol{v}_1^2),
\end{array}
\end{equation}
where $W(\boldsymbol{v}_1^2)=O(1)$, $W(\boldsymbol{v}_1^1,\boldsymbol{v}_1^2)=O(1)$ and
\begin{equation} \label{E:coeff1}
\boldsymbol{\mathcal{C}}_1^{ij}= \frac{1}{2}\pi\mu R^{1/2}d_{ij}^{-1/2}.
\end{equation}

Due to parallelogram identity for any vector fields $\boldsymbol{\varphi}$ and $\boldsymbol{\psi}$:
\[
4 W(\boldsymbol{\varphi},\boldsymbol{\psi})=W(\boldsymbol{\varphi}+\boldsymbol{\psi})-W(\boldsymbol{\varphi}-\boldsymbol{\psi}),
\]
we obtain
\begin{equation} \label{E:diff-1}
\begin{array}{ r l l}
& \displaystyle |W(\boldsymbol{v}_1)-W(\boldsymbol{u}_{sh})|
=|W(\xi\boldsymbol{v}^1_1+\zeta\boldsymbol{v}^2_1)-W(\xi\boldsymbol{u}_{sh}^1+\zeta\boldsymbol{u}_{sh}^2)|\\[5pt]
\leq & \displaystyle \xi^2\left|W(\boldsymbol{v}^1_1)-W(\boldsymbol{u}_{sh}^1)\right|+\zeta^2\left|W(\boldsymbol{v}^2_1)-W(\boldsymbol{u}^2_{sh})\right|\\[5pt]
+ & \displaystyle \frac{\xi\zeta}{2}\left|W(\boldsymbol{v}^1_1+\boldsymbol{v}^2_1)-W(\boldsymbol{u}_{sh}^1+\boldsymbol{u}_{sh}^2)\right|
+\frac{\xi\zeta}{2}\left| W(\boldsymbol{v}^1_1-\boldsymbol{v}^2_1) - W(\boldsymbol{u}_{sh}^1-\boldsymbol{u}_{sh}^2)\right|.
\end{array}
\end{equation}
In order to evaluate these differences in \eqref{E:diff-1} we construct four trial tensors $\boldsymbol{\mathcal{S}}_1^k \in F$, $k=1,\ldots,4$, such that,
\[
\begin{array}{ r l l}
& W^*_1(\boldsymbol{\mathcal{S}}_{1}^{1})\leq W^*_1\left(\boldsymbol{\sigma}(\boldsymbol{u}_{sh}^1)\right),
& W^*_1(\boldsymbol{\mathcal{S}}_{1}^{2})\leq W^*_1\left(\boldsymbol{\sigma}(\boldsymbol{u}_{sh}^2)\right),\\[3pt]
& W^*_1(\boldsymbol{\mathcal{S}}_{1}^{3})\leq W^*_1\left(\boldsymbol{\sigma}(\boldsymbol{u}_{sh}^1+\boldsymbol{u}_{sh}^2)\right),
& W^*_1(\boldsymbol{\mathcal{S}}_{1}^{4})\leq W^*_1\left(\boldsymbol{\sigma}(\boldsymbol{u}_{sh}^1-\boldsymbol{u}_{sh}^2)\right),
\end{array}
\]
obtaining
\[
\begin{array}{ r l l}
\displaystyle |W(\boldsymbol{v}_1)-W(\boldsymbol{u}_{sh})|
\leq & \displaystyle \xi^2\left|W(\boldsymbol{v}_1^1)-W^*_1(\boldsymbol{\mathcal{S}}_{1}^{1})\right|+
\zeta^2[W(\boldsymbol{v}_1^2)-W^*_1(\boldsymbol{\mathcal{S}}_{1}^{4})]\\ [5pt]
+ & \displaystyle \frac{\xi\zeta}{2}\left|W(\boldsymbol{v}_1^1+\boldsymbol{v}_1^2)-W^*_1(\boldsymbol{\mathcal{S}}_1^2)\right|+
\frac{\xi\zeta}{2}\left|W(\boldsymbol{v}_1^1-\boldsymbol{v}_1^2)-W^*_1(\boldsymbol{\mathcal{S}}_1^3)\right|.
\end{array}
\]

To construct e.g. $\boldsymbol{\mathcal{S}}_{1}^{1}$ we use the following observation. If $\boldsymbol{v}_1^1$ is a ``good" approximation of $\boldsymbol{u}^1_{sh}$ then we expect that
$\boldsymbol{\sigma}(\boldsymbol{v}_1^1)=2\mu\boldsymbol{D}(\boldsymbol{v}_1^1)-p_{v_1^1}\boldsymbol{I}$, where
\[
\boldsymbol{D}(\boldsymbol{v}_1^1)=
\begin{pmatrix} yG' & \frac{1}{2}(G+F'-\frac{y^2}{2}G'')\\
\frac{1}{2}(G+F'-\frac{y^2}{2}G'') & -yG'  \end{pmatrix}, \quad p_{v_1^1}=0,
\]
with $G(x)$ and $F(x)$ defined in \eqref{E:F,G,K,M-1}, is a ``good" approximation of $\boldsymbol{\sigma}(\boldsymbol{u}^1_{sh})$.
Then we correct $\boldsymbol{\sigma}(\boldsymbol{v}_1^1)$ so that it belongs to the set $F$, thus,
obtaining $\boldsymbol{\mathcal{S}}_1^1$.

Now define
\begin{equation}   \label{E:S1(1)-new}
\boldsymbol{\mathcal{S}}_{1}^{1} =
\mu \begin{pmatrix} \displaystyle 0 & \displaystyle G(x)-C \\ \displaystyle G(x)-C &  \displaystyle -yG'(x) \end{pmatrix} \quad \mbox{in} \quad \Pi_{ij},
\end{equation}
with $G(x)$ given by \eqref{E:F,G,K,M-1}, and the constant $\displaystyle C=\frac{1}{H(\gamma_{ij})}$. Clearly, $\boldsymbol{\mathcal{S}}_{1}^{1} \in F$.

When we evaluate the difference between $W(\boldsymbol{v}_{1}^1)$ defined by \eqref{E:W(v1)-new} and
$W^*_1(\boldsymbol{\mathcal{S}}_{1}^1)$, defined by \eqref{E:dual-fnl-1} we obtain the following estimate:
\[
\left|W(\boldsymbol{v}_1^1)-W^*_1(\boldsymbol{\mathcal{S}}_{1}^{1})\right|=O(1), \quad \mbox{as } \delta\rightarrow 0.
\]

Now choose
\[
\boldsymbol{\mathcal{S}}_1^2=\begin{pmatrix} 0 & 0 \\ 0 & 0 \end{pmatrix} \in F,
\]
hence, $W^*_1(\boldsymbol{\mathcal{S}}_1^2)=0$.

We also choose $\boldsymbol{\mathcal{S}}_1^3=\boldsymbol{\mathcal{S}}_1^1$ and $\boldsymbol{\mathcal{S}}_1^4=\boldsymbol{\mathcal{S}}_1^1$ and compute
the corresponding functionals obtaining
\[
\left|W(\boldsymbol{v}_1^1+\boldsymbol{v}_1^2)-W^*_1(\boldsymbol{\mathcal{S}}_1^1)\right|=O(1), \quad
\left|W(\boldsymbol{v}_1^1-\boldsymbol{v}_1^2)-W^*_1(\boldsymbol{\mathcal{S}}_1^1)\right|=O(1),
\]
as $\delta\rightarrow 0$. Therefore, $W(\boldsymbol{v}_1)$ determines the asymptotics of $W(\boldsymbol{u}_{sh})$:
\begin{equation} \label{E:difference-1}
\left|W(\boldsymbol{v}_1)-W(\boldsymbol{u}_{sh})\right|\leq E_1(\boldsymbol{U}^i,\boldsymbol{U}^j,\omega^i,\omega^j), \quad \mbox{as }\delta\rightarrow 0,
\end{equation}
where the error term $E_1$ is defined in \eqref{E:error-in-1}.
Hence,
\begin{equation} \label{E:W(u1)}
W(\boldsymbol{u}_{sh})=[(\boldsymbol{U}^i-\boldsymbol{U}^j)\cdot\boldsymbol{p}^{ij}+R\omega^i+R\omega^j]^2 \boldsymbol{\mathcal{C}}_1^{ij}\delta^{-1/2}+
E_1(\boldsymbol{U}^i,\boldsymbol{U}^j,\omega^i,\omega^j) \mbox{ as }\delta\rightarrow 0.
\end{equation}

\subsubsection{Verification of \eqref{E:error-in-2}}

Consider the minimizer $\displaystyle \boldsymbol{u}_{sq}$ of $W$ over $V_{sq}$. Choose $\boldsymbol{v}_2$ by \eqref{E:v2}-\eqref{E:G&F-v2}.
Introduce an approximate pressure
\begin{equation}  \label{E:p(v2)}
p_{v_2}=\mu\left(6\int_{\gamma_0}^x F(t)dt+G'-3y^2F'\right),
\end{equation}
where $\gamma_0$ is some constant and compute
\begin{equation}  \label{E:S(v2)}
\boldsymbol{\sigma}(\boldsymbol{v}_2)=(U^i_2-U^j_2)\begin{pmatrix} 2\mu(G'+3y^2F')-p_{v_2}
& \mu(6yF-yG''-y^3F'') \\  \mu(6yF-yG''-y^3F'') & -2\mu(G'+3y^2F')-p_{v_2}
\end{pmatrix}.
\end{equation}
Now we construct a corrector $\boldsymbol{\mathcal{S}}_{2}=(\boldsymbol{\mathcal{S}}_{2})_{i,j=1,2}$ to \eqref{E:S(v2)} so that $\boldsymbol{\mathcal{S}}_{2} \in F$.
Define $\boldsymbol{\mathcal{S}}_{2}$ by \eqref{E:S2(1)}
where constants $C_i$ ($i=1..4$) are chosen so that $\boldsymbol{\mathcal{S}}_{2}\boldsymbol{n}=\boldsymbol{0}$ on $\partial \Pi_{ij}^\pm$, that is,
\[
\begin{array}{l l l}
& \displaystyle C_1=\frac{G''(\gamma_{ij})}{\gamma_{ij}}-6\frac{F(\gamma_{ij})}{\gamma_{ij}}, & \quad \displaystyle C_2=3\frac{F''(\gamma_{ij})}{\gamma_{ij}}, \\[5pt]
& \displaystyle C_3=-3G'(\gamma_{ij})+C_1\frac{\gamma_{ij}^2}{2}, & \quad \displaystyle C_4=-9F'(\gamma_{ij})+\frac{3}{2}C_2\gamma_{ij}^2.
\end{array}
\]

Then we compute the functional $W(\boldsymbol{v}_{2})$ and obtain
\begin{equation}   \label{E:W(v2)}
\begin{array}{r l}
W(\boldsymbol{v}_{2})
= & \displaystyle
(U^i_2-U^j_2)^2\left(\boldsymbol{\mathcal{C}}_{2}^{ij}\delta^{-3/2}+\boldsymbol{\mathcal{C}}_{3}^{ij}\delta^{-1/2}\right)+\mathbf{C}_4\mu(U^i_2-U^j_2)^2,
\end{array}
\end{equation}
where
\begin{equation} \label{E:coeff2}
\boldsymbol{\mathcal{C}}_{2}^{ij}= \frac{3}{8}\pi\mu R^{3/2}d_{ij}^{-3/2}, \quad
\boldsymbol{\mathcal{C}}_{3}^{ij}=\frac{207}{320}\pi\mu R^{1/2}d_{ij}^{-1/2}.
\end{equation}
When we evaluate the difference between $W(\boldsymbol{v}_{2})$, defined by \eqref{E:W(v2)}, and
$W^*_2(\boldsymbol{\mathcal{S}}_{2})$, defined by \eqref{E:dual-fnl-1} we obtain the following estimate:
\begin{equation}   \label{E:difference-2}
\left|W(\boldsymbol{v}_{2})-W^*_2(\boldsymbol{\mathcal{S}}_{2})\right|\leq E_2(\boldsymbol{U}^i,\boldsymbol{U}^j), \quad \mbox{ as } \delta\rightarrow 0,
\end{equation}
where $E_2$ is defined by \eqref{E:error-in-2}.
Therefore,  $W(\boldsymbol{v}_{2})$ determines the asymptotics of $W(\boldsymbol{u}_{sq})$ as $\delta\rightarrow 0$:
\begin{equation}   \label{E:W(u2)}
W(\boldsymbol{u}_{sq}) = \left([\boldsymbol{U}^i-\boldsymbol{U}^j]\cdot\boldsymbol{q}^{ij}\right)^2\left(\boldsymbol{\mathcal{C}}_{2}^{ij}\delta^{-3/2}+\boldsymbol{\mathcal{C}}_{3}^{ij}\delta^{-1/2}\right)+E_2(\boldsymbol{U}^i,\boldsymbol{U}^j).
\end{equation}

\subsubsection{Verification of \eqref{E:error-in-3}}

Consider the minimizer $\displaystyle \boldsymbol{u}_{per}$ of $W$ over $V_{per}$.
The test function $\displaystyle \boldsymbol{v}_{3} \in V_{per}$ may be written in the form:
\begin{equation} \label{E:v3-decomp}
\boldsymbol{v}_{3}=R(\omega^i-\omega^j)\boldsymbol{v}^1_{3}+\beta_{ij}R\boldsymbol{v}^2_{3},
\end{equation}
where
\[
\left.\boldsymbol{v}^1_3\right|_{\partial B^i}=\boldsymbol{G}^i_3, \quad \left.\boldsymbol{v}^1_3\right|_{\partial B^j}=\boldsymbol{G}^j_3, \quad\quad\quad
\left.\boldsymbol{v}^2_3\right|_{\partial B^i}=\boldsymbol{0}, \quad \left.\boldsymbol{v}^2_3\right|_{\partial B^j}=\boldsymbol{0}.
\]
and the parameter $\beta_{ij}$ is defined as follows. The flux of $\boldsymbol{v}_{3}$ through $\ell_{ij}$ should satisfy
$\displaystyle \int_{\ell_{ij}}\boldsymbol{v}_{3}\cdot\boldsymbol{n}ds=\beta_{ij}^*-\frac{\delta_{ij}}{2R}(U^i_1+U^j_1)$, then
\begin{equation} \label{E:beta0}
\beta_{ij}=\beta_{ij}^*-\frac{\omega^i-\omega^j}{R}\int_{\ell_{ij}}\boldsymbol{v}^1_{3}\cdot\boldsymbol{n}ds -\frac{\delta_{ij}}{2R}(U^i_1+U^j_1),
\end{equation}
where $\displaystyle  \frac{\delta_{ij}}{2R}(U^i_1+U^j_1)=\int_{\ell_{ij}}\boldsymbol{u}_{t}\cdot\boldsymbol{n}ds$.

Define $\boldsymbol{v}_3$ by \eqref{E:v3}-\eqref{E:P,Q,K,M-3}. Computing $\displaystyle \int_{\ell_{ij}}\boldsymbol{v}^1_{3}\cdot\boldsymbol{n}ds=\frac{\delta_{ij}}{2}+\frac{\delta_{ij}^2}{8R}$
we obtain that $\beta_{ij}$ is defined by \eqref{E:beta0}.

Now we rewrite \eqref{E:v3} as one vector:
\begin{equation}  \label{E:v3-new}
\begin{array}{l l}
\boldsymbol{v}_3 & \displaystyle =\begin{pmatrix} \displaystyle G(x)+3y^2F(x) \\ \displaystyle -yG'(x)-y^3F'(x) \end{pmatrix} \quad \mbox{ in } \Pi_{ij},
\end{array}
\end{equation}
where
\begin{equation}  \label{E:G,F-33}
G(x)=R(\omega^i-\omega^j)P(x)+\beta_{ij}RK(x), \quad  F(x)=R(\omega^i-\omega^j)Q(x)+\beta_{ij}RM(x),
\end{equation}
with $P$, $Q$, $K$ and $M$ defined by \eqref{E:P,Q,K,M-3}.

As before to choose a trial tensor $\boldsymbol{\mathcal{S}}_{3} \in F_{per}$ we look at the stress tensor
$\boldsymbol{\sigma}(\boldsymbol{v}_3)=2\mu\boldsymbol{D}(\boldsymbol{v}_3)-p_{v_3} $, where
an approximate pressure $p_{v_3}$ can be chosen as in \eqref{E:p(v2)} for $G$, $F$ given by
\eqref{E:G,F-33}.

Then for such functions $G$ and $F$ define the trial tensor $\boldsymbol{\mathcal{S}}_{3} \in F_{per}$ by \eqref{E:S3(1)} where constants $C_i$ as $\boldsymbol{\mathcal{S}}_{2}$.

Next we compute the functional $W(\boldsymbol{v}_3)$:
\begin{equation} \label{E:W(v3)}
\begin{array}{l l}
W(\boldsymbol{v}_3) & \displaystyle
= \beta_{ij}^2\left(\boldsymbol{\mathcal{C}}_{4}^{ij}\delta^{-5/2} +\boldsymbol{\mathcal{C}}_{5}^{ij}\delta^{-3/2} +\boldsymbol{\mathcal{C}}_{6}^{ij}\delta^{-1/2}\right)\\[7pt]
& \displaystyle + R(\omega^i-\omega^j)\beta_{ij}\left(\boldsymbol{\mathcal{C}}_{7}^{ij}\delta^{-3/2}+\boldsymbol{\mathcal{C}}_{8}^{ij}\delta^{-1/2}\right)\\[7pt]
& \displaystyle + R^2(\omega^i-\omega^j)^2 \boldsymbol{\mathcal{C}}_{9}^{ij}\delta^{-1/2} \\[5pt]
& \displaystyle +\textbf{C}_5R^2(\omega^i-\omega^j)^2+\textbf{C}_6R^2(\omega^i-\omega^j)\beta_{ij}+\textbf{C}_7R^2\beta_{ij}^2,
\end{array}
\end{equation}
where
\begin{equation} \label{E:coeff3}
\begin{array}{l l l l}
& \displaystyle \boldsymbol{\mathcal{C}}_{4}^{ij}= \frac{9}{4}\pi\mu R^{5/2}d_{ij}^{-5/2},
& \displaystyle \boldsymbol{\mathcal{C}}_{5}^{ij}=\frac{99}{160}\pi\mu R^{3/2}d_{ij}^{-3/2},
& \displaystyle \boldsymbol{\mathcal{C}}_{6}^{ij}= \frac{29241}{17920}\pi\mu R^{1/2}d_{ij}^{-1/2}\\[7pt]
& \displaystyle \boldsymbol{\mathcal{C}}_{7}^{ij}= -3\pi\mu R^{3/2}d_{ij}^{-3/2},
& \displaystyle \boldsymbol{\mathcal{C}}_{8}^{ij}=\frac{9}{40}\pi\mu R^{1/2}d_{ij}^{-1/2},
& \displaystyle \boldsymbol{\mathcal{C}}_{9}^{ij}= \frac{3}{2}\pi\mu R^{1/2}d_{ij}^{-1/2}.
\end{array}
\end{equation}

When we evaluate the difference between $W(\boldsymbol{v}_{3})$, defined by \eqref{E:W(v3)}-\eqref{E:coeff3}, and
$W^*_3(\boldsymbol{\mathcal{S}}_{3})$, defined by \eqref{E:dual-fnl-2} we obtain the following estimate:
\begin{equation}   \label{E:difference-3}
\left|W(\boldsymbol{v}_{3})-W^*_3(\boldsymbol{\mathcal{S}}_{3})\right|\leq E_3(\omega^i,\omega^j,\beta_{ij}), \quad \mbox{ as } \delta\rightarrow 0,
\end{equation}
where $E_3$ is defined by \eqref{E:error-in-3}.
Therefore, $W(\boldsymbol{v}_{3})$ determines the asymptotics of $W(\boldsymbol{u}_{per})$:
\begin{equation}   \label{E:W(u3)}
\begin{array}{l l}
W(\boldsymbol{u}_{per}) & \displaystyle
= \beta_{ij}^2\left(\boldsymbol{\mathcal{C}}_{4}^{ij}\delta^{-5/2} +\boldsymbol{\mathcal{C}}_{5}^{ij}\delta^{-3/2} +\boldsymbol{\mathcal{C}}_{6}^{ij}\delta^{-1/2}\right)\\[7pt]
& \displaystyle + R(\omega^i-\omega^j)\beta_{ij}\left(\boldsymbol{\mathcal{C}}_{7}^{ij}\delta^{-3/2}+\boldsymbol{\mathcal{C}}_{8}^{ij}\delta^{-1/2}\right)\\[7pt]
& \displaystyle + R^2(\omega^i-\omega^j)^2 \boldsymbol{\mathcal{C}}_{9}^{ij}\delta^{-1/2} +E_3(\omega^i,\omega^j,\beta_{ij}).
\end{array}
\end{equation}

\subsection{Construction trial fields in boundary neck $\Pi_{ij}$, ($i\in\mathbb{I}$, $j\in\mathcal{N}_i\cap\mathbb{B}$)}  \label{A:bound}

\begin{figure}[!ht]
  \centering
  \includegraphics[scale=1.05]{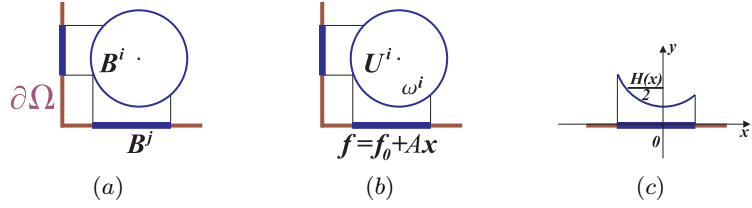}\\
  $(a)$ \hspace{85pt} $(b)$ \hspace{85pt} $(c)$
  \caption{$(a)$ The disk $B^i$ and quasidisk $B^j$, $(b)$ boundary conditions on $B^i$ and $B^j$, $(c)$ the local coordinate system}\label{F:quasi}
\end{figure}

The construction of the trial vector and tensor fields in neck $\Pi_{ij}$ connecting the disk $B^i$ and the quasidisk $B^j$ as in Fig. \ref{F:quasi} can be done similarly to the case of the neck connecting two interior disks.
Let $B^i$ be a disk near the boundary and $B^j$ a quasidisk, neighboring to $B_i$.
We call $\boldsymbol{U}^i,\omega^i$ of the quasidisk $B^j$ the ``translational'' and the ``angular" velocities of the quasidisk $B^j$, respectively.
Now consider the neck $\Pi_{ij}$ connecting the disk $B^i$ and the quasidisk $B^j$ lying on the lower part of $\partial\Omega$. The other cases are treated similarly.

Recall that in general necks are not symmetric. To construct a test function $\boldsymbol{v} \in V_{ij}$ in the neck $\Pi_{ij}$
connecting the disk $B^i$ and quasidisk $B^j$ (Fig. \ref{F:quasi}$(a)$) we extend $\Pi_{ij}$ to
\begin{equation}  \label{E:sym-neck}
\widehat{\Pi}_{ij}=\left\{(x,y) \in \mathbb{R}^2: \, -\gamma_{ij}<x < \gamma_{ij}, \,\,-H_{ij}(x)/2<y<H_{ij}(x)/2\right\},
\end{equation}
where
\[
\gamma_{ij}=
\begin{cases}
\gamma_{ij}^-, & \mbox{if } \gamma_{ij}^->\gamma_{ij}^+,\\
\gamma_{ij}^+, & \mbox{if } \gamma_{ij}^+>\gamma_{ij}^-,
\end{cases}
\]
and construct $\boldsymbol{\widehat{v}}\in H^1(\widehat{\Pi}_{ij})$ such that
$\boldsymbol{\widehat{v}}=\boldsymbol{U}^i+R\omega^i(n_1\boldsymbol{e}_2-n_2\boldsymbol{e}_1)$ on $\partial B^i \cap \partial \widehat{\Pi}_{ij}$,
$\boldsymbol{\widehat{v}}=A\boldsymbol{x}$ on $\partial B^j \cap \partial \widehat{\Pi}_{ij}$ and
$\displaystyle \frac{1}{R}\int_{\ell_{ij}}\boldsymbol{\widehat{v}}\cdot\boldsymbol{n}ds=\beta_{ij}^*$. Then the restriction of this function on $\Pi_{ij}$:
$\boldsymbol{v}=\boldsymbol{\widehat{v}}|_{\Pi_{ij}}$ would be a function from $V_{ij}$ (since this would be $H^1(\Pi_{ij})$ function and
all boundary conditions would be satisfied). Hence, $W_{\Pi_{ij}}(\boldsymbol{u})\leq W_{\Pi_{ij}}(\boldsymbol{v})$,
where $\boldsymbol{u}$ is the actual minimizer of $W_{\Pi_{ij}}$ in $V_{ij}$.

We can decompose a test function $\boldsymbol{v}\in V_{ij}$ ($i\in\mathbb{I},j\in\mathcal{N}_i\cap\mathbb{B}$) in to the sum
$\boldsymbol{v}=\boldsymbol{u}_t+\boldsymbol{v}_1+\boldsymbol{v}_2$ so that $W_{\Pi_{ij}}(\boldsymbol{v})=W_{\Pi_{ij}}(\boldsymbol{v}_1)+W_{\Pi_{ij}}(\boldsymbol{v}_2)$, and $W_{\Pi_{ij}}(\boldsymbol{v}_3)=0$ where
functions $\boldsymbol{v}_1$, $\boldsymbol{v}_2$ satisfy
\begin{equation} \label{E:boundary-12}
\begin{array}{r l l l l}
(a) & \boldsymbol{u}_t=\frac{1}{2}(\boldsymbol{U}^i+\boldsymbol{U}^j) \mbox{ for } \boldsymbol{x} \in \partial B^i \mbox{ and } \boldsymbol{x} \in  \partial B^j,\\[10pt]
(b) & \boldsymbol{v}_1=\begin{pmatrix} U^i_1 \\ 0 \end{pmatrix}+\omega^i\begin{pmatrix} \sqrt{R^2-x^2} \\ x \end{pmatrix}
\mbox{for } \boldsymbol{x}\in \partial B^i, \boldsymbol{v}_1=\begin{pmatrix} U_1^j \\ \omega^j x \end{pmatrix}  \mbox{for } \boldsymbol{x}\in \partial B^j,\\[10pt]
&\mbox{with } \displaystyle \int_{\ell_{ij}}\boldsymbol{v}_1\cdot\boldsymbol{n}ds=\beta_{ij}^*R-\int_{\ell_{ij}}\boldsymbol{u}_t\cdot\boldsymbol{n}ds
= \beta_{ij}^*R-\frac{1}{4}[(\boldsymbol{U}^i+\boldsymbol{U}^j)\cdot\boldsymbol{p}^{ij}]\delta_{ij},\\[5pt]
(c) & \boldsymbol{v}_2=\begin{pmatrix} 0 \\ U_2^i \end{pmatrix} \mbox{for } \boldsymbol{x}\in \partial B^i,
\boldsymbol{v}_2=\begin{pmatrix} ax \\ U_2^j \end{pmatrix} \mbox{for } \boldsymbol{x}\in \partial B^j.
\end{array}
\end{equation}
Then the corresponding decomposition of the minimizer is: $\boldsymbol{u}=\boldsymbol{u}_t+\boldsymbol{u}_1+\boldsymbol{u}_2$.

In the local coordinate system we choose these functions $\boldsymbol{v}_1$, $\boldsymbol{v}_2$ to be as follows:
\begin{equation} \label{E:v1-bound}
\displaystyle \boldsymbol{v}_1=\begin{pmatrix} U_1^j+yK(x)+y^2M(x)  \\  \omega^j x-\frac{y^2}{2}K'(x)-\frac{y^3}{3}M'(x) \end{pmatrix} \quad \mbox{in } \Pi_{ij},
\end{equation}
where
\[
\begin{array}{l l}
\displaystyle K(x)=-\frac{12(\omega^i-\omega^j)x^2}{H(x)^2}+\frac{8(U_1^i-U_1^j+R\omega^i)}{H(x)}-\frac{8\omega^i}{H(x)}(R-\sqrt{R^2-x^2})+\frac{24R\beta_{ij}}{H(x)^2}, \quad \\[10pt]
\displaystyle M(x)=\frac{24(\omega^i-\omega^j)x^2}{H(x)^3}-\frac{12(U_1^i-U_1^j+R\omega^i)}{H(x)^2}+\frac{12\omega^i}{H(x)^2}(R-\sqrt{R^2-x^2})-\frac{48R\beta_{ij}}{H(x)^3},
\end{array}
\]
with
\[
\beta_{ij}= \beta_{ij}^*-\frac{3\delta_{ij}}{4R}U_1^i-\frac{\delta_{ij}}{4R}U_1^j-\frac{\delta_{ij}}{2}\omega^i
= \beta_{ij}^*-\frac{3\delta_{ij}}{4R}\boldsymbol{U}^i\cdot\boldsymbol{p}^{ij}-\frac{\delta_{ij}}{4R}\boldsymbol{U}^j\cdot\boldsymbol{p}^{ij}-\frac{\delta_{ij}}{2}\omega^i,
\]
and
\begin{equation} \label{E:v2-bound}
\boldsymbol{v}_2=\begin{pmatrix} ax+yF(x)+y^2G(x)  \\  U_2^j-ay-\frac{y^2}{2}F'(x)-\frac{y^3}{3}G'(x) \end{pmatrix} \quad \mbox{in } \Pi_{ij},
\end{equation}
where
\[
F(x)=-\frac{8ax}{H(x)}-\frac{24(U_2^i-U_2^j)x}{H(x)^2}, \quad G(x)=\frac{12ax}{H(x)^2}+\frac{48(U_2^i-U_2^j)x}{H(x)^3}.
\]

We construct a lower bound for the viscous dissipation rate in the neck near the boundary we consider a dual variational problem in a symmetric neck
$\widehat{\Pi}_{ij}$ \eqref{E:sym-neck} with
\[
\gamma_{ij}=
\begin{cases}
\gamma_{ij}^-, & \mbox{if } \gamma_{ij}^-<\gamma_{ij}^+,\\
\gamma_{ij}^+, & \mbox{if } \gamma_{ij}^+<\gamma_{ij}^-,
\end{cases}
\]
that is, the problem of the maximization of the dual functional \eqref{E:dual-fnl-1} with $\boldsymbol{G}_3^i$ and $\boldsymbol{G}_3^j$ given by (\ref{E:boundary-12}a)
in the dual space $F_{per}$ given by \eqref{E:dual-set-2}, all defined over the neck $\widehat{\Pi}_{ij}$. Indeed,
\[
W_{\Pi_{ij}}(\boldsymbol{u})\geq W_{\widehat{\Pi}_{ij}}(\boldsymbol{u})\geq W^*_{\widehat{\Pi}_{ij}}(\boldsymbol{\mathcal{S}}).
\]
Here the test tensor $\boldsymbol{\mathcal{S}}$ is decomposed into the sum $\boldsymbol{\mathcal{S}}=\boldsymbol{\mathcal{S}}_1+\boldsymbol{\mathcal{S}}_2$,
with $\boldsymbol{\mathcal{S}}_1 \in F_{per}$, $\boldsymbol{\mathcal{S}}_2 \in F$ corresponding
to two trial vector fields $\boldsymbol{v}_1$, $\boldsymbol{v}_2$ defined above.

We choose
\begin{equation} \label{E:S1-bound}
\boldsymbol{\mathcal{S}}_1=\begin{pmatrix} -2\int_{\gamma_0}^x M(t)dt+C_1x & & K+2yM-C_1y-C_2  \\  K+2yM-C_1y-C_2 & & -yK'-y^2M'\end{pmatrix} \quad \mbox{in } \widehat{\Pi}_{ij},
\end{equation}
where constants $C_1=2M(\gamma_{ij})$, $C_2=K(\gamma_{ij})$. This tensor is divergence free in $\widehat{\Pi}_{ij}$, and on the lateral boundary $\partial \widehat{\Pi}_{ij}^\pm$
of the neck $\widehat{\Pi}_{ij}$ it takes values:
\[
\chi^+=-2\int_{\gamma_0}^{\gamma_{ij}}M(x)dx+C_1\gamma_{ij}, \quad \chi^+=-2\int_{\gamma_0}^{-\gamma_{ij}}M(t)dt-C_1\gamma_{ij},
\]
hence $\boldsymbol{\mathcal{S}}_1 \in F_{per}$.

Now we compute the functional $W(\boldsymbol{v}_{1})$ in the neck connecting a disk and a quasidisk:
\begin{equation*}
\begin{array}{r l}
W(\boldsymbol{v}_{1})
=  & \displaystyle \beta_{ij}^2\left[18\pi\mu R^{5/2}\delta_{ij}^{-5/2}+\frac{51}{20}\pi\mu R^{3/2}\delta_{ij}^{-3/2}+\frac{20889}{2240}\pi\mu R^{1/2}\delta_{ij}^{-1/2}\right] \\[7pt]
+  & \displaystyle 6\pi\mu R^{3/2}\beta_{ij}\left(U^i_1-U^j_1+R\omega^i\right)\delta_{ij}^{-3/2}-3\pi\mu R^{3/2}\beta_{ij} R(\omega^i-\omega^j)\delta_{ij}^{-3/2}\\[7pt]
+  & \displaystyle \frac{19}{20}\pi\mu R^{1/2}\beta_{ij} \left(U^i_1-U^j_1+R\omega^i\right)\delta_{ij}^{-1/2}-\frac{3}{8}\pi\mu R^{1/2}\beta_{ij} R(\omega^i-\omega^j)\delta_{ij}^{-1/2}\\[7pt]
-  & \displaystyle 3\pi\mu R^{1/2}\beta_{ij} R\omega^i\delta_{ij}^{-1/2}+4\pi\mu R^{1/2}\left(U^i_1-U^j_1+R\omega^i\right)^2\delta_{ij}^{-1/2}\\[7pt]
-  & \displaystyle 3\pi\mu R^{1/2}\left(U^i_1-U^j_1+R\omega^i\right)R(\omega^i-\omega^j)\delta_{ij}^{-1/2}\\[7pt]
+  & \displaystyle \frac{9}{2}\pi\mu R^{1/2}R^2(\omega^i-\omega^j)^2\delta_{ij}^{-1/2}+E_1^{b}(\boldsymbol{U}^i,\omega^i,\boldsymbol{f}, \beta_{ij})
\end{array}
\end{equation*}
\begin{equation}   \label{E:W(v1)-bound}
\begin{array}{r l}
=:  & \displaystyle \beta_{ij}^2\left[\boldsymbol{\mathcal{C}}_{1}^{ij}\delta^{-5/2}+\boldsymbol{\mathcal{C}}_{2}^{ij}\delta^{-3/2}+\boldsymbol{\mathcal{C}}_{3}^{ij}\delta^{-1/2}\right] \\[7pt]
+  & \displaystyle \boldsymbol{\mathcal{C}}_{4}^{ij}\beta_{ij}\left(U^i_1-U^j_1+R\omega^i\right)\delta^{-3/2}+\boldsymbol{\mathcal{C}}_{5}^{ij}\beta_{ij} R(\omega^i-\omega^j)\delta^{-3/2}\\[7pt]
+  & \displaystyle \boldsymbol{\mathcal{C}}_{6}^{ij}\beta_{ij}\left(U^i_1-U^j_1+R\omega^i\right)\delta^{-1/2}+\boldsymbol{\mathcal{C}}_{7}^{ij}\beta_{ij} R(\omega^i-\omega^j)\delta^{-1/2}\\[7pt]
+  & \displaystyle \boldsymbol{\mathcal{C}}_{8}^{ij}\beta_{ij}R\omega^i\delta^{-1/2}+\boldsymbol{\mathcal{C}}_{9}^{ij}\left(U^i_1-U^j_1+R\omega^i\right)^2\delta^{-1/2}\\[7pt]
+  & \displaystyle \boldsymbol{\mathcal{C}}_{10}^{ij}\left(U^i_1-U^j_1+R\omega^i\right)R(\omega^i-\omega^j)\delta^{-1/2}\\[7pt]
+  & \displaystyle \boldsymbol{\mathcal{C}}_{11}^{ij}R^2(\omega^i-\omega^j)^2\delta^{-1/2}+E_1^{b}(\boldsymbol{U}^i,\omega^i,\boldsymbol{f}, \beta_{ij}),
\end{array}
\end{equation}
where the error term
\begin{equation} \label{E:error-b-1}
\begin{array}{r l}
E_1^{b}(\boldsymbol{U}^i,\omega^i,\boldsymbol{f}, \beta_{ij})= & \displaystyle \mu \left(C_7\beta_{ij}^2+
C_8[(\boldsymbol{U}^i-\boldsymbol{U}^j)\cdot\boldsymbol{p}^{ij}]^2+C_{9}R^2(\omega^i)^2+C_{10}R^2(\omega^j)^2\right)
\end{array}
\end{equation}
with universal constants $C_k$, $k=7,\ldots,10$ and
\[
\begin{array}{r l l l}
& \displaystyle \boldsymbol{\mathcal{C}}_{1}^{ij} = 18\pi\mu R^{5/2}d_{ij}^{-5/2},
& \displaystyle \boldsymbol{\mathcal{C}}_{2}^{ij} = \frac{51}{20}\pi\mu R^{3/2}d_{ij}^{-3/2},
& \displaystyle \boldsymbol{\mathcal{C}}_{3}^{ij} = \frac{20889}{2240}\pi\mu R^{1/2}d_{ij}^{-1/2},\\[7pt]
& \displaystyle \boldsymbol{\mathcal{C}}_{4}^{ij} = 6\pi\mu R^{3/2}d_{ij}^{-3/2},
& \displaystyle \boldsymbol{\mathcal{C}}_{5}^{ij} = - 3\pi\mu R^{3/2}d_{ij}^{-3/2},
& \displaystyle \boldsymbol{\mathcal{C}}_{6}^{ij} = \frac{19}{20}\pi\mu R^{1/2}d_{ij}^{-1/2},\\[7pt]
& \displaystyle \boldsymbol{\mathcal{C}}_{7}^{ij} = -\frac{3}{8}\pi\mu R^{1/2}d_{ij}^{-1/2},
& \displaystyle \boldsymbol{\mathcal{C}}_{8}^{ij} = -3\pi\mu R^{1/2}\delta_{ij}^{-1/2},
& \displaystyle \boldsymbol{\mathcal{C}}_{9}^{ij} = 4\pi\mu R^{1/2}d_{ij}^{-1/2},\\[7pt]
& \displaystyle \boldsymbol{\mathcal{C}}_{10}^{ij} = -3\pi\mu R^{1/2}d_{ij}^{-1/2},
& \displaystyle \boldsymbol{\mathcal{C}}_{11}^{ij} = \frac{9}{2}\pi\mu R^{1/2}d_{ij}^{-1/2}.
\end{array}
\]

When we evaluate the difference between $W(\boldsymbol{v}_{1})$, defined by \eqref{E:W(v1)-bound}, and
$W^*_3(\boldsymbol{\mathcal{S}}_{1})$, defined by \eqref{E:dual-fnl-2} with the functions $\boldsymbol{G}_k^i$, $\boldsymbol{G}_k^j$ given by (\ref{E:boundary-12}a) we obtain the following estimate:
\begin{equation}   \label{E:difference-1-bound}
\left|W(\boldsymbol{v}_{1})-W^*_2(\boldsymbol{\mathcal{S}}_{1})\right|\leq E_1^{b}(\boldsymbol{U}^i,\omega^i,\boldsymbol{f}, \beta_{ij}), \quad \mbox{ as } \delta\rightarrow 0.
\end{equation}
Therefore,  $W(\boldsymbol{v}_{1})$ determines the asymptotics of $W(\boldsymbol{u}_{1})$ in the neck connecting the disk and quasidisk:
\begin{equation}   \label{E:W(u1)-bound}
\begin{array}{r l}
W(\boldsymbol{u}_{1})
=  & \displaystyle \beta_{ij}^2\left[\boldsymbol{\mathcal{B}}_{1}^{ij}\delta^{-5/2}+\boldsymbol{\mathcal{B}}_{2}^{ij}\delta^{-3/2}+\boldsymbol{\mathcal{B}}_{3}^{ij}\delta^{-1/2}\right] \\[7pt]
+  & \displaystyle \boldsymbol{\mathcal{B}}_{4}^{ij}\beta_{ij}\left[(\boldsymbol{U}^i-\boldsymbol{U}^j)\cdot\boldsymbol{p}^{ij}+R\omega^i\right]\delta^{-3/2}+\boldsymbol{\mathcal{B}}_{5}^{ij}\beta_{ij} R(\omega^i-\omega^j)\delta^{-3/2}\\[7pt]
+  & \displaystyle \boldsymbol{\mathcal{B}}_{6}^{ij}\beta_{ij}\left[(\boldsymbol{U}^i-\boldsymbol{U}^j)\cdot\boldsymbol{p}^{ij}+R\omega^i\right]\delta^{-1/2}+\boldsymbol{\mathcal{B}}_{7}^{ij}\beta_{ij} R(\omega^i-\omega^j)\delta^{-1/2}\\[7pt]
+  & \displaystyle \boldsymbol{\mathcal{B}}_{8}^{ij}\beta_{ij}R\omega^i\delta^{-1/2}+\boldsymbol{\mathcal{B}}_{9}^{ij}\left[(\boldsymbol{U}^i-\boldsymbol{U}^j)\cdot\boldsymbol{p}^{ij}+R\omega^i\right]^2\delta^{-1/2}\\[7pt]
+  & \displaystyle \boldsymbol{\mathcal{B}}_{10}^{ij}\left[(\boldsymbol{U}^i-\boldsymbol{U}^j)\cdot\boldsymbol{p}^{ij}+R\omega^i\right]R(\omega^i-\omega^j)\delta^{-1/2}\\[7pt]
+  & \displaystyle \boldsymbol{\mathcal{B}}_{11}^{ij}R^2(\omega^i-\omega^j)^2\delta^{-1/2}+E_1^{b}(\boldsymbol{U}^i,\omega^i,\boldsymbol{f}, \beta_{ij}).
\end{array}
\end{equation}

Now we choose $\boldsymbol{\mathcal{S}}_{2} \in F$ as follows:
\begin{equation}   \label{E:S2(1)-bound}
\begin{array}{r l l }
& \displaystyle (\boldsymbol{\mathcal{S}}_{2})_{11} = & \displaystyle \mu \left(3yF'+3y^2G'-2\int_{\gamma_0}^x G+\frac{x^2}{2}\left(C_2+C_3y+C_4y^2\right)\right.\\[5pt]
& & \displaystyle \left.-(C_5+C_6y+C_7y^2)\right),\\[5pt]
& \displaystyle (\boldsymbol{\mathcal{S}}_{2})_{12} = & \displaystyle (\boldsymbol{\mathcal{S}}_{2})_{21} =
\mu \left(F+2yG-\frac{3}{2}y^2F''-y^3G''-C_1x-C_2xy-C_3x\frac{y^2}{2}-C_4x\frac{y^3}{3}\right),\\[5pt]
& \displaystyle (\boldsymbol{\mathcal{S}}_{2})_{22} = & \displaystyle \mu\left(-yF'-y^2G'-2\int_{\gamma_0}^x G+\frac{1}{2}y^3F'''+\frac{1}{4}y^4G'''\right.\\[5pt]
& & \displaystyle \left.+C_1y+C_2\frac{y^2}{2}+C_3\frac{y^3}{6}+C_4\frac{y^4}{12}\right),
\end{array}
\end{equation}
where constants $C_i$ ($i=1..7$) are chosen so that $\boldsymbol{\mathcal{S}}_{2}\boldsymbol{n}=\boldsymbol{0}$ on $\partial \Pi_{ij}^\pm$, that is,
\[
\begin{array}{l l}
& \displaystyle C_1=\frac{1}{\gamma_{ij}}F(\gamma_{ij}), \quad \displaystyle C_2=\frac{2}{\gamma_{ij}}G(\gamma_{ij}), \quad
\displaystyle C_3=-\frac{3}{\gamma_{ij}}F''(\gamma_{ij}), \quad \displaystyle C_4=-\frac{3}{\gamma_{ij}}G''(\gamma_{ij}), \\[7pt]
& \displaystyle C_5=C_2\frac{\gamma_{ij}^2}{2}-2\int_{\gamma_0}^{\gamma_{ij}} G, \quad \displaystyle C_6=C_3\frac{\gamma_{ij}^2}{2}+3F'(\gamma_{ij}),
\quad \displaystyle C_7=C_4\frac{\gamma_{ij}^2}{2}+3G'(\gamma_{ij}).
\end{array}
\]

We now compute the functional $W(\boldsymbol{v}_{2})$:
\begin{equation}   \label{E:W(v2)-bound}
\begin{array}{r l}
W(\boldsymbol{v}_{2})
= & \displaystyle (U^i_2-U^j_2)^2\left[6\pi\mu R^{3/2}\delta_{ij}^{-3/2}+\frac{63}{20}\pi\mu R^{1/2}\delta_{ij}^{-1/2}\right]\\[7pt]
+ & \displaystyle 6\pi\mu R^{1/2}(U^i_2-U^j_2)Ra\delta_{ij}^{-1/2}+E_2^{b}(\boldsymbol{U}^i,\boldsymbol{f})\\[7pt]
=: & \displaystyle (U^i_2-U^j_2)^2\left(\boldsymbol{\mathcal{C}}_{12}^{ij}\delta^{-3/2}+\boldsymbol{\mathcal{C}}_{13}^{ij}\delta^{-1/2}\right)\\[7pt]
+ & \displaystyle
\boldsymbol{\mathcal{C}}_{14}^{ij}(U^i_2-U^j_2)Ra\delta^{-1/2}+E_2^{b}(\boldsymbol{U}^i,\boldsymbol{f}).
\end{array}
\end{equation}
where the error term
\begin{equation} \label{E:error-b-2}
E_2^{b}(\boldsymbol{U}^i,\boldsymbol{f})=\mu C_{11}[(\boldsymbol{U}^i-\boldsymbol{U}^j)\cdot\boldsymbol{q}^{ij}]^2,
\end{equation}
with universal constant $C_{11}$ and
\[
\begin{array}{r l l l}
& \displaystyle \boldsymbol{\mathcal{C}}_{12}^{ij} = 6\pi\mu R^{3/2}d_{ij}^{-3/2},
& \displaystyle \boldsymbol{\mathcal{C}}_{13}^{ij} = \frac{63}{20}\pi\mu R^{1/2}d_{ij}^{-1/2},
& \displaystyle \boldsymbol{\mathcal{C}}_{14}^{ij} = 6\pi\mu R^{1/2}d_{ij}^{-1/2}.
\end{array}
\]

When we evaluate the difference between $W(\boldsymbol{v}_{2})$, defined by \eqref{E:W(v2)-bound}, and
$W^*_2(\boldsymbol{\mathcal{S}}_{2})$, defined by \eqref{E:dual-fnl-1} we obtain the following estimate:
\begin{equation}   \label{E:difference-2-bound}
\left|W(\boldsymbol{v}_{2})-W^*_2(\boldsymbol{\mathcal{S}}_{2})\right|\leq E_2^{b}(\boldsymbol{U}^i,\boldsymbol{f}), \quad \mbox{ as } \delta\rightarrow 0.
\end{equation}
Therefore,  $W(\boldsymbol{v}_{2})$ determines the asymptotics of $W(\boldsymbol{u}_{2})$ in the neck connecting the disk and quasidisk:
\begin{equation}   \label{E:W(u2)-bound}
\begin{array}{r l}
W(\boldsymbol{u}_{2})
= & \displaystyle \left[(\boldsymbol{U}^i-\boldsymbol{U}^j)\cdot\boldsymbol{q}^{ij}\right]^2\left(\boldsymbol{\mathcal{B}}_{12}^{ij}\delta^{-3/2}+\boldsymbol{\mathcal{B}}_{13}^{ij}\delta^{-1/2}\right)\\[7pt]
+ & \displaystyle \boldsymbol{\mathcal{B}}_{14}^{ij}\left[(\boldsymbol{U}^i-\boldsymbol{U}^j)\cdot\boldsymbol{q}^{ij}\right]Ra\delta^{-1/2}+E_2^{b}(\boldsymbol{U}^i,\boldsymbol{f}).
\end{array}
\end{equation}

\Section{Appendix}

\subsection{Positive definiteness of $Q_{ij}$ in \eqref{E:Q_ij}} \label{A:Q}

\begin{lemma} \label{L:quadr_form}
The quadratic form $Q_{ij}$, defined by \eqref{E:Q_ij}, \eqref{E:coefficients}, is positive definite.
\end{lemma}

We show the positive definiteness of $Q_{ij}$ with coefficients \eqref{E:coefficients} for $i\in\mathbb{I}$, $j\in\mathcal{N}_i\cap\mathbb{I}$ (the case when $j\in\mathcal{N}_i\cap\mathbb{B}$ can be treated similarly).
To this end we just need to show that $Q_{ij}(\boldsymbol{U}^i,\boldsymbol{U}^j,\omega^i,\omega^j,\beta_{ij})$ can be written as a sum of squares. Indeed,
\[
\begin{array}{r l l}
\displaystyle
Q(\boldsymbol{U}^i,\boldsymbol{U}^j,\omega^i,\omega^j,\beta_{ij})=
\frac{\pi\mu\sqrt{R}}{\sqrt{d_{ij}}\sqrt{\delta}} &
\displaystyle \left\{\frac{1}{2}\left[(\boldsymbol{U}^i-\boldsymbol{U}^j)\cdot\boldsymbol{p}^{ij}+ R\omega^i+R\omega^j\right]^2\right.\\[7pt]
+& \displaystyle \left[(\boldsymbol{U}^i-\boldsymbol{U}^j)\cdot\boldsymbol{q}^{ij}\right]^2\left(\frac{3}{8}R\delta^{-1}+\frac{207}{320}\right)\\[7pt]
+& \displaystyle \beta_{ij}^2\left(2\delta^{-2}+\frac{249}{160 R}\delta^{-1}+\frac{13491}{17920 R^2} \right) \\[7pt]
+& \displaystyle \left. \left[\frac{3}{4}R(\omega^i-\omega^j)-\beta_{ij}\left(\frac{1}{2}\delta^{-1}-\frac{15}{16 R}\right)\right]^2\right\} \\[5pt]
+& \displaystyle E_{ij}(\boldsymbol{U}^i,\boldsymbol{U}^j,\omega^i,\omega^j,\beta_{ij})\\[5pt]
= \displaystyle \frac{\pi\mu\sqrt{R}}{\sqrt{d_{ij}}\sqrt{\delta}} & \displaystyle \left\{\frac{1}{2}\left[(\boldsymbol{U}^i-\boldsymbol{U}^j)\cdot\boldsymbol{p}^{ij}+ R\omega^i+R\omega^j\right]^2\right.\\[7pt]
+& \displaystyle \left[(\boldsymbol{U}^i-\boldsymbol{U}^j)\cdot\boldsymbol{q}^{ij}\right]^2\left(\frac{3}{8}R\delta^{-1}+\frac{207}{320}\right)\\[7pt]
+& \displaystyle \beta_{ij}^2\left( \sqrt{2}\delta^{-1}+\frac{249}{320\sqrt{2}R} \right)^2+\frac{645273}{1433600 R^2}\beta_{ij}^2\\[7pt]
+& \displaystyle \left. \left[\frac{3}{4}R(\omega^i-\omega^j)-\beta_{ij}\left(\frac{1}{2}\delta^{-1}-\frac{15}{16R}\right)\right]^2\right\}\\[5pt]
+& \displaystyle E_{ij}(\boldsymbol{U}^i,\boldsymbol{U}^j,\omega^i,\omega^j,\beta_{ij})
\end{array}
\]

\subsection{Existence and uniqueness of the solution to \eqref{E:v_form-var}} \label{A:ex-un-1}

\begin{lemma} \label{L:exist-uniq1}
There exists a unique minimizer of $W_{\Omega_F}(\cdot)$ over $V$.
\end{lemma}

\begin{proof} First of all we note that $V$ is not empty.
Indeed, let $\boldsymbol{w}$ be a solution to Stokes equation
satisfying the Dirichlet boundary conditions, that is, of
\begin{equation}   \label{E:Dir-pr}
\left\{
\begin{array}{r l l}
(a) & \displaystyle \mu\triangle \boldsymbol{u}=\nabla p,
& \boldsymbol{x} \in \Omega_F \\[5pt]
(b) & \displaystyle \nabla\cdot \boldsymbol{u}=0, &
 \boldsymbol{x} \in \Omega_F\\[5pt]
(c) & \displaystyle
\boldsymbol{u}=\boldsymbol{U}^{i}_0+R\omega^i_0(n_1^i\boldsymbol{e}_2-n_2^i\boldsymbol{e}_1),
& \boldsymbol{x} \in \partial B^{i}, \quad i=1,\ldots, N \\[5pt]
(d) & \displaystyle \boldsymbol{u}=\boldsymbol{f}, &
\boldsymbol{x} \in
\partial\Omega
\end{array}
\right.
\end{equation}
with some given $\boldsymbol{U}^{i}_0$ and $\omega^i_0$,
$i=1,\ldots, N$, which is known to exist. Then $\boldsymbol{w} \in
V$.

The existence of the minimizer $\boldsymbol{u}$ would be shown in
few steps.

\textit{Step 1.} Set $\displaystyle
m:=\inf_{V}W_{\Omega_F}(\cdot)$. If $m=+\infty$ then we are done.
So we henceforth assume that $m$ is finite. Select a minimizing
sequence $\{\boldsymbol{u}_k\}_{k=1}^{\infty}$ of the functional
$W_{\Omega_F}(\cdot)$, that is,
$W_{\Omega_F}(\boldsymbol{u}_k)\rightarrow m$ as
$k\rightarrow\infty$.

\textit{Step 2.} Choose a sequence
$\{\boldsymbol{\bar{u}}_k\}_{k=1}^{\infty}$, where each
$\boldsymbol{\bar{u}}_k$ solves a Stokes equation \eqref{E:Dir-pr}
with boundary conditions (\ref{E:Dir-pr}$c$) on $\partial B^i$
which are the same as the corresponding boundary conditions of
$\boldsymbol{u}_k$. Hence,
\[
m\leq W_{\Omega_F}(\boldsymbol{\bar{u}}_k) \leq
W_{\Omega_F}(\boldsymbol{u}_k)\rightarrow m \quad \mbox{as} \quad
k\rightarrow\infty.
\]

Now we claim that there exists a convergent subsequence
$\{\boldsymbol{\varphi}^i_{kj}\}$ of values of velocities on each
disk, that is, there exists $\boldsymbol{\varphi}^i$ such that
$\boldsymbol{\varphi}^i_{kj}\rightarrow \boldsymbol{\varphi}^i$ as
$j\rightarrow\infty$ on each $\partial B^i$, $i=1,\ldots, N$

We prove this assertion by contradiction. Assume there exists such
a disk $B^{i_0}$ and a subsequence
$\{\boldsymbol{\varphi}^{i_0}_{kj}\}=\{\boldsymbol{U}^{i_0}_{kj}+\boldsymbol{\Omega}^{i_0}_{kj}\}$
on it that $\boldsymbol{\varphi}^{i_0}_{kj}\rightarrow\infty$ as
$j\rightarrow\infty$, where $\boldsymbol{U}^{i_0}_{kj}$ is the
translational velocity and $\boldsymbol{\Omega}^{i_0}_{kj}$ is the
rotational part. We select a subsequence
$\{\boldsymbol{\bar{u}}_{kj}\} \in V$ of the above sequence
$\{\boldsymbol{\bar{u}}_{k}\}$, that is,
$\{\boldsymbol{\bar{u}}_{kj}\}$ solves the Stokes problem with
fixed values $\boldsymbol{\varphi}^{i}_{kj}$ on the disks $B^i$.
Since $W_{\Omega_F}(\boldsymbol{\bar{u}}_{kj})$ is bounded and
$\boldsymbol{\varphi}^{i_0}_{kj}\rightarrow\infty$ we have
\[
W_{\Omega_F}\left(\frac{\boldsymbol{\bar{u}}_{kj}}{\boldsymbol{\varphi}^{i_0}_{kj}}\right)=
\frac{1}{(|\boldsymbol{\varphi}^{i_0}_{kj}|)^2}W_{\Omega_F}(\boldsymbol{\bar{u}}_{kj})\leq
\varepsilon_{kj}^{i_0},
\]
where $\varepsilon_{kj}^{i_0}$ is a small arbitrary number.

Without loss of generality, we assume that the disk $B^{i_0}$ lies
in the corner of the domain $\Omega$, that is, near two parts of
the boundary $\partial\Omega$ (see Fig. \ref{F:bound}), where
values of $\boldsymbol{\bar{u}}_{kj}$ are prescribed and equal to
$\boldsymbol{f}$.

\begin{figure}[h!]
  \centering
  \includegraphics[scale=0.95]{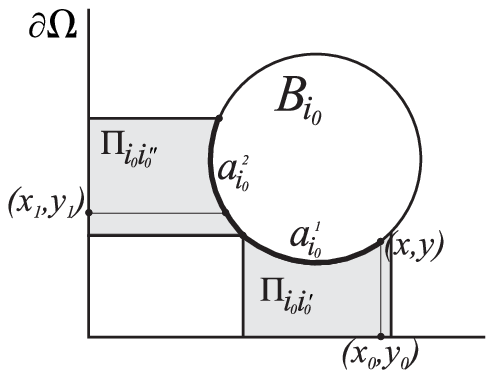}
  \caption{} \label{F:bound}
\end{figure}

First consider a neck $\Pi_{i_{0} i'_{0}}$ and choose a symmetric
neck $\widetilde{\Pi}_{i_{0} i'_{0}}$ containing $\Pi_{i_{0}
i'_{0}}$, that is, $\widetilde{\Pi}_{i_{0} i'_{0}}\supset
\Pi_{i_{0} i'_{0}}$. For simplicity we drop indices $kj$
indicating the subsequence and $i_0$ indicating the disk number.
Then denote
$\boldsymbol{\varphi}=(\varphi_1,\varphi_2)=(\boldsymbol{U},\boldsymbol{\Omega})=(U_1+\Omega_1,U_2+\Omega_2)$
and $\displaystyle
\frac{\boldsymbol{\bar{u}}}{|\boldsymbol{\varphi}|}=(\bar{u},\bar{v})$.
Choose a point $(x_0,y_0)$ on the lower part of $\partial\Omega$
and $(x,y)$ on $\partial B^{i_0}$ (Fig. \ref{F:bound}). By
integrating along the line connecting $x$ and $x_0$ we obtain:
\[
\bar{v}(x,y)=\frac{f_2(x_0,y_0)}{|\boldsymbol{\varphi}|}+\int_{y_0}^{y}\frac{\partial
\bar{v}}{\partial y}(x,y)dy,
\]
or
\[
|\bar{v}(x,y)|\leq \eta_2+\int_{0}^{H(\gamma_{i_0
i'_0})/2}\left|\frac{\partial \bar{v}}{\partial y}(x,y)\right|dy,
\]
where $\bar{v}(x,y)=\varphi_{2}$ and $\eta_2$ is some small number
since $\frac{f_2(x_0,y_0)}{|\boldsymbol{\varphi}|}\rightarrow 0$
as $j\rightarrow\infty$ and $H(\gamma_{i_0 i'_0})/2$ is a
half-hight of the neck $\widetilde{\Pi}_{i_0 i'_0}$. Then
integrate over the arc $a_{i_0}^1$ of $\partial B^{i_0}$ which is
an upper boundary boundary of the neck $\Pi_{i_0 i'_0}$:
\[
\left|\int_{a_{i_0}^1} \bar{v}ds\right|\leq \int_{a_{i_0}^1}\left|
\bar{v}\right|ds \leq |a_{i_0}^1|\eta_2+\int_{\Pi_{i_0
i'_0}}\left|\frac{\partial \bar{v}}{\partial
y}(x,y)\right|d\boldsymbol{x}
\]
and
\[
\left|\int_{a_{i_0}^1} \bar{v}ds\right|^2\leq
C_1\eta_2+C_2\int_{\Pi_{i_0 i'_0}}\left|\frac{\partial
\bar{v}}{\partial y}(x,y)\right|^2d\boldsymbol{x}\leq
C_1\eta_2+C_2W_{\Omega_F}(\boldsymbol{\bar{u}})\leq C_1\eta_2+C_2
\varepsilon.
\]

\begin{figure}[h!]
  \centering
  \includegraphics[scale=0.95]{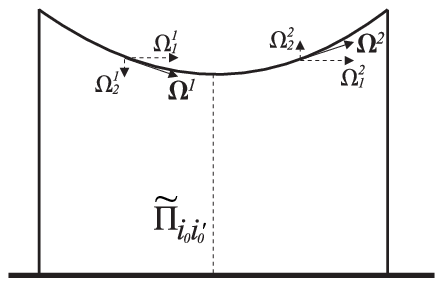}
  \caption{} \label{F:bound1}
\end{figure}

Consider
\[
\int_{a_{i_0}^1} \bar{v}ds=\int_{a_{i_0}^1} (U_2+\Omega_2)ds=
|a_{i_0}^1|U_2+\int_{a_{i_0}^1}\Omega_2 ds.
\]
Note that $\int_{a_{i_0}^1}\Omega_2 ds=0$ since the second
components of the vector $\boldsymbol{\Omega}$ cancel each other
when we write an integral as a sum over infinitesimal elements of
the symmetric arc $a_{i_0}^1$(see Fig. \ref{F:bound1}). Hence,
\[
C_3U_2^2\leq C_1\eta_2+C_2 \varepsilon,
\]
where $\eta_2=(\eta_2)^{i_0}_{kj}$ and
$\varepsilon=\varepsilon^{i_0}_{kj}$ are small arbitrary numbers.
So the second component of the vector $\boldsymbol{U}^{i_0}_{kj}$
tends to zero as $j\rightarrow\infty$.

Similarly, one can show that the first component of the vector
$\boldsymbol{U}^{i_0}_{kj}$ tends to zero as $j\rightarrow\infty$
considering the neck $\Pi_{i_0 i''_0}$.

Thus, we prove that $\boldsymbol{U}^{i_0}_{kj}\rightarrow 0$ as
$j\rightarrow\infty$.

To show that $\boldsymbol{\Omega}^{i_0}_{kj}\rightarrow 0$ as
$j\rightarrow\infty$ we consider only the left-hand side of the
neck $\Pi_{i_0 i'_0}$. Repeating all computations described above,
we obtain that the integral of $\Omega_2$ over the left half of
the symmetric arc $a_{i_0}^{1\ell}$ tends to zero as
$j\rightarrow\infty$:
\[
\int_{a_{i_0}^{1\ell}}\Omega_2ds \leq C_1\eta_2+C_2 \varepsilon,
\]
But $\Omega_2=\omega^{i_0}_{kj} \cos \theta$, where
$\omega^{i_0}_{kj}$ is a constant and $\theta$ is an angle of the
vector $\boldsymbol{\Omega}^{i_0}_{kj}$ with vertical axis. Thus,
$\theta>0$ everywhere in that half of the neck, hence
$\Omega_2\rightarrow 0$. Similarly, one can show that
$\Omega_1\rightarrow 0$, hence, the whole subsequence
$\boldsymbol{\varphi}^{i_0}_{kj}$ is convergent what contradicts
our assumption.

\begin{figure}[ht!]
  \centering
  \includegraphics[scale=0.95]{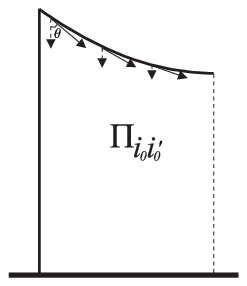}
  \caption{} \label{F:bound2}
\end{figure}

Note that if $B^{i_0}$ is not near $\partial\Omega$ one can use
the same argument using two nearest neighboring disks $B^{m_1}$
and $B^{m_2}$ of $B^{i_0}$ where there exist convergent
subsequences $\boldsymbol{\varphi}^{m_1}_{kj}$ and
$\boldsymbol{\varphi}^{m_2}_{kj}$, respectively. We can always
assume that since our graph is connected, that is there exist
paths connecting the given disk with all parts of the boundary
$\partial\Omega$.

\textit{Step 3.} Denote the limit of convergent subsequences on
each disk $B^i$ by $\boldsymbol{v}^i$, $i=1,\ldots,N$ and consider
the subsequence $\{\boldsymbol{\bar{u}}_{kj}\} \in V$ such that on
the boundary of each disk $\partial B^i$:
$\boldsymbol{\bar{u}}_{kj}\rightarrow \boldsymbol{v}^i$ as
$j\rightarrow\infty$.

\textit{Step 4.} Choose a sequence
$\{\boldsymbol{\bar{\bar{u}}}_{j}\} \in V$ such that each
$\boldsymbol{\bar{u}}_{j}$ solves the Stokes problem
\eqref{E:EL-boldneck} with
$\boldsymbol{\bar{\bar{u}}}_{j}=\boldsymbol{v}^i$ on $\partial
B^i$, $i=1,\ldots,N$.

\textit{Step 5.} Note that
\begin{equation} \label{E:step5}
W_{\Omega_F}(\boldsymbol{\bar{\bar{u}}}_{j})\leq
W_{\Omega_F}(\boldsymbol{\bar{u}}_{kj}) +
W_{\Omega_F}(\boldsymbol{\hat{u}}_j)+
2\sqrt{W_{\Omega_F}(\boldsymbol{\bar{u}}_{kj})W_{\Omega_F}(\boldsymbol{\hat{u}}_j)},
\end{equation}
where
$\boldsymbol{\hat{u}}_j=\boldsymbol{\bar{u}}_{kj}-\boldsymbol{\bar{\bar{u}}}_{j}$,
that is solves
\[
\left\{
\begin{array}{r l l}
(a) & \displaystyle \mu\triangle \boldsymbol{u}=\nabla p,
& \boldsymbol{x} \in \Omega_F \\[5pt]
(b) & \displaystyle \nabla\cdot \boldsymbol{u}=0, &
 \boldsymbol{x} \in \Omega_F\\[5pt]
(c) & \displaystyle
\boldsymbol{u}=\boldsymbol{v}^{i}_{kj}-\boldsymbol{v}^i,
& \boldsymbol{x} \in \partial B^{i}, \quad i=1,\ldots, N \\[5pt]
(d) & \displaystyle \boldsymbol{u}=\boldsymbol{0}, &
\boldsymbol{x} \in
\partial\Omega
\end{array}
\right.
\]
Hence, (see e.g. \cite{galdi})
$\|\boldsymbol{\hat{u}}_j\|_{H^1(\Omega_F)}\leq
C\|\boldsymbol{v}^{i}_{kj}-\boldsymbol{v}^i\|_{H^{1/2}(\partial
B^i)}\rightarrow 0$ as $j\rightarrow\infty$. Then from
\eqref{E:step5} we obtain
\[
W_{\Omega_F}(\boldsymbol{\bar{\bar{u}}}_{j})\rightarrow m \quad
\mbox{as} \quad j\rightarrow\infty.
\]

\textit{Step 6.} Choose
$\boldsymbol{\bar{\bar{u}}}=\boldsymbol{\bar{\bar{u}}}_{1}$ and
$\boldsymbol{u}^0_{j}=\boldsymbol{\bar{\bar{u}}}_{j}-\boldsymbol{\bar{\bar{u}}}$
where $\boldsymbol{u}^0_{j} \in H^1_0(\Omega_F)$ such that
$W_{\Omega_F}(\boldsymbol{u}^0_{j}
+\boldsymbol{\bar{\bar{u}}})\rightarrow m$ as
$j\rightarrow\infty$.

\textit{Step 7.} Consider
\begin{equation} \label{E:step7}
\begin{array}{l l}
\|\boldsymbol{\bar{\bar{u}}}_{j}\|_{L^2(\Omega_F)} & \displaystyle
\leq
\|\boldsymbol{\bar{\bar{u}}}_{j}-\boldsymbol{\bar{\bar{u}}}\|_{L^2(\Omega_F)}+\|\boldsymbol{\bar{\bar{u}}}\|_{L^2(\Omega_F)}\\[5pt]
& \displaystyle =
\|\boldsymbol{u}^0_{j}\|_{L^2(\Omega_F)}+\|\boldsymbol{\bar{\bar{u}}}\|_{L^2(\Omega_F)}
\leq C_1 \|\nabla\boldsymbol{u}^0_{j}\|_{L^2(\Omega_F)}+C_2
\end{array}
\end{equation}
Note that for the $H^1_0(\Omega_F)$-function
$\boldsymbol{u}^0_{j}$ one has
\[
\begin{array}{l l}
 &
\displaystyle
W_{\Omega_F}(\boldsymbol{u}^0_{j})=\mu\int_{\Omega_F}\left\{\left(\frac{\partial
u^0_{j1}}{\partial x}\right)^2+\frac{1}{2}\left(\frac{\partial
u^0_{j1}}{\partial y}+\frac{\partial u^0_{j2}}{\partial
x}\right)^2+\left(\frac{\partial u^0_{j2}}{\partial
y}\right)^2\right\} d\boldsymbol{x}\\[10pt]
 & \displaystyle = \mu\int_{\Omega_F}\left\{\left(\frac{\partial
u^0_{j1}}{\partial x}\right)^2+\frac{1}{2}\left(\frac{\partial
u^0_{j1}}{\partial y}\right)^2+\frac{1}{2}\left(\frac{\partial
u^0_{j2}}{\partial x}\right)^2+\frac{\partial u^0_{j1}}{\partial
y}\frac{\partial u^0_{j2}}{\partial x}+\left(\frac{\partial
u^0_{j2}}{\partial y}\right)^2\right\} d\boldsymbol{x}
\\[10pt]
& \displaystyle = \mu\int_{\Omega_F}\left\{\left(\frac{\partial
u^0_{j1}}{\partial x}\right)^2+\frac{1}{2}\left(\frac{\partial
u^0_{j1}}{\partial y}\right)^2+\frac{1}{2}\left(\frac{\partial
u^0_{j2}}{\partial x}\right)^2+\frac{\partial u^0_{j1}}{\partial
x}\frac{\partial u^0_{j2}}{\partial y}+\left(\frac{\partial
u^0_{j2}}{\partial y}\right)^2\right\} d\boldsymbol{x}\\[10pt]
& \displaystyle =
\mu\int_{\Omega_F}\left\{\frac{1}{2}\left(\frac{\partial
u^0_{j1}}{\partial x}\right)^2+\frac{1}{2}\left(\frac{\partial
u^0_{j1}}{\partial
y}\right)^2+\frac{1}{2}\left(\frac{1}{2}\frac{\partial
u^0_{j2}}{\partial x}\right)^2+\frac{1}{2}\left(\frac{\partial
u^0_{j2}}{\partial y}\right)^2\right\} d\boldsymbol{x}\\[10pt]
& \displaystyle = \frac{\mu}{2}\int_{\Omega_F} |\nabla
\boldsymbol{u}^0_{j} |^2d\boldsymbol{x}= \frac{\mu}{2} \|\nabla
\boldsymbol{u}^0_{j}\|^2_{L^2(\Omega_F)}
\end{array}
\]
Here we used the integration by parts (all boundary terms
disappear since $\boldsymbol{u}^0_{j}\in H^1_0(\Omega_F)$ and the
divergence free property of $\boldsymbol{u}^0_{j}$. Thus,
continuing \eqref{E:step7} we obtain
\[
\begin{array}{l l}
& \displaystyle C_1
\|\nabla\boldsymbol{u}^0_{j}\|_{L^2(\Omega_F)}+C_2\leq C_3
W_{\Omega_F}(\boldsymbol{u}^0_{j})+ C_2\\[10pt]
& \displaystyle \leq
C_3\left(W_{\Omega_F}(\boldsymbol{\bar{\bar{u}}}_{j})+W_{\Omega_F}(\boldsymbol{\bar{\bar{u}}})+
2\sqrt{W_{\Omega_F}(\boldsymbol{\bar{\bar{u}}}_{j})W_{\Omega_F}(\boldsymbol{\bar{\bar{u}}})}\right)+
C_2 \leq C_4
\end{array}
\]
So the sequence $\boldsymbol{\bar{\bar{u}}}_{j}$ is bounded in
$L^2(\Omega_F)$, as well as
\[
\|\nabla\boldsymbol{\bar{\bar{u}}}_{j}\|_{L^2(\Omega_F)}\leq
\|\nabla\boldsymbol{u}^0_{j}\|_{L^2(\Omega_F)}+\|\nabla\boldsymbol{\bar{\bar{u}}}\|_{L^2(\Omega_F)}\leq
C_5.
\]
Note that $\boldsymbol{\bar{\bar{u}}}\in H^1(\Omega_F)$ as a
solution to the Stokes problem \cite{galdi}.

Therefore, from all the above $\boldsymbol{\bar{\bar{u}}}_{j}$ is
bounded sequence in $H^1(\Omega_F)$.

\textit{Step 8.} From the previous step it follows that there
exists a vector field $\boldsymbol{u} \in H^1(\Omega_F)$ such that
$\boldsymbol{\bar{\bar{u}}}_{j}\rightharpoonup \boldsymbol{u}$
(weakly) in $H^1(\Omega_F)$ as $j\rightarrow\infty$.

We have to show that $\boldsymbol{u} \in V$. Recall that for
$\boldsymbol{\bar{\bar{u}}}\in V$:
$\boldsymbol{\bar{\bar{u}}}_j-\boldsymbol{\bar{\bar{u}}}\in
H^1_0(\Omega_F)$, which is closed linear subspace of
$H^1(\Omega_F)$, hence, $H^1_0(\Omega_F)$ is weakly closed
(Mazur's theorem), thus
\[
\boldsymbol{\bar{\bar{u}}}_j-\boldsymbol{\bar{\bar{u}}}\rightharpoonup
\boldsymbol{u}-\boldsymbol{\bar{\bar{u}}},
\]
therefore, $\boldsymbol{u}\in V$.

\textit{Step 9.} Finally, we show that
$W_{\Omega_F}(\boldsymbol{u})=m$. The functional
$W_{\Omega_F}(\cdot)$ is convex in $\nabla\boldsymbol{u}$ hence it
is weakly lower semi-continues, that is,
\[
W_{\Omega_F}(\boldsymbol{u})\leq\liminf_{j\rightarrow\infty}\boldsymbol{\bar{\bar{u}}}_j=m.
\]
Since $\boldsymbol{u}\in V$, we have
\[
W_{\Omega_F}(\boldsymbol{u})=m=\min_{V}W_{\Omega_F}(\cdot).
\]

\end{proof}

\subsection{Existence and uniqueness of the solution to \eqref{E:W-fict}} \label{A:ex-un-2}

\begin{lemma} \label{L:exist-uniq2}
There exists a unique minimizer of $W_{\boldsymbol{\Pi}}(\cdot)$ over $V_{\boldsymbol{\Pi}}$.
\end{lemma}

\begin{proof} To show the existence of a minimizer of functional
$\displaystyle W_{\boldsymbol{\Pi}}(\boldsymbol{v})$ for
$\boldsymbol{v} \in V_{\boldsymbol{\Pi}}$ we have to prove that it is coercive and
convex in the gradient of $\boldsymbol{v}$ and that $V_{\boldsymbol{\Pi}}$ is
nonempty convex and closed set.

\textit{Step 1.} Set $m:=\inf_{V_{\boldsymbol{\Pi}}} W_{\boldsymbol{\Pi}}$. If
$m=+\infty$, we are done. So we hereafter assume that $m$ is
finite. Select a minimizing sequence
$\{\boldsymbol{u}_n\}_{n=1}^\infty \subset V_{\boldsymbol{\Pi}}$ of
$W_{\boldsymbol{\Pi}}(\cdot)$, that is,
$W_{\boldsymbol{\Pi}}(\boldsymbol{u}_n)\rightarrow m$ as
$n\rightarrow\infty$.

\textit{Step 2.} By the argument similar to Step 2 of Lemma
\ref{L:exist-uniq1} one can choose a bounded subsequence of the
velocities (which are constant vectors) on the disks $B^i$,
$i=1,\ldots,\,N$. Denote such a subsequence by
$\{\boldsymbol{v}_n^i\}_{n=1}^{\infty}$, $i=1,\ldots,\,N$, that
is, there exists a positive $C_0$ such that for all $n$ and every
$i=1,\ldots,\,N$: $|\boldsymbol{v}_n^i|<C_0$. Denote the limits of
those subsequences by $\boldsymbol{v}^i$.

Note that by the Korn's inequality (see e.g. \cite{brenner}) we
have:
\begin{equation} \label{E:korn-cont}
\|\nabla \boldsymbol{u}_n\|^2_{L^2(\boldsymbol{\Pi})}\leq
C_{\boldsymbol{\Pi}}\left(W_{\boldsymbol{\Pi}}(\boldsymbol{u}_n)+\|\boldsymbol{u}_n\|^2_{L^2(\boldsymbol{\Pi})}\right),
\end{equation}
where $C_{\boldsymbol{\Pi}}$ is some constant. Therefore, we first
show that the sequence $\boldsymbol{u}_n$ is bounded in
$L^2(\boldsymbol{\Pi})$.

\textit{Step 3.} Select a sequence $\boldsymbol{w}_n$ of solutions
to the Stokes problem \eqref{E:Dir-pr} with
$\boldsymbol{w}_n=\boldsymbol{v}_n^i$ on each $\partial B^i$, and
the functions $\boldsymbol{w}$ that solves \eqref{E:Dir-pr} with
$\boldsymbol{w}=\boldsymbol{v}^i$ on $\partial B^i$,
$i=1,\ldots,N$. Then,
\[
\|\boldsymbol{w}_n\|_{H^1(\Omega_F)}\leq
\|\boldsymbol{w}_n-\boldsymbol{w}\|_{H^1(\Omega_F)}+\|\boldsymbol{w}\|_{H^1(\Omega_F)},
\]
where both norms of the right-hand side are bounded \cite{galdi}
as solutions to Dirichlet problem, thus,
$\|\boldsymbol{w}_n\|_{H^1(\boldsymbol{\Pi})}\leq \mathcal{C}$.

Define a sequence
$\boldsymbol{\varphi}_n=\boldsymbol{w}_n|_{\boldsymbol{\Pi}}-\boldsymbol{u}_n$
in $\boldsymbol{\Pi}$, which is zero vectors on each sphere
$\partial B^i$ and zero vector on the external boundary $\partial
\Omega$. Using \eqref{E:korn-cont} for the sequence
$\boldsymbol{\varphi}_n$ and the fact that
\[
\begin{array}{l l}
W_{\boldsymbol{\Pi}}(\boldsymbol{\varphi}_n) & \displaystyle \leq
W_{\boldsymbol{\Pi}}(\boldsymbol{w}_n)+W_{\boldsymbol{\Pi}}(\boldsymbol{u}_n)+
2\sqrt{W_{\boldsymbol{\Pi}}(\boldsymbol{w}_n)W_{\boldsymbol{\Pi}}(\boldsymbol{u}_n)}\\[5pt]
 & \displaystyle \leq
\mathcal{C}+(m+\varepsilon_n)+2\sqrt{\mathcal{C}(m+\varepsilon_n)}\leq
C_1,
\end{array}
\]
we prove that $\boldsymbol{\varphi}_n$ is bounded in
$L^2(\boldsymbol{\Pi})$.

Consider the neck $\Pi_{ij}$ (Fig. \ref{F:neck1}) between adjacent
disks $B^i$ and $B^j$. Choose a point $(x_0,y_0)$ on the sphere
$\partial B^j$ and $(x,y)\in \Pi_{ij}$, and consider
\[
\varphi_{n}^{1}(x,y)=\int_{x_0}^{x}\frac{\partial
\varphi_{n}^{1}}{\partial \xi}(\xi,y)d\xi,
\]
for almost all $x$, where the superscript here indicates a
component of the the vector
$\boldsymbol{\varphi}_n=(\boldsymbol{\varphi}_n^1,\boldsymbol{\varphi}_n^2)$.
Then we obtain the following estimate:
\[
|\varphi_{n}^{1}(x,y)|^2\leq C_2
\int_{x_0}^{H(\gamma)/2}\left|\frac{\partial
\varphi_{n}^{1}}{\partial \xi}(\xi,y)\right|^2 d\xi,
\]
and $\gamma=\max(\gamma_{ij}^-,\gamma_{ij}^+)$. Now we integrate
the last inequality with respect to $y$ between $\gamma_{ij}^-$
and $\gamma_{ij}^+$, and with respect to $x$ along the arc $a_i$
obtaining
\[
\int_{\Pi_{ij}}|\varphi_{n}^{1}(x,y)|^2
d\boldsymbol{x}=\|\varphi_{n}^{1}(x,y)\|^2_{L^2(\Pi_{ij})}\leq
C_2|\gamma_{ij}^+-\gamma_{ij}^-|H(\gamma)/2 = C_3.
\]
Here we recall that $\int_{\Pi_{ij}}|\frac{\partial
\varphi_{n}^{1}}{\partial \xi}|^2d\boldsymbol{x}$ is bounded from
above by $W_{\boldsymbol{\Pi}}(\boldsymbol{\varphi}_n)$. So the
first components of $\boldsymbol{\varphi}_n$ form a bounded
sequence in $L^2(\Pi_{ij})$ for any neck $\Pi_{ij}\in
\boldsymbol{\Pi}$.

To show that the second components of $\boldsymbol{\varphi}_n$ is
a bounded sequence we consider another neck $\Pi_{ki}$ (see Fig.
\ref{F:neck1}(b)) that forms a small angle $\alpha$ with the
direction of the second component $\varphi_{n}^{2}$. Then by
triangle inequality we have
\[
\|\varphi_{n}^{2}\|_{L^{2}(\Pi_{ik})}\leq
C(\alpha)\left(\|\varphi_{n}^{1}\|_{L^{2}(\Pi_{ik})}+\|\varphi_{n}^{3}\|_{L^{2}(\Pi_{ik})}\right),
\]
where the constant $C(\alpha)$ depends on the small angle
$\alpha$. Since the viscous dissipation rate is invariant under
rotation we evaluate the norm of $\varphi_{n}^{3}$ to obtain an
estimate for $\varphi_{n}^{2}$.

In the system of coordinate with $y$-axis along the neck
$\Pi_{ki}$ we consider
\[
\varphi_{n}^{3}(x,y)=\int_{y_0}^{y}\frac{\partial
\varphi_{n}^{3}}{\partial \eta}(x,\eta)d\eta,
\]
where the point $(x_0,y_0)$ lies on $\partial B^k$. As before,
integrating with respect to $x$ over
$(\gamma_{ki}^-,\gamma_{ki}^+)$ and with respect to $y$ over the
arc $a_k$ and using that $\int_{\Pi_{ij}}|\frac{\partial
\varphi_{n}^{3}}{\partial \eta}|^2d\boldsymbol{x}$ is bounded from
above by $W_{\boldsymbol{\Pi}}(\boldsymbol{\varphi}_n)$, we obtain
the following estimate:
\[
\int_{\Pi_{ij}}|\varphi_{n}^{3}(x,y)|^2
d\boldsymbol{x}=\|\varphi_{n}^{3}(x,y)\|^2_{L^2(\Pi_{ij})}\leq
C_4.
\]

Hence, the sequence $\boldsymbol{\varphi}_n$ is bounded in
$L^2(\boldsymbol{\Pi})$. In the view of the Korn's inequality
\eqref{E:korn-cont} we obtain that $\boldsymbol{\varphi}_n$ is bounded
in $H^1(\boldsymbol{\Pi})$. By triangle inequality we obtain that
sequence $\boldsymbol{u}_n$ is bounded in $H^1(\boldsymbol{\Pi})$
as well.

The boundedness of $\boldsymbol{u}_n$ in $H^1(\boldsymbol{\Pi})$
implies there exists a subsequence $\boldsymbol{u}_{nn'}$ and a
function $\boldsymbol{\hat{u}}\in H^1(\Pi_{ij})$ such that
$\boldsymbol{u}_{nn'}\rightharpoonup \boldsymbol{\hat{u}}$
(weakly) in $H^1(\Pi_{ij})$ as $n'\rightarrow\infty$. It remains
to show that $\boldsymbol{\hat{u}} \in V_{\boldsymbol{\Pi}}$ and
$W_{\boldsymbol{\Pi}}(\boldsymbol{\hat{u}})=m$. To do so we need
to prove that $V_{\boldsymbol{\Pi}}$ is convex and closed set of
$H^1(\boldsymbol{\Pi})$, and $W_{\boldsymbol{\Pi}}(\cdot)$ is a
convex functional.

\begin{figure}[!ht]
  \centering
  \includegraphics[scale=0.95]{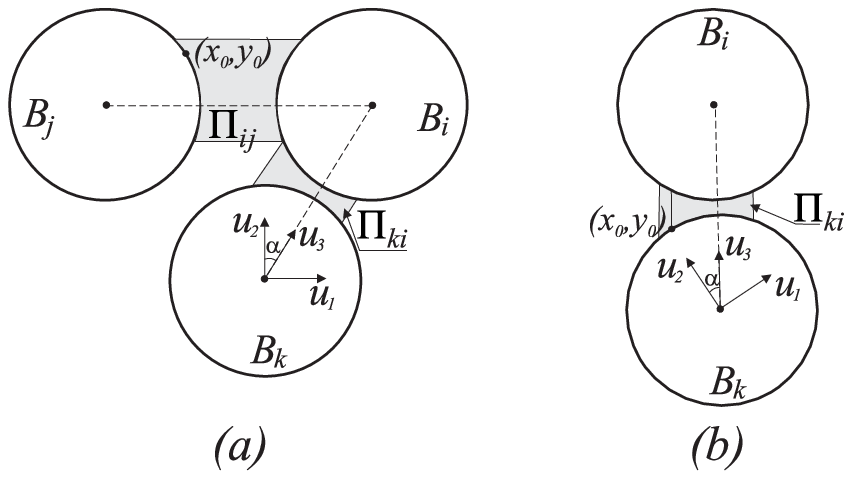}
  \caption{}\label{F:neck1}
\end{figure}

\textit{Step 4.} $V_{\boldsymbol{\Pi}}$ is convex if for any
$\boldsymbol{u},\,\boldsymbol{v} \in V_{\boldsymbol{\Pi}}$ their convex combination
is from  $V_{ij}$:
$\boldsymbol{w}=(1-\lambda)\boldsymbol{u}+\lambda\boldsymbol{v}
\in V_{\boldsymbol{\Pi}}$, $0<\lambda <1$. To show that $\boldsymbol{w} \in V_{\boldsymbol{\Pi}}$ we
have to prove $\nabla\cdot\boldsymbol{w}=0$ in $\boldsymbol{\Pi}$,
$\boldsymbol{w}=\boldsymbol{f}$ on $\partial\boldsymbol{\Pi}\cap
\partial\Omega$,
$\boldsymbol{w}=\boldsymbol{U}^{i}+R\omega^i(n_1^i\boldsymbol{e}_2-n_2^i\boldsymbol{e}_1)$
on $\partial B^{i}$ and $\displaystyle
\int_{\partial\triangle_{ijk}}\boldsymbol{w}\cdot\boldsymbol{n}ds=0$,
which all clearly hold.

\textit{Step 5.} To prove that $V_{\boldsymbol{\Pi}}$ is closed we consider a
sequence $\{\boldsymbol{v}_{n}\}_{n=1}^{\infty}$,
$\boldsymbol{v}_n \in V_{\boldsymbol{\Pi}}$ for all $n$ and
$\boldsymbol{v}_n\rightarrow \boldsymbol{v}$ in
$H^1(\boldsymbol{\Pi})$, and show that $\boldsymbol{v} \in V_{\boldsymbol{\Pi}}$ as
well. For this we need to check that $\nabla\cdot\boldsymbol{v}=0$
in $\boldsymbol{\Pi}$ in the weak sense, and
$\boldsymbol{v}=\boldsymbol{U}^{i}+R\omega^i(n_1^i\boldsymbol{e}_2-n_2^i\boldsymbol{e}_1)$
on $\partial B^{i}$, $i=1,\ldots,N$, $\displaystyle
\int_{\partial\triangle_{ijk}}\boldsymbol{v}\cdot\boldsymbol{n}ds=0$,
and $\boldsymbol{v}=\boldsymbol{f}$ if
$\partial\boldsymbol{\Pi}\cap
\partial\Omega$, which understood in the sense of trace. Any
convergent sequence in $H^1(\boldsymbol{\Pi})$ is weakly
convergent, hence, $\displaystyle
\int_{\boldsymbol{\Pi}}\boldsymbol{v}_n\cdot\nabla \varphi
d\boldsymbol{x}\rightarrow
\int_{\boldsymbol{\Pi}}\boldsymbol{v}\cdot\nabla \varphi
d\boldsymbol{x}$ for every $\varphi \in
C_0^\infty(\boldsymbol{\Pi})$. Thus the incompressibility of
$\boldsymbol{v}$ holds.

The boundary conditions of the definition of $V_{\boldsymbol{\Pi}}$ are considered
in the sense of trace, that is, for the $H^1$ function
$\boldsymbol{v}_n$ there exist a linear bounded map
$\gamma(\boldsymbol{v}_n)$ and a constant $C$ such that
$\|\gamma(\boldsymbol{v}_n)\|_{L_{2}(\partial
\boldsymbol{\Pi})}\leq
C\|\boldsymbol{v}_n\|_{H^1(\boldsymbol{\Pi})}$. Therefore, we have
\[
\begin{array}{l l}
& \displaystyle \sum_{i\in \mathbb{I},j,k\in
\mathcal{N}_i}\int_{\partial\triangle_{ijk}}\left(\gamma(\boldsymbol{v}_k)-\gamma(\boldsymbol{v})\right)\cdot\boldsymbol{n}\,ds=
\sum_{i\in \mathbb{I},j,k\in
\mathcal{N}_i}\int_{\partial\triangle_{ijk}}\gamma(\boldsymbol{v}_n-\boldsymbol{v})\cdot\boldsymbol{n}\,ds\\[13pt]
\leq & \displaystyle \sum_{i\in \mathbb{I},j,k\in
\mathcal{N}_i}\left|\partial\triangle_{ijk}\right|^{1/2}\|\gamma(\boldsymbol{v}_n-\boldsymbol{v})\|^{1/2}_{L_{2}(\partial
\triangle_{ijk})}\leq
C_1\|\gamma(\boldsymbol{v}_n-\boldsymbol{v})\|^{1/2}_{L_{2}(\partial
\boldsymbol{\Pi})}
\\[13pt]
\leq & \displaystyle
C_2\|\boldsymbol{v}_n-\boldsymbol{v}\|^{1/2}_{H^1(\boldsymbol{\Pi})}\rightarrow
0 \,\, \mbox{as}\,\, n\rightarrow\infty,
\end{array}
\]
where $\displaystyle C_2=\left(C\sum_{i\in \mathbb{I},j,k\in
\mathcal{N}_i}\left|\partial\triangle_{ijk}\right|\right)^{1/2}$,
and $\displaystyle
\partial\boldsymbol{\Pi}=\bigcup_{i\in\mathbb{I}}\partial
B^i\cup \bigcup_{i\in\mathbb{I},i,k\in\mathcal{N}_i}\partial
\triangle_{ijk}$ The remaining boundary conditions can be proved
similarly.

Finally, we show that $W_{\boldsymbol{\Pi}}(\cdot)$ is convex in
$\nabla \boldsymbol{u}$. For this we write
\[
W_{\boldsymbol{\Pi}}(\nabla\boldsymbol{u})=\mu\int_{\boldsymbol{\Pi}}L(\nabla\boldsymbol{u})d\boldsymbol{x},
\]
where
$L(\boldsymbol{p})=p_{11}^2+\frac{1}{2}(p_{12}+p_{21})^2+p_{22}^2$
and $\boldsymbol{p}=(p_{ij})_{i,j=1,2}$ is $2\times 2$ matrix.
Recall that by definition convexity means that
\[
\sum_{i,j,k,l=1}^{2}L_{p_{ij}p_{kl}}(p)a_{ij}a_{kl}\geq 0, \quad
\mbox{ for any } A=(a_{ij})_{i,j=1,2}.
\]

Therefore, we have
\[
\begin{array}{l l l}
& \displaystyle L_{p_{11}p_{11}}=L_{p_{22}p_{22}}=2, \quad L_{p_{12}p_{12}}=L_{p_{21}p_{21}}=L_{p_{12}p_{21}}=L_{p_{21}p_{12}}=1,\\[9pt]
& \displaystyle L_{p_{11}p_{12}}=L_{p_{11}p_{21}}=L_{p_{11}p_{22}}=L_{p_{12}p_{11}}=L_{p_{12}p_{22}}=L_{p_{21}p_{11}}=L_{p_{21}p_{22}}=0,\\[9pt]
& \displaystyle L_{p_{22}p_{11}}=L_{p_{22}p_{12}}=L_{p_{22}p_{21}}=0.
\end{array}
\]
Hence,
\[
2a_{11}^2+2a_{22}^2+ a_{12}^2+a_{21}^2 +a_{12}a_{12}+a_{21}a_{21}=2a_{11}^2+2a_{22}^2+ (a_{12}+a_{21})^2\geq 0,
\]
thus, $W_{\boldsymbol{\Pi}}(\cdot)$ is convex in $\nabla\boldsymbol{u}$. This implies that $W_{\boldsymbol{\Pi}}$ is convex, it is weakly lower semicontinuous:
\[
W_{\boldsymbol{\Pi}}(\boldsymbol{\hat{u}})\leq \liminf_n W_{\boldsymbol{\Pi}}(\boldsymbol{u}_n)=m,
\]
and $W_{\boldsymbol{\Pi}}(\boldsymbol{\hat{u}})=m$.

By Mazur's Theorem $V_{\boldsymbol{\Pi}}$ is weakly closed, hence, $\boldsymbol{\hat{u}} \in V_{\boldsymbol{\Pi}}$.

\textit{Step 6.} To show a uniqueness of the minimizer we proceed by contradiction. Assume that there exist two minimizers of this
problem $\boldsymbol{\hat{u}}$ and $\boldsymbol{\hat{u}}$. Consider their difference $\boldsymbol{\hat{w}}=\boldsymbol{\hat{u}}-\boldsymbol{\hat{v}}$
that solves \eqref{E:EL-boldneck} with $\boldsymbol{f}=\boldsymbol{0}$ on $\partial\Omega$ and some
$\boldsymbol{\hat{W}}^i+R\omega_{\hat{w}}^{i}(n_1^i\boldsymbol{e}_2-n_2^i\boldsymbol{e}_1)$ on $\partial B^i$, $i=1,\ldots,N$. By multiplying equation
(\ref{E:EL-boldneck}a) by $\boldsymbol{\hat{w}}$ and integrating by parts over $\boldsymbol{\Pi}$ one obtains that $W_{\boldsymbol{\Pi}}(\boldsymbol{\hat{w}})=0$.

Using the argument given by Step 2 of Lemma \ref{L:exist-uniq1} we
can show that $\boldsymbol{\hat{w}}=\boldsymbol{0}$ on each disk
$B^i$, $i=1,\ldots,N$ and the argument of Step 3 of the current
lemma we can prove that $\boldsymbol{\hat{w}}=\boldsymbol{0}$ in
$\boldsymbol{\Pi}$.

\end{proof}

\subsection{Existence and uniqueness of the solution to \eqref{E:EL-lower-beta}} \label{A:ex-un-3}
Suppose there are two solutions to \eqref{E:EL-lower-beta}. Then their difference satisfies
\begin{equation}   \label{A:EL-lower-beta}
\left\{
\begin{array}{r l l}
(a) & \displaystyle \mu\triangle \boldsymbol{u}=\nabla p, & \boldsymbol{x} \in \Pi_{ij}, \\
(b) & \displaystyle \nabla\cdot \boldsymbol{u}=0, & \boldsymbol{x} \in \Pi_{ij},\\
(c') & \displaystyle \boldsymbol{u}=0, & \boldsymbol{x} \in \partial B^{i}, \\
(c'') & \displaystyle \boldsymbol{u}=0, & \boldsymbol{x} \in \partial B^{j}, \\
(d) & \displaystyle \frac{1}{R}\int_{\ell_{ij}}\boldsymbol{u}\cdot\boldsymbol{n}ds=0,\\
(e) & \displaystyle \boldsymbol{\sigma}(\boldsymbol{u})\boldsymbol{n}=-p^{\pm}_{ij}\boldsymbol{n}, & \boldsymbol{x} \in \partial\Pi^{\pm}_{ij},\\
(f) & \displaystyle \boldsymbol{u}=0,& \boldsymbol{x} \in \partial\Pi_{ij}\cap \partial\Omega.
\end{array}
\right.
\end{equation}
Integrating by parts the incompressibility condition   \eqref{A:EL-lower-beta}(b)  over  half-necks
and using  \eqref{A:EL-lower-beta}(b),(c'),(c''),(d),(e), we obtain:
\[
\int_{ \partial\Pi^{+}_{ij}}  \boldsymbol{u} \cdot \boldsymbol{n} ds =
\int_{ \partial\Pi^{-}_{ij}}  \boldsymbol{u} \cdot \boldsymbol{n} ds=0.
\]
Multiplying  \eqref{A:EL-lower-beta}(a) by $ \boldsymbol{u}$ and integrating by parts,
we obtain
\[
W_{\Pi_{ij}}( \boldsymbol{u}) = \int_{\partial \Pi_{ij}}
\boldsymbol{u} \cdot \boldsymbol{\sigma}(\boldsymbol{u})\boldsymbol{n} ds =
p^+ \int_{ \partial\Pi^{+}_{ij}} \boldsymbol{u} \cdot \boldsymbol{n} ds +
p^- \int_{ \partial\Pi^{-}_{ij}} \boldsymbol{u} \cdot \boldsymbol{n} ds=0.
\]
Hence $ \boldsymbol{u}$ must represent a rigid-body rotation. Due to our homogeneous boundary
conditions, this motion must be identically zero.

Note, that since $\boldsymbol{u}$ the solution to \eqref{A:EL-lower-beta} must be identically zero, it implies that $p^{\pm}_{ij}=0$. Hence we showed uniqueness of $p^{\pm}_{ij}$ as well.

\subsection{Dual problem to \eqref{E:W-fict}} \label{A:dual-fnl}

\begin{lemma} \label{L:maximizer}
Let $\boldsymbol{u}_{per}=(u_1^{per},u_2^{per})$ is the minimizer of $W_{\Pi_{ij}}(\cdot)$ over $V_{per}$. Then the stress tensor
$\boldsymbol{\sigma}(\boldsymbol{u}_{per})$ is the maximizer of $W^*_3(\cdot)$, defined by \eqref{E:dual-fnl-2}, over the set $F_{per}$, defined by \eqref{E:dual-set-2}.
\end{lemma}

\textit{Proof.} For simplicity we use a notation $\boldsymbol{\sigma}$ for $\boldsymbol{\sigma}(\boldsymbol{u}_{per})$.
To prove that $\boldsymbol{\sigma}$ is a maximizer of $W^*_3(\cdot)$ we consider an arbitrary tensor $\boldsymbol{T} \in F_{per}$, such that
$\boldsymbol{T}\boldsymbol{n}=\chi_T^\pm\boldsymbol{n}$ on $\partial \Pi_{ij}^\pm$, where $\chi_T^+$, $\chi_T^-$ are arbitrary numbers,
then evaluate $W^*_3(\boldsymbol{\sigma}+\boldsymbol{T})$ and show that
$W^*_{3}(\boldsymbol{\sigma}) \geq W^*_{3}(\boldsymbol{\sigma}+\boldsymbol{T})$.

Denote the first linear part of \eqref{E:dual-fnl-2} by
\[
\begin{array}{l l}
\mathcal{L}(\boldsymbol{\sigma})& \displaystyle
=\frac{R\beta_{ij}^*}{H(\gamma_{ij}^+)}\int_{\partial\Pi_{ij}^+}\boldsymbol{n}\cdot\boldsymbol{\sigma}\boldsymbol{n}ds-
\frac{R\beta_{ij}^*}{H(\gamma_{ij}^-)}\int_{\partial\Pi_{ij}^-}\boldsymbol{n}\cdot\boldsymbol{\sigma}\boldsymbol{n}ds\\[9pt]
& \displaystyle =
\frac{R\beta_{ij}^*}{H(\gamma_{ij}^+)}\int_{\partial\Pi_{ij}^+}\boldsymbol{n}\cdot(p_{ij}^+)\boldsymbol{n}ds-
\frac{R\beta_{ij}^*}{H(\gamma_{ij}^-)}\int_{\partial\Pi_{ij}^-}\boldsymbol{n}\cdotp_{ij}^-\boldsymbol{n}ds \\[9pt]
& \displaystyle =
R\beta_{ij}^*\left(\frac{1}{H(\gamma_{ij}^+)}p_{ij}^+|\partial\Pi_{ij}^+|- \frac{1}{H(\gamma_{ij}^-)}p_{ij}^-|\partial\Pi_{ij}^-| \right)\\[9pt]
& \displaystyle =R\beta_{ij}^*(p_{ij}^+-p_{ij}^-),
\end{array}
\]
while
\[
\begin{array}{l l}
\mathcal{L}(\boldsymbol{T})& \displaystyle
=\frac{R\beta_{ij}^*}{H(\gamma_{ij}^+)}\int_{\partial\Pi_{ij}^+}\boldsymbol{n}\cdot\boldsymbol{T}\boldsymbol{n}ds-
\frac{R\beta_{ij}^*}{H(\gamma_{ij}^-)}\int_{\partial\Pi_{ij}^-}\boldsymbol{n}\cdot\boldsymbol{T}\boldsymbol{n}ds\\[7pt]
& \displaystyle =R\beta_{ij}^*(\chi_T^+-\chi_T^-).
\end{array}
\]
Also denote
\[
\mathcal{M}(\boldsymbol{\sigma})=R(\omega^i-\omega^j)\left[\int_{\partial B^i}\boldsymbol{G}^i_{3}\cdot\boldsymbol{\sigma}\boldsymbol{n}ds
+\int_{\partial B^j}\boldsymbol{G}^j_{3}\cdot\boldsymbol{\sigma}\boldsymbol{n}ds\right].
\]

Now we consider the quadratic term of \eqref{E:dual-fnl-2} and notice that for $\boldsymbol{\sigma}+\boldsymbol{T}$ one has
\[
\begin{array}{l l}
& \displaystyle
-\frac{1}{4\mu}\int_{\Pi_{ij}}\left[(\boldsymbol{\sigma}+\boldsymbol{T}):(\boldsymbol{\sigma}+\boldsymbol{T})
-\frac{\mbox{tr}(\boldsymbol{\sigma}+\boldsymbol{T})^2}{2}\right]d\boldsymbol{x}\\[9pt]
& \displaystyle =
-\frac{1}{4\mu}\int_{\Pi_{ij}}\left[\boldsymbol{\sigma}:\boldsymbol{\sigma}+2\boldsymbol{\sigma}:\boldsymbol{T}+
\boldsymbol{T}:\boldsymbol{T}
-\frac{(\mbox{tr}\boldsymbol{\sigma})^2}{2}-
\mbox{tr}\,\boldsymbol{\sigma}\,\mbox{tr}\,\boldsymbol{T}
-\frac{(\mbox{tr}\,\boldsymbol{T})^2}{2}\right]d\boldsymbol{x}\\[9pt]
& \displaystyle =
\boldsymbol{w}(\boldsymbol{\sigma})+2\textit{w}(\boldsymbol{\sigma},\boldsymbol{T})+\boldsymbol{w}(\boldsymbol{T}),
\end{array}
\]
where
\begin{equation} \label{E:quadr-form-last}
\textit{w}(\boldsymbol{\sigma},\boldsymbol{T})=
-\frac{1}{4\mu}\int_{\Pi_{ij}}\left[(\boldsymbol{\sigma}:\boldsymbol{T})-\frac{\mbox{tr}\,\boldsymbol{\sigma}\,
\mbox{tr}\,\boldsymbol{T}}{2}\right]d\boldsymbol{x},
\end{equation}
and the corresponding quadratic form:
\[
\boldsymbol{w}(\boldsymbol{\sigma})=
-\frac{1}{4\mu}\int_{\Pi_{ij}}\left[(\boldsymbol{\sigma}:\boldsymbol{\sigma})-\frac{(\mbox{tr}\,\boldsymbol{\sigma})^2
}{2}\right]d\boldsymbol{x}.
\]
Hence, $W^*_{3}(\boldsymbol{\sigma}+\boldsymbol{T})$ can be
written as:
\[
W^*_{3}(\boldsymbol{\sigma}+\boldsymbol{T})= \mathcal{L}(\boldsymbol{\sigma})+\mathcal{L}(\boldsymbol{T})+\mathcal{M}(\boldsymbol{\sigma})+\mathcal{M}(\boldsymbol{T})
+\boldsymbol{w}(\boldsymbol{\sigma})+2\textit{w}(\boldsymbol{\sigma},\boldsymbol{T})+\boldsymbol{w}(\boldsymbol{T}),
\]
where
$W^*_{3}(\boldsymbol{\sigma})=\mathcal{L}(\boldsymbol{\sigma})+\mathcal{M}(\boldsymbol{\sigma})+\boldsymbol{w}(\boldsymbol{\sigma})$
in these notations.

Taking into account that
\begin{equation*}
\boldsymbol{\sigma}=\boldsymbol{\sigma}(\boldsymbol{u}_{per})=
\begin{pmatrix} 2\mu\frac{\partial u_1^{per}}{\partial x}-p_\sigma & \mu\left(\frac{\partial u_1^{per}}{\partial y}+\frac{\partial u_2^{per}}{\partial x}\right)\\
\mu\left(\frac{\partial u_1^{per}}{\partial y}+\frac{\partial u_2^{per}}{\partial x}\right) & 2\mu\frac{\partial u_2^{per}}{\partial y}-p_\sigma \end{pmatrix},
\end{equation*}
and
\[
p_\sigma=-\frac{\mbox{tr}\,\boldsymbol{\sigma}}{2},
\]
we consider
\begin{equation} \label{E:quad-form-1}
\begin{array}{l l}
  \displaystyle 2\textit{w}(\boldsymbol{\sigma},\boldsymbol{T}) &  \displaystyle =-\frac{1}{2\mu}\int_{\Pi_{ij}}\left[(\boldsymbol{\sigma}:\boldsymbol{T})-\frac{\mbox{tr}\,\boldsymbol{\sigma}\, \mbox{tr}\,\boldsymbol{T}}{2}\right]d\boldsymbol{x}\\[9pt]
& \displaystyle = -\frac{1}{2\mu}\int_{\Pi_{ij}}\left[\left(2\mu\frac{\partial u_1^{per}}{\partial x}-p_\sigma\right)T_{11}+
2\mu\left(\frac{\partial u_1^{per}}{\partial y}+\frac{\partial u_2^{per}}{\partial x}\right)T_{12}\right.\\[9pt]
&  \displaystyle \left. + \left(2\mu\frac{\partial u_2^{per}}{\partial y}-p_\sigma\right)T_{22}+p_\sigma(T_{11}+T_{22}) \right]d\boldsymbol{x}\\[9pt]
&  \displaystyle = -\frac{1}{2\mu}\int_{\Pi_{ij}}\left[2\mu\frac{\partial u_1^{per}}{\partial x}T_{11}+ 2\mu\frac{\partial u_1^{per}}{\partial y}T_{12}+
2\mu\frac{\partial u_2^{per}}{\partial x} T_{12}+ 2\mu\frac{\partial u_2^{per}}{\partial y}T_{22} \right]d\boldsymbol{x}  \\[9pt]
&  \displaystyle = -\int_{\Pi_{ij}}\left[\frac{\partial u_1^{per}}{\partial x}T_{11}+ \frac{\partial u_1^{per}}{\partial y}T_{12}+
\frac{\partial u_2^{per}}{\partial x} T_{12}+ \frac{\partial u_2^{per}}{\partial y}T_{22} \right]d\boldsymbol{x}  \\[9pt]
&  \displaystyle = \int_{\Pi_{ij}}\left[u_1^{per}\frac{\partial T_{11}}{\partial x}+ u_1^{per}\frac{\partial T_{12}}{\partial y}+
u_2^{per}\frac{\partial T_{12}}{\partial x} + u_2^{per}\frac{\partial T_{22}}{\partial y} \right]d\boldsymbol{x} \\[9pt]
&  \displaystyle \quad - \int_{\partial\Pi_{ij}}(u_1^{per} n_1 T_{11}+u_1^{per} n_2 T_{12}+ u_2^{per} n_1 T_{12} + u_2^{per} n_2 T_{22})ds \\[9pt]
&  \displaystyle =  - \int_{\partial\Pi_{ij}^+} \boldsymbol{u}_{per}\cdot \boldsymbol{T}\boldsymbol{n}\,ds-
\int_{\partial\Pi_{ij}^-} \boldsymbol{u}_{per}\cdot \boldsymbol{T}\boldsymbol{n}\,ds\\[9pt]
&  \displaystyle -R(\omega^i-\omega^j)\left(\int_{\partial B^i}\boldsymbol{G}^i_{3}\cdot\boldsymbol{T}\boldsymbol{n}ds+
\int_{\partial B^j}\boldsymbol{G}^j_{3}\cdot\boldsymbol{T}\boldsymbol{n}ds\right),
\end{array}
\end{equation}
where the volume integral over $\Pi_{ij}$ disappears due to
divergence-free property of the tensor $\boldsymbol{T}$.
Continuing \eqref{E:quad-form-1} and taking into account
$\boldsymbol{T}\boldsymbol{n}=\chi_T^\pm\boldsymbol{n}$ on
$\partial\Pi_{ij}^\pm$ and the fact that $\boldsymbol{u}_{per}$ satisfies
$\displaystyle \frac{1}{R}\int_{\partial\Pi_{ij}^\pm}\boldsymbol{u}_{per}\cdot\boldsymbol{n}ds=\pm\beta_{ij}^*$
(this is from the divergence free condition of $\boldsymbol{u}_{per}$ in $\Pi_{ij}^\pm$) we obtain:
\[
\begin{array}{l l}
2\textit{w}(\boldsymbol{\sigma},\boldsymbol{T})
& \displaystyle =-\chi_T^+ \int_{\partial\Pi_{ij}^+} \boldsymbol{u}_{per}\cdot \boldsymbol{n}\,ds-\chi_T^- \int_{\partial\Pi_{ij}^-} \boldsymbol{u}_{per}\cdot \boldsymbol{n}\,ds\\[7pt]
& \displaystyle -R(\omega^i-\omega^j)\left(\int_{\partial B^i}\boldsymbol{G}^i_{3}\cdot\boldsymbol{T}\boldsymbol{n}ds+\int_{\partial B^j}\boldsymbol{G}^j_{3}\cdot\boldsymbol{T}\boldsymbol{n}ds\right)\\[10pt]
& \displaystyle =-\chi_T^+R\beta_{ij}^*+\chi_T^-R\beta_{ij}^*- \mathcal{M}(\boldsymbol{T})
=-\mathcal{L}(\boldsymbol{T})-\mathcal{M}(\boldsymbol{T}).
\end{array}
\]

Finally, we observe that
\[
\boldsymbol{w}(\boldsymbol{T})=
-\frac{1}{4\mu}\int_{\Pi_{ij}}\left[(\boldsymbol{T}:\boldsymbol{T})-\frac{(\mbox{tr}\,\boldsymbol{T})^2
}{2}\right]d\boldsymbol{x}=-\frac{1}{4\mu}\int_{\Pi_{ij}}\left[\frac{1}{2}(T_{11}-T_{22})^2+2T_{12}^2
\right]d\boldsymbol{x}\leq 0.
\]

Therefore, from all the above we conclude
\[
\begin{array}{l l}
W^*_{3}(\boldsymbol{\sigma}+\boldsymbol{T})&  \displaystyle =
\mathcal{L}(\boldsymbol{\sigma})+\mathcal{L}(\boldsymbol{T})+
\boldsymbol{w}(\boldsymbol{\sigma})+2\textit{w}(\boldsymbol{\sigma},\boldsymbol{T})+\boldsymbol{w}(\boldsymbol{T})\\[9pt]
&  \displaystyle =
W^*_{3}(\boldsymbol{\sigma})+\mathcal{L}(\boldsymbol{T})-\mathcal{M}(\boldsymbol{T})+\mathcal{M}(\boldsymbol{T})-\mathcal{L}(\boldsymbol{T})+\boldsymbol{w}(\boldsymbol{T})\\[9pt]
&  \displaystyle = W^*_{3}(\boldsymbol{\sigma})+\boldsymbol{w}(\boldsymbol{T}) \leq W^*_{3}(\boldsymbol{\sigma}),
\end{array}
\]
and this inequality holds for any $\boldsymbol{T} \in F_{per}$. Thus, $\boldsymbol{\sigma}(\boldsymbol{u}_{per})$ is a maximizer of the functional \eqref{E:dual-fnl-2}.

\end{document}